\RequirePackage{fix-cm}
\documentclass{svjour3}                     
\smartqed  
\usepackage{fancyhdr}
\usepackage{booktabs}
\usepackage{graphicx}
\usepackage{amsfonts}
\usepackage{color}
\usepackage{commath}
\usepackage{amsmath}              
\usepackage{latexsym}
\usepackage{subfigure}
\usepackage{url}
\usepackage{natbib}
\setcitestyle{authoryear,open={(},close={)}} 
\begin{document}

\title{Fractional Order Optimal Control Model For Malaria Infection 
}


\author{E. Okyere, F. T. Oduro, S. K. Amponsah, I. K. Dontwi
}


\institute{Eric Okyere \at
              Department of Basic Sciences, University of Health and Allied Sciences, PMB 31, Ho, Volta Region, Ghana
\\
              Tel.: +233-508389028\\
              \email{eokyere@uhas.edu.gh}           
          \and
           F. T. Oduro, S. K. Amponsah, I. K. Dontwi\at
           Department of Mathematics, Kwame Nkrumah University of Science and Technology, Kumasi, Ghana
}

\date{July, 2016}

\maketitle

\begin{abstract}
We propose and study an optimal control model for malaria infection described by system of fractional differential equations. The model is formulated in terms of the left and right Caputo fractional derivatives. We determine the necessary conditions for the optimality of the controlled dynamical system.
The forward-backward sweep method with generalized euler scheme is applied to numerically compute the solutions of the optimality system.
\end{abstract}\\


\section{Introduction}
Mathematical modeling of malaria infection using differential equations with integer order and non-integer derivatives have been considered by many authors
\citep[see, e.g,][]{wow17, Okyere2016, Ngonghala2012, Ross1911, Abdullahi2013, Keegan2013, Chitnis2006}.  A comprehensive survey on malaria infection modeling has been captured in the work by \cite{Mandal2011}. Infectious diseases prevention and control strategies are very important factors and if properly addressed can go a long way in reducing disease infections in a population. Optimal control theory is one of the useful and efficient modeling tools that can help in the modeling framework to better understand the dynamics of the disease and suggest possible prevention and control measures. In recent times, malaria modeling formulation based on optimal control theory  with time dependent controls has attracted much attention \citep[see, e.g,][]{Okosun2013bk, Otieno2016, MiniGosh2015, Agusto2012, Mwanga2014, Torres2013, makinde2014, Okosun2011}.\\

The paper by \cite{Lashari2012a} proposed and numerically analyzed an optimal control model to investigate the effectiveness of time dependent controls on vector-borne diseases. From the numerical results, they \citep{Lashari2012a} observed that, optimal control strategies are effective in reducing the infected populations. \cite{Lashari2012} proposed a time dependent multiple controls to study malaria infection. The authors in \citep{Blayneh2009} formulated a time dependent treatment and prevention optimal control model to study vector-borne diseases. \cite{Okosun2013} presented a controlled dynamical malaria model based on non-linear incidence. \cite{BNKIM2013} developed a dynamic control problem to study malaria transmission using two time dependent control functions. \\

The Ebola infection that claimed many lives in Africa has also been considered with optimal control theory \citep[see, e.g,][]{Rachah2016, Bonyah2016, Grigorieva2015}. Compartmental models for infectious diseases such as TB, HIV and Cholera with time dependent controls have been studied \citep[see, e.g,][]{Moualeu2015, Magombedze2011, ToniBakhtiar2016, Choi2015}. \\

Fractional order derivative has an important property called memory effect and because of that, the theory and application of fractional calculus have widely been used to model dynamical processes in the field of science, engineering and many others \citep[see, e.g,][]{wow2, wow}. \cite{Agrawal2004} developed general optimal control problems that are driven by Riemann-Liouville fractional derivatives. In another research by the same author \citep{Agrawal2008}, he formulated similar optimal control problems with Caputo derivatives and developed a stable numerical scheme for the mathematical model. \\

Several researchers have developed and analyzed optimal control problems whereby the controlled dynamic problems are characterised by Riemann-Liouville derivatives \citep[see, e.g,][]{Frederico, Tricad, Agrawal2010, Agrawal2007}. An enzyme kinetic optimal control problem was proposed and numerically studied by \cite{Basir2015}. The authors in \citep{Ding2012} studeid HIV-Immune System model using fractional optimal control. They deduced the necessary conditions for the optimality system for a general control problem with Caputo derivatives using the concept of calculus of variation and fractional integration by parts formula. Necessary optimality conditions have been derived for an optimal control problem whereby the state model had both first order derivative and non-integer derivative \citep{Pooseh2014}.\\

It has been demonstrated in literature that optimal control problems are appropriate for biological systems and other processes that follow non-linear behavior. Therefore modeling infectious diseases using fractional derivatives together with optimal control theory is a huge advantage. \cite{Okyere2016} extended the malaria infection model  by \cite{Tumwiine2007a, Tumwiine2007b} to include fractional derivatives. They \citep{Okyere2016} investigated asymptotic analysis of the model equilibria and solved the dynamical problem with the numerical scheme developed by \cite{wow4, wow3}. In this work, we propose a fractional optimal control model for malaria infection and numerically simulate the solutions of the constructed control problem. We extend the malaria epidemic model studied by the authors in \citep{Tumwiine2007a, Tumwiine2007b} to become an optimal control problem of which the state and costate equations are characterised by differential equations of fractional orders. Our model formulation is different from previous works on malaria disease modeling since we have incorporated fractional order derivatives into the controlled dynamical system. \\

The paper is organized as follows. Section~\ref{onesec} deals with brief review on the mathematical formulation of fractional optimal control problem. We formulate the controlled malaria model and deduce the necessary conditions for the optimality of the model problem in section~\ref{malaria4}. Section~\ref{twosec} deals with numerical solutions and discussions of simulations results. We then conclude the paper in section~\ref{threesec}.


\section{Optimal Control Problem With Fractional Derivatives}\label{onesec}
This section deals with brief review on the mathematical formulation of fractional optimal control model. We will then apply the modeling concept to formulate our new malaria model in section~\ref{malaria4}. There are many types of fractional derivatives but the most used ones in mathematical modeling and engineering applications are Riemann-Liouville derivative and Caputo derivative \citep[see, e.g,][]{wow2, wow}. For this work, we develop the optimization model with Caputo fractional derivatives..

\begin{definition} \citep{wow, Pooseh2014}
The left Caputo fractional derivative is defined as
\[ _{a}^{c}D^{\alpha}_{t}w(t)=\frac{1}{\Gamma(p-\alpha)}\int_a^t \! \frac{w^{(p)}(x)}{(t-x)^{\alpha+1-p}} \, \mathrm{d}x.\]
\end{definition}
for $p-1<\alpha<p.$\\

\begin{definition} \citep{wow, Pooseh2014}
The right Caputo fractional derivative is defined as
\[ _{t}^{c}D^{\alpha}_{b}w(t)=\frac{(-1)^{p}}{\Gamma(p-\alpha)}\int_t^b \! \frac{w^{(p)}(x)}{(x-t)^{\alpha+1-p}} \, \mathrm{d}x.\]
\end{definition}
for $p-1<\alpha<p,$\\

where $\alpha$ is the order of the Caputo derivatives.

\newpage
In the works by \cite{Agrawal2008, Agrawal2008a}, he formulated an optimal control problem with Caputo fractional derivatives as follows. Find an optimal control $u$ that minimise the objective functional or performance index:

\begin{equation}\label{FOCP1}
J(u)=\int_0^1 \! W(x, u, t) \, \mathrm{d}t.
\end{equation}

Subject to the fractional dynamic constraint

\begin{equation}\label{FOCP2}
_{0}^{c}D^{\alpha}_{t}x(t)=M(x, u, t)
\end{equation}

with initial condition
\begin{equation}\label{FOCP3}
x(0)=x_{o},
\end{equation}

where $W(x, u, t)$ and $M(x, u, t)$ are two arbitrary functions and $x(t)$ is the state variable \citep{Agrawal2004, Agrawal2008}. By using the method of lagrange multipliers, calculus of variation and the integration by parts formula for fractional order derivatives and also following a similar derivation by the same author in \citep{Agrawal2004}, he obtained the necessary conditions or equations (\ref{eric1}-\ref{eric3})  for the fractional order controlled problem as follows:

\begin{align}\label{eric1}
 _{0}^{c}D^{\alpha}_{t}x(t) & =M(x, u, t)
\end{align}

\begin{align}\label{eric2}
_{t}^{c}D^{\alpha}_{1}x(t) & =\frac{\partial W}{\partial x}+\lambda \frac{\partial M}{\partial x}
\end{align}

\begin{align}\label{eric3}
\frac{\partial W}{\partial u}+\lambda \frac{\partial M}{\partial u} & =0
\end{align}

with

\begin{equation}\label{eric4}
\begin{cases}
x(0)=x_o &\\ \\
\lambda (1)=0,
\end{cases}
\end{equation}

where $\lambda$ is the co-state variable called lagrange multiplier.
%

\section{Controlled Malaria Model With Fractional Derivatives}\label{malaria4}
In this section, we will use the obtained results in the optimal control problem formulation briefly reviewed in section~\ref{onesec} to deduce our necessary conditions or equations for the optimality of the optimal control malaria model. Our model formulation seeks to find an optimal control that would minimise the infected human populations and the cost of implementing the optimal control and prevention strategies. In this work, we consider three time dependent controls namely treated bednets, treatment and insecticide spray.\\

The main aim in our present formulation is to find optimal control $U$ that minimise the objective functional:

\begin{equation}\label{FOCP1a}
J(U)=\int_0^b \! W(X, U, t) \, \mathrm{d}t.
\end{equation}

Subject to the fractional dynamic constraints

\begin{equation}\label{FOCP2a}
_{0}^{c}D^{\alpha}_{t}X(t)=M(X, U, t)
\end{equation}

with initial condition
\begin{equation}\label{FOCP3a}
X(0)=X_{o}
\end{equation}

where
\begin{align*}
X(t)& =(S_{H}, I_{H}, R_{H}, S_{V}, I_{V})^{T}\\ \\
X(0)& =(S_{H}(0), I_{H}(0), R_{H}(0), S_{V}(0), I_{V}(0))^{T}\\ \\
U(t)& =(u_{1}, u_{1}, u_{3})^{T}\\ \\
W(X, U, t)& =AI_{H}+\frac{d_{1} u^2_{1}}{2}+\frac{d_{2} u^2_{2}}{2}+\frac{d_{3}u^2_{3}}{2}
\end{align*}

\begin{align}
M(X, U, t)=\left[
        \begin{array}{c}
          \lambda^{\alpha}_{h}N_H-\frac{(1-u_{1})a^{\alpha}bS_H I_V}{N_{H}}+\nu^{\alpha} I_{H}+\gamma^{\alpha} R_H-\mu^{\alpha}_{h}S_{H} \\ \\
          \frac{(1-u_{1})a^{\alpha}bS_H I_V}{N_{H}}-\nu^{\alpha}I_{H}-(r^{\alpha}+\rho^{\alpha}u_{2})I_{H}-\delta^{\alpha}I_{H}-\mu^{\alpha}_{h}I_{H} \\ \\
          (r^{\alpha}+\rho^{\alpha}u_{2})I_{H}-(\gamma^{\alpha}+\mu^{\alpha}_{h})R_{H} \\ \\
          (1-u_{3})\lambda^{\alpha}_{v}N_V-\frac{(1-u_{1})a^{\alpha}cS_V I_H}{N_{H}}-\mu^{\alpha}_{v} S_{V}-\eta^{\alpha} u_{3}S_{V} \\ \\
          \frac{(1-u_{1})a^{\alpha}cS_V I_H}{N_{H}}-\mu^{\alpha}_{v} I_{V}-\eta^{\alpha}u _{3}I_{V}
        \end{array}
      \right]
\end{align}

The variable $X(t)$ is the state vector and $U$ is the control vector. As in \cite{Agrawal2004, Agrawal2008a}, when $\alpha=1$, the formulated malaria controlled dynamical model becomes classical optimal control problem for the uncontrolled malaria model proposed and analyzed by \cite{Tumwiine2007a, Tumwiine2007b}. Also note that, when $\alpha=1$, with control functions $u_1=u_2=u_3=0$, then the dynamic constraint equation~(\ref{FOCP2a}) with the initial condition~(\ref{FOCP3a}) becomes the malaria infection model developed by the authors in \citep{Tumwiine2007a, Tumwiine2007b}\\

The meaning of the model variables and time dependent optimal control functions are defined in table~\ref{tab:table1}. The definitions for the model parameters can be found in table~\ref{tab:table2}

\begin{table}[h!]
  \centering
  \caption{Description of model variables and control functions.}
  \label{tab:table1}
  \begin{tabular}{llr}
    \toprule
    Model variables and control functions & Description \\
    \midrule
    $S_{H}$ & Susceptible human class \\
    $I_{H}$ & Infected human class \\
    $R_{H}$     & partially immune human class \\
    $S_{V}$     & Susceptible mosquito class \\
    $I_{V}$          & Infected mosquito class \\
    $N_{H}$     & Total human population \\
    $N_{V}$          & Total mosquito population\\
    $u_{1}$          & The control on the use of treated bednets\\
    $u_{2}$     & The control on treatment of infected individuals \\
    $u_{3}$          & The control on the use of insecticide spray\\
  \end{tabular}
\end{table}

\begin{table}[h!]
  \centering
  \caption{Description of model parameters.}
  \label{tab:table2}
  \begin{tabular}{llr}
    \toprule
    Model Parameters & Description \\
    \midrule
    $\lambda_{h}$ & natural birth rate of humans\\
    $\lambda_{v}$ & natural birth rate of mosquitoes\\
    $\mu_{h}$     & natural death rate of humans\\
    $\mu_{v}$     & natural death rate of mosquitoes\\
    $a$          &average daily biting rate on human by a single mosquito\\
    $b$          &proportion of bites on human that produce an infection\\
    $\delta$     & death rate due to the disease\\
    $c$          &probability that a mosquito becomes infectious\\
    $\nu$        & recovery rate of human hosts from the disease\\
    $\gamma$      &    rate of loss of immunity in human hosts\\
    $r$         &rate at which human hosts acquire immunity\\
    $\rho$      &  rate constant due to treatment   \\
    $\eta$         & rate constant due to insecticide spray\\
    $A,\ d_1, \ d_2,\ d_3$ & weight constants\\
    \bottomrule
  \end{tabular}
\end{table}

\newpage
It is important to mention that, rigorous proofs on the derivation of necessary conditions for the optimality of various fractional dynamical systems can be found in literature \citep[see, e.g,][]{Pooseh2014, Ding2012, Agrawal2004}. Therefore in this paper, using the ideas by the authors in \citep[see, e.g,][]{Agrawal2008, Agrawal2008a, Pooseh2014, Ding2012, Agrawal2004}, we state the necessary conditions (\ref{eric1a}-\ref{eric3a}) for our optimization problem as follows:

\begin{align}\label{eric1a}
 _{0}^{c}D^{\alpha}_{t}X(t) & = M(X, U, t)
\end{align}

\begin{align}\label{eric2a}
_{t}^{c}D^{\alpha}_{b}\lambda (t) & =\frac{\partial W}{\partial X}+\lambda^{T} \frac{\partial M}{\partial X}
\end{align}

\begin{align}\label{eric3a}
\frac{\partial W}{\partial U}+\lambda^{T} \frac{\partial M}{\partial U} & =0
\end{align}
and
\begin{equation}\label{eric4a}
\begin{cases}
X(0)=X_o &\\ 
\lambda (b)=0
\end{cases}
\end{equation}

where $\lambda(t) =(\lambda_{1}, \lambda_{2}, \lambda_{3}, \lambda_{4}, \lambda_{5})^{T}$ is the co-state vector.\\

Now using the compact form of the necessary conditions stated above, we obtain the optimality system for the fractional optimal control malaria model in an expanded form as follows:

\begin{equation}\label{SIRmalaria2}
\begin{array}{lllll}
\displaystyle _{0}^{c}D^{\alpha}_{t}S_{H} &=& \displaystyle {\lambda^{\alpha}_{h}N_H-\frac{(1-u_{1})a^{\alpha}bS_H I_V}{N_{H}}+\nu^{\alpha} I_{H}+\gamma^{\alpha} R_H-\mu^{\alpha}_{h}S_{H}} ,\qquad \\[15pt]
\displaystyle _{0}^{c}D^{\alpha}_{t}I_{H} &=& \displaystyle {\frac{(1-u_{1})a^{\alpha}bS_H I_V}{N_{H}}-\nu^{\alpha}I_{H}-(r^{\alpha}+\rho^{\alpha}u_{2})I_{H}-\delta^{\alpha}I_{H}-\mu^{\alpha}_{h}I_{H}}, \qquad  \\[15pt]
\displaystyle _{0}^{c}D^{\alpha}_{t}R_{H} &=& \displaystyle {(r^{\alpha}+\rho^{\alpha}u_{2})I_{H}-(\gamma^{\alpha}+\mu^{\alpha}_{h})R_{H}}, \qquad  \\[15pt]
\displaystyle _{0}^{c}D^{\alpha}_{t}S_{V}&=& \displaystyle {(1-u_{3})\lambda^{\alpha}_{v}N_V-\frac{(1-u_{1})a^{\alpha}cS_V I_H}{N_{H}}-\mu^{\alpha}_{v} S_{V}-\eta^{\alpha} u_{3}S_{V}} ,\qquad \\[15pt]
\displaystyle _{0}^{c}D^{\alpha}_{t}I_{V} &=& \displaystyle {\frac{(1-u_{1})a^{\alpha}cS_V I_H}{N_{H}}-\mu^{\alpha}_{v} I_{V}-\eta^{\alpha}u _{3}I_{V}}, \qquad
\end{array}
\end{equation}

\begin{equation}\label{eric5c}
\begin{cases}
_{t}^{c}D^{\alpha}_{b}\lambda_{1}= & \lambda_{1}\left[\lambda^{\alpha}_{h}-\dfrac{(1-u_{1})(I_{H}+R_{H})a^{\alpha}bI_{V}}{N^2_{H}}-\mu^{\alpha}_{h}\right]
     +\lambda_{2}\left[\dfrac{(1-u_{1})(I_{H}+R_{H})a^{\alpha}bI_{V}}{N^2_{H}}\right]\\ \\
     &-(\lambda_{5}-\lambda_{4})\left[\dfrac{(1-u_{1}))a^{\alpha}cS_{H}I_{H}}{N^2_{H}}\right]\\ \\
  _{t}^{c}D^{\alpha}_{b}\lambda_{2}= & A+\lambda_{1}\left[\lambda^{\alpha}_{h}+\dfrac{(1-u_{1})a^{\alpha}bS_{H}I_{V}}{N^2_{H}}+\nu^{\alpha}\right]
     -(\lambda_{4}-\lambda_{5})\left[\dfrac{(1-u_{1})(S_{H}+R_{H})a^{\alpha}cS_{V}}{N^2_{H}}\right]\\\\
     &-\lambda_{2}\left[\dfrac{(1-u_{1})a^{\alpha}bS_{H}I_{V}}{N^2_{H}}+\nu^{\alpha}+(r^{\alpha}+\rho^{\alpha}u_{2})+\delta^{\alpha}+\mu^{\alpha}_{h}\right]+\lambda_{3}\left[r^{\alpha}+\rho^{\alpha}u_{2}\right]\\\\    _{t}^{c}D^{\alpha}_{b}\lambda_{3}= & \lambda_{1}\left[\lambda^{\alpha}_{h}+\dfrac{(1-u_{1})a^{\alpha}bS_{H}I_{V}}{N^2_{H}}+\gamma^{\alpha}_{h}\right]
     -\lambda_{2}\left[\dfrac{(1-u_{1})a^{\alpha}bS_{H}I_{V}}{N^2_{H}}\right]-\lambda_{3}(\gamma^{\alpha}+\mu^{\alpha}_{h})\\\\
     &-(\lambda_{5}-\lambda_{4})\left[\dfrac{(1-u_{1})a^{\alpha}cS_{V}I_{H}}{N^2_{H}}\right]\\ \\
     _{t}^{c}D^{\alpha}_{b}\lambda_{4}= & \lambda_{4}\left[(1-u_{3})\lambda^{\alpha}_{v}+\dfrac{(1-u_{1})a^{\alpha}cI_{H}}{N_{H}}-\mu^{\alpha}_{v}-\eta^{\alpha}u_{3}\right]
     + \lambda_{5}\left[\dfrac{(1-u_{1})a^{\alpha}cI_{H}}{N_{H}}\right]\\ \\
     _{t}^{c}D^{\alpha}_{b}\lambda_{5}= & (\lambda_{2}-\lambda_{1})\left[\dfrac{(1-u_{1})a^{\alpha}bS_{H}}{N_{H}}\right]
     + \lambda_{4}\left[(1-u_{3})\lambda^{\alpha}_{v}\right]-\lambda_{5}\left[\mu^{\alpha}_{v}+\eta^{\alpha}u_{3}\right]
     \end{cases}
\end{equation}\\


\begin{equation}\label{eric5}
\begin{cases}
u_{1}=\dfrac{(\lambda_{2}-\lambda_{1})\frac{a^{\alpha}bS_H I_V}{N_{H}}+(\lambda_{5}-\lambda_{4})\frac{a^{\alpha}cS_V I_H}{N_{H}}}{B_{1}} &\\ \\
u_{2} = \frac{(\lambda_{2}-\lambda_{3})r^{\alpha}_{o}I_{H}}{B_2}&\\ \\
u_{3} =\frac{\lambda_{4}(\lambda^{\alpha}_{v}N_{V}+\eta^{\alpha}S_{V})+\lambda_{5}\eta^{\alpha}I_{V}}{B_3}
\end{cases}
\end{equation}\\

\begin{equation}\label{eric6}
\begin{cases}
S_{H}(0)=S_{H_{O}}, I_{H}(0)=I_{H_{O}}, R_{H}(0)=R_{H_{O}}, S_{V}(0)=S_{V_{O}}, I_{V}(0)=I_{V_{O}}&\\ \\
\lambda_1 (b)=0, \lambda_2 (b)=0, \lambda_3 (b)=0, \lambda_4 (b)=0, \lambda_5 (b)=0,
\end{cases}
\end{equation}


%
%

%
%



\section{Numerical Simulations and Discussion}\label{twosec}
The book by \cite{Lenhart2007} described the numerical solutions of optimal control models for biological systems. In their work, they constructed the forward-backward sweep algorithm with Runge-Kutta method to solve the optimality system of biological models. In this work, the forward-backward algorithm with the generalized euler method is applied to numerically compute the solutions of the optimality system. The generalized euler method is one of the numerical schemes for fractional order models and has recently been applied together with the forward-backward algorithm in a TB model \citep{Sweilam2016}.\\

The following initial conditions and model parameters values are used for the numerical simulations:
$ S_{H}(0)=800,\ I_{H}(0)=200, \ R_{H}(0)=20, \ S_{V}(0)=1000,\ I_{V}(0)=500$, $A=100,\ d_1=70, \ d_2=130,\ d_3=40,\ a=0.29,\ b=0.75,\ c=0.75,$ $r=0.05\, \gamma=0.000017,\ \rho=0.7,\ \delta=0.02,\ \nu=0.0022,\ \lambda_h=0.0015875,\ \lambda_v=0.071,\ \mu_h=0.00004,\ \mu_v=0.1429,\ \eta=0.25$  and $\alpha=1,\ 0.99, \ 0.95,\ 0.90.$ \\

We considered the following optimal control strategies for our numerical illustrations: Treated bednets control only ($u_1\neq 0,\ u_2= 0,\ u_3=0$), treatment control only ($u_1=0,\ u_2\neq 0,\ u_3=0$), insecticide spray control only ($u_1= 0,\ u_2= 0,\ u_3\neq0$), treated bednets and treatment controls ($u_1\neq 0,\ u_2\neq 0,\ u_3=0$),
 treated bednets and insecticide spray controls ($u_1\neq 0,\ u_2= 0,\ u_3\neq 0$), treatment and insecticide spray controls ($u_1= 0,\ u_2\neq0,\ u_3\neq 0$),
  and all the three controls ($u_1\neq 0,\ u_2\neq 0,\ u_3\neq 0$).\\

 From the plots generated, it is clear that reduction of infected human and mosquito populations occurs when $\alpha=1$ (integer order) with time dependent controls but that of non-integer orders ($\alpha=0.99, 0.95, 0.90$) with time dependent controls are much more significant.

\begin{figure}[!ht]
\centering
\subfigure[]{
\includegraphics[scale=0.51]{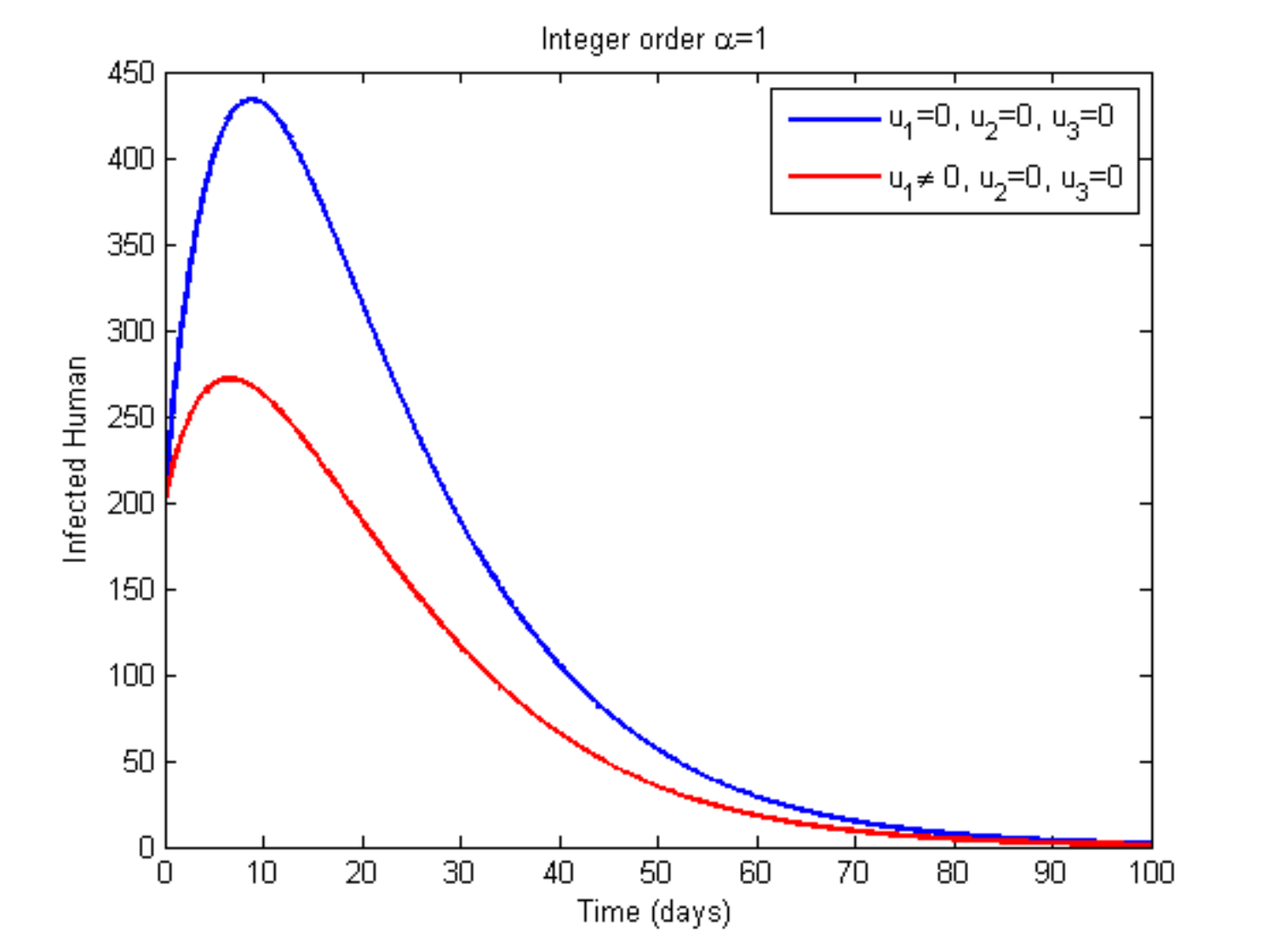}}\hfil
\subfigure[]{
\includegraphics[scale=0.51]{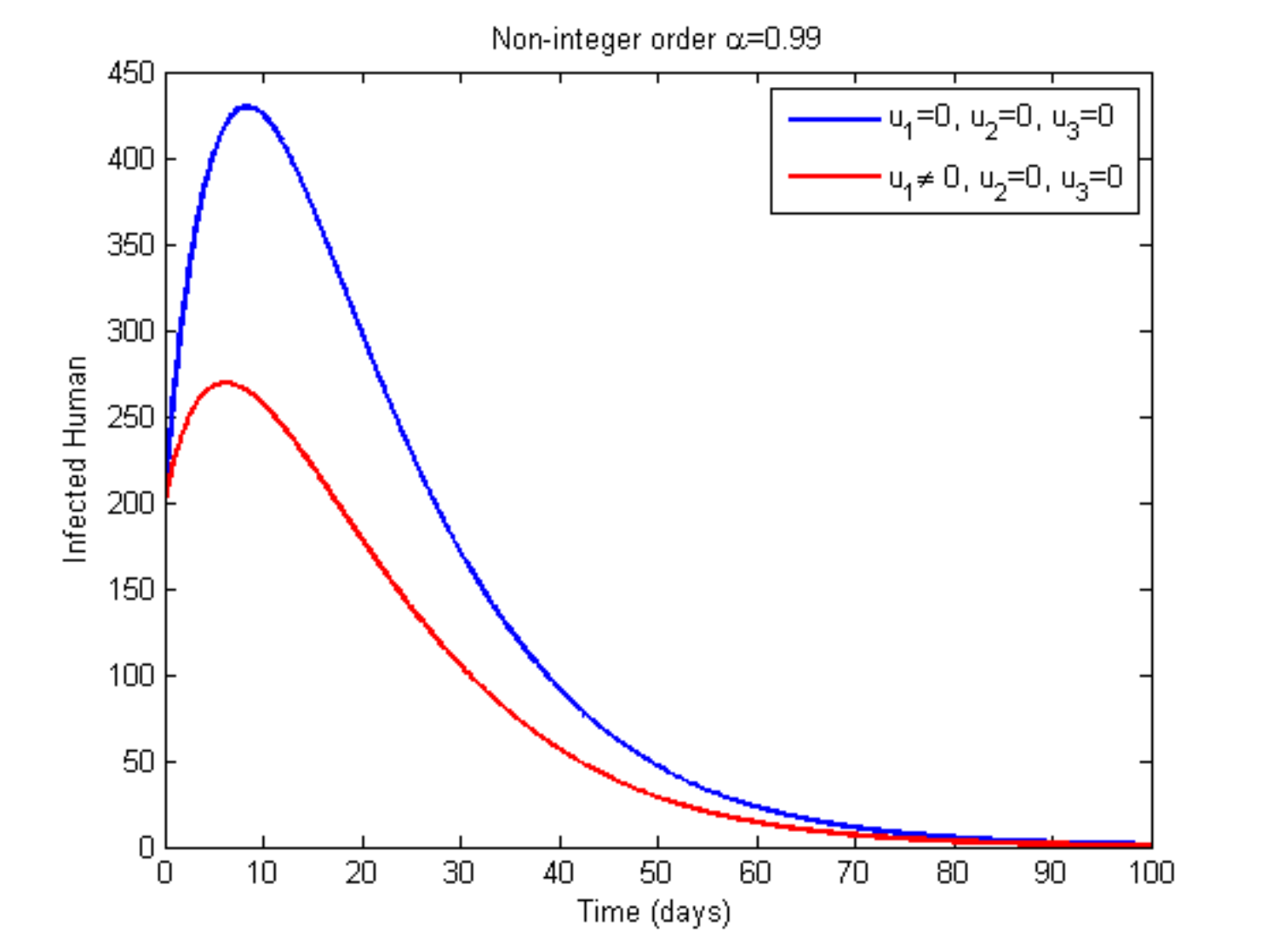}}\hfil
\subfigure[]{
\includegraphics[scale=0.51]{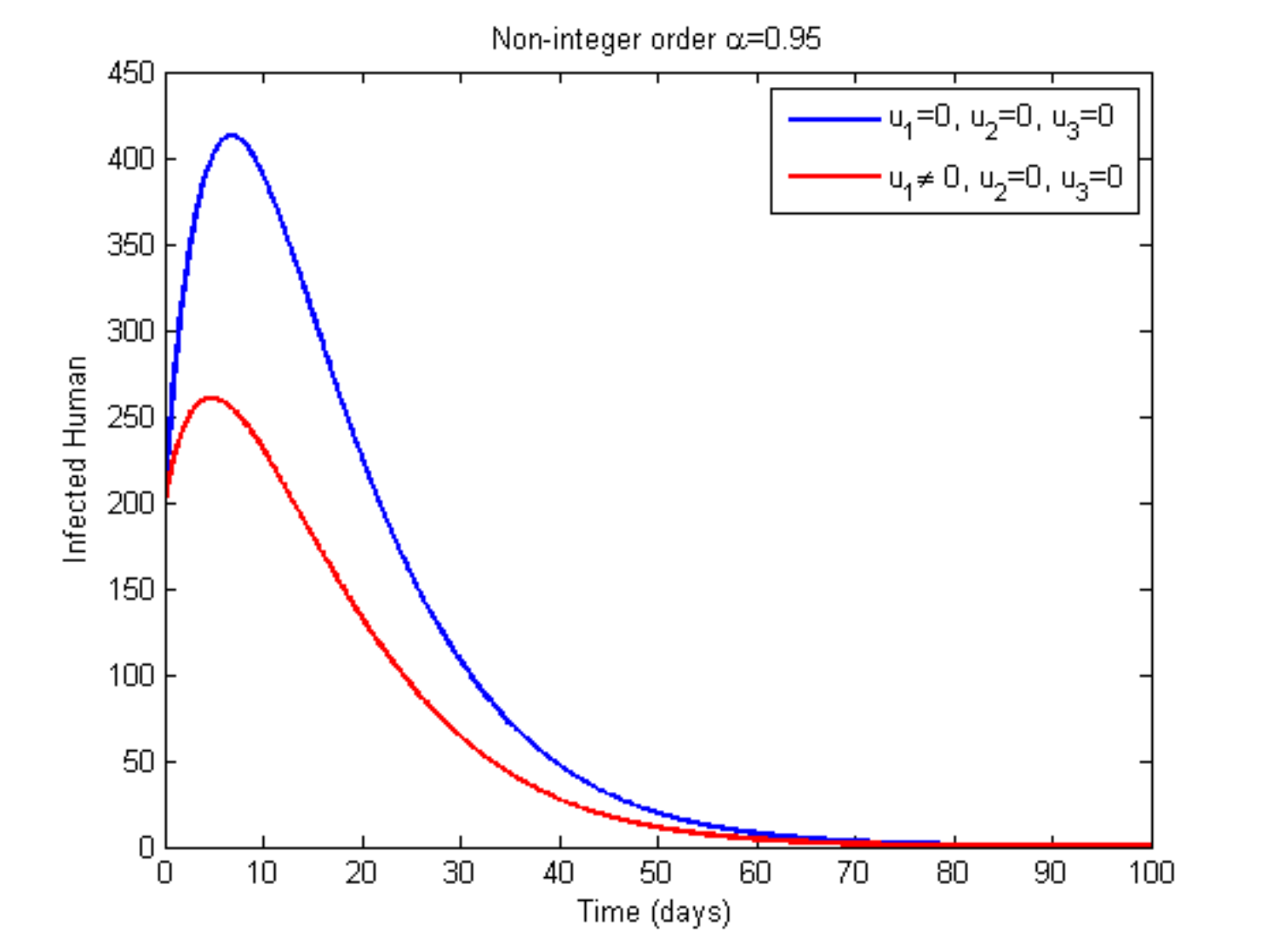}}\hfil
\subfigure[]{\includegraphics[scale=0.51]{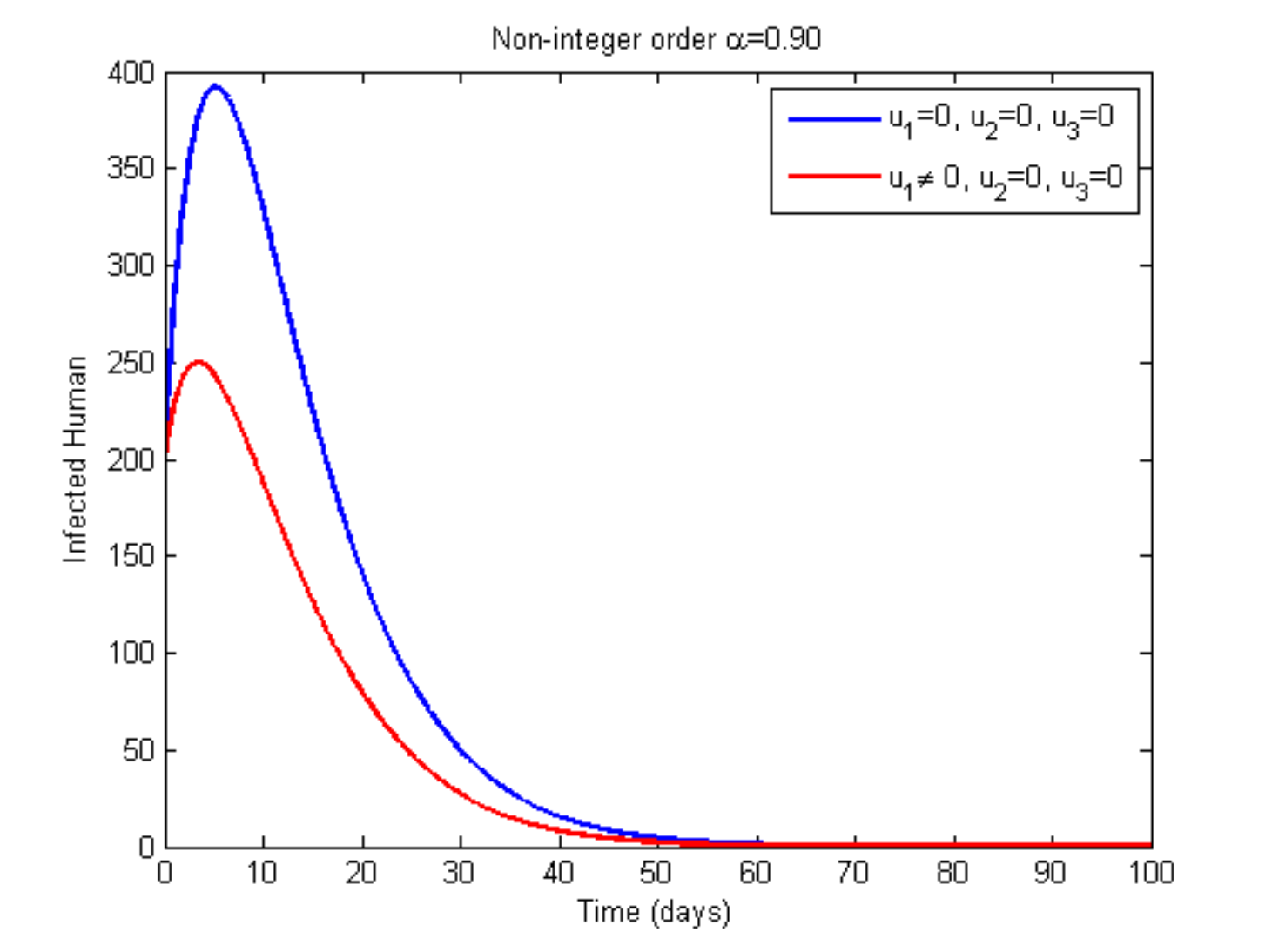}}
\caption{Numerical solutions of infected human host with treated bednets control for $\alpha=1, 0.99, 0.95, 0.90$}
\label{fg1}
\end{figure}

\begin{figure}[!ht]
\centering
\subfigure[]{
\includegraphics[scale=0.51]{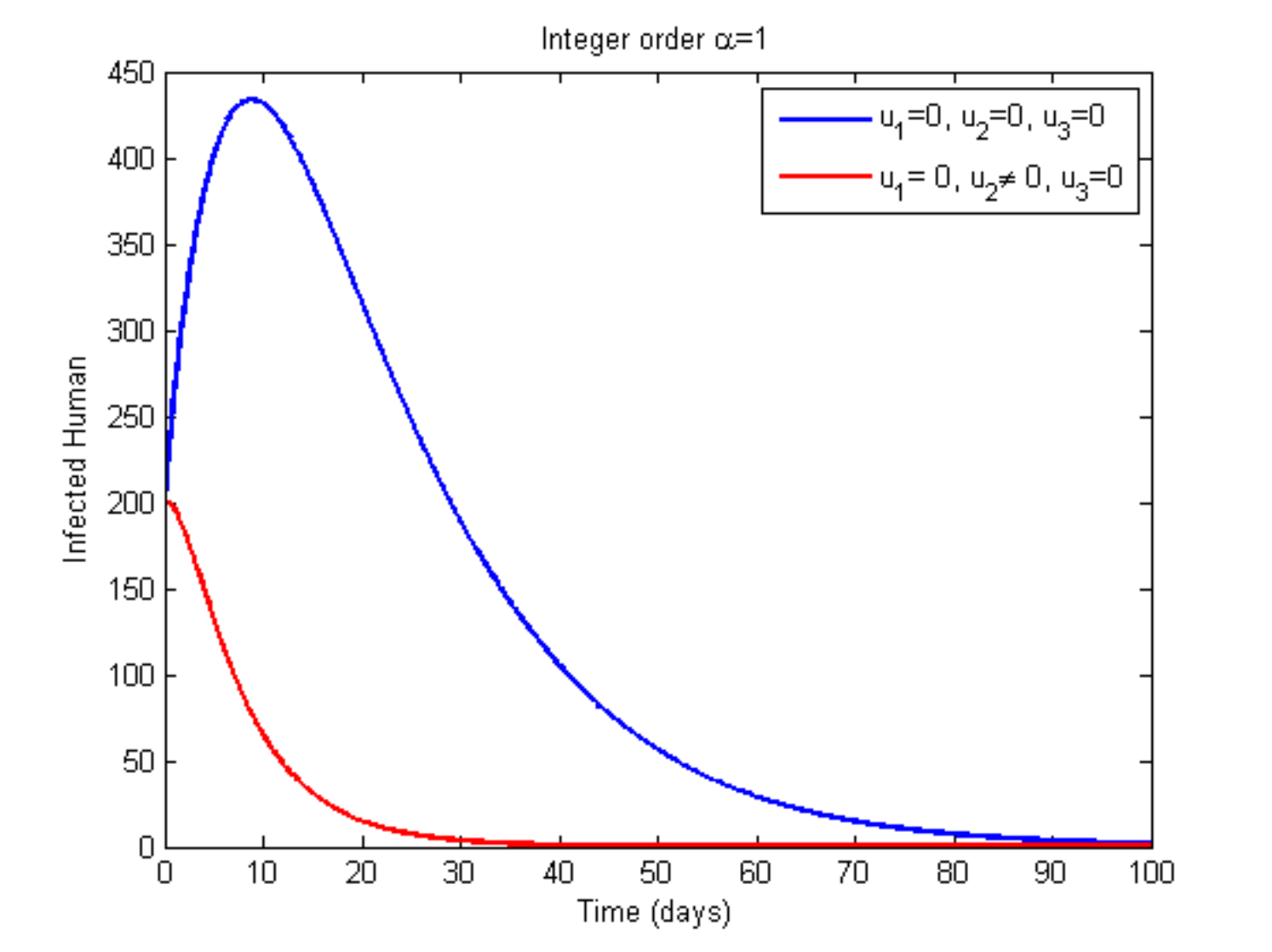}}\hfil
\subfigure[]{
\includegraphics[scale=0.51]{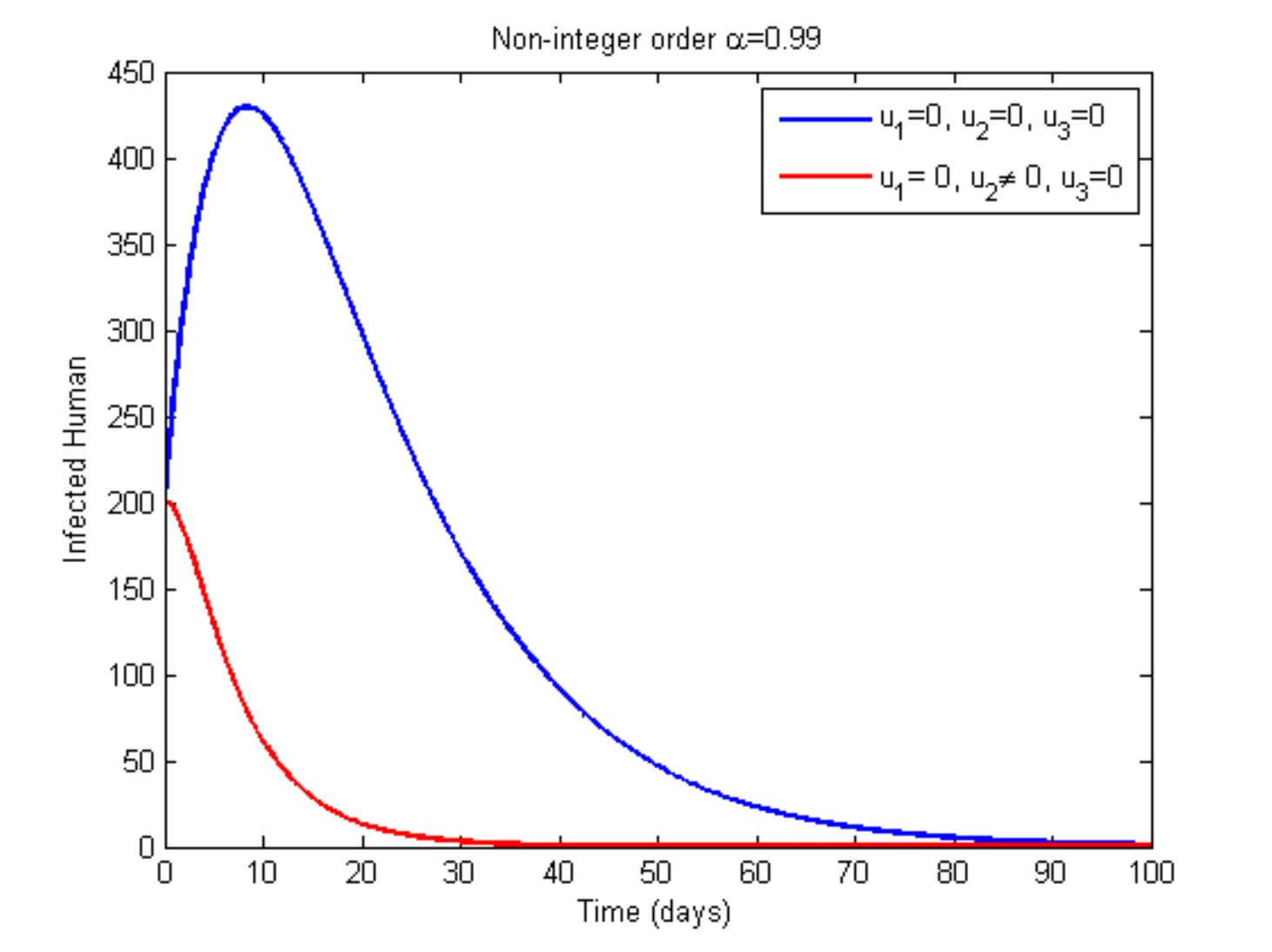}}\hfil
\subfigure[]{
\includegraphics[scale=0.51]{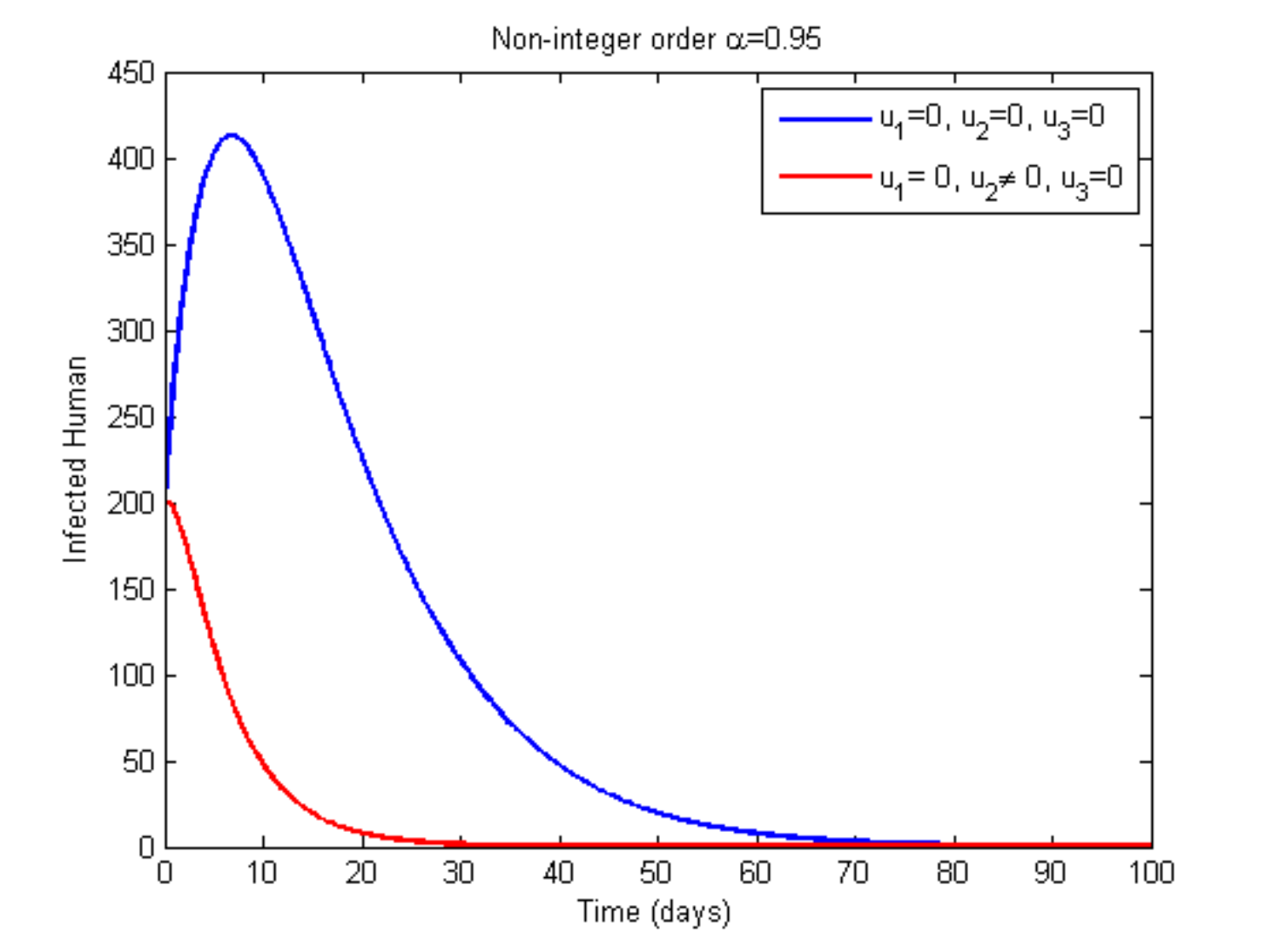}}\hfil
\subfigure[]{\includegraphics[scale=0.51]{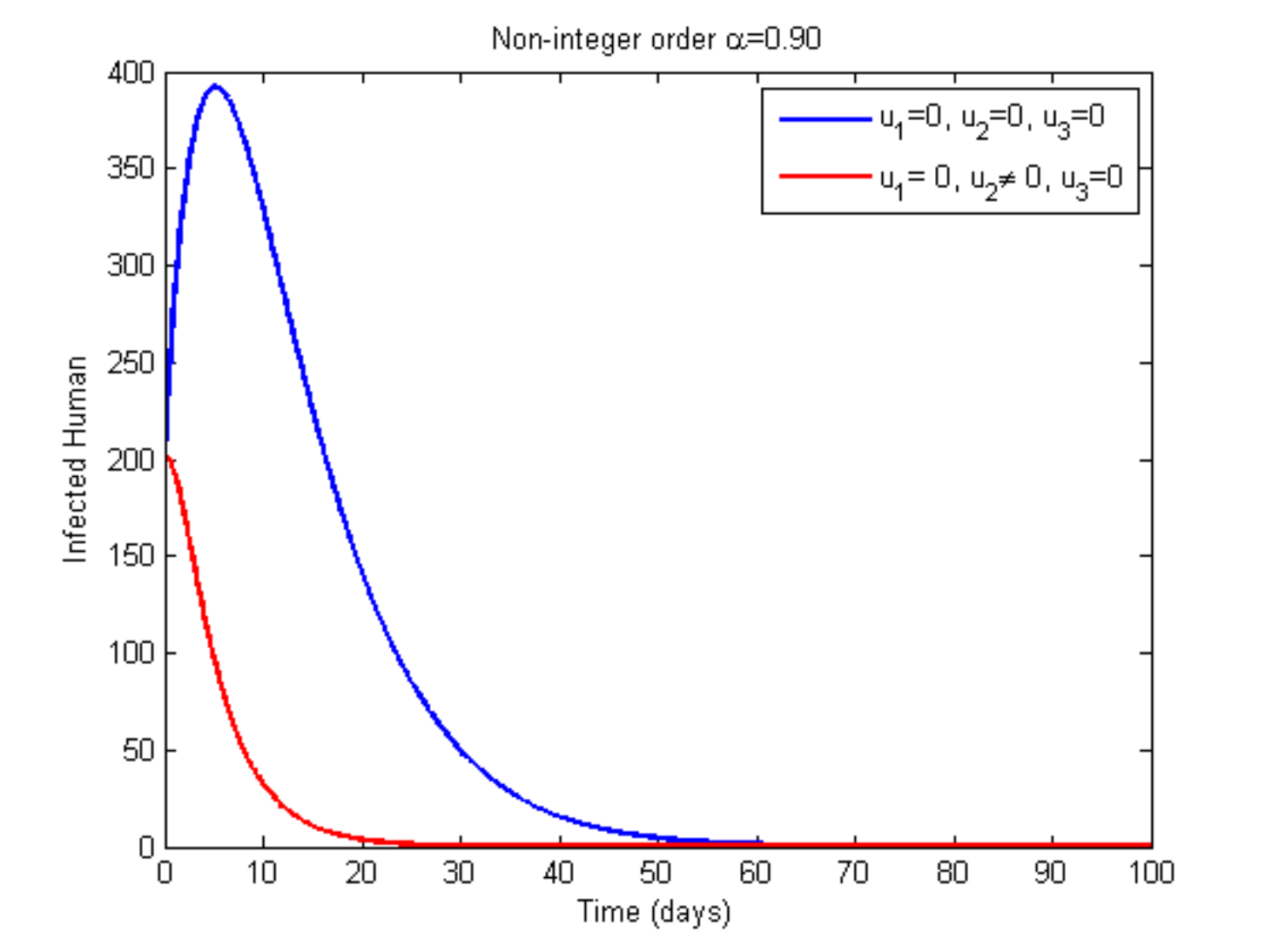}}
\caption{Numerical solutions of infected human host with treatment control for $\alpha=1, 0.99, 0.95, 0.90$}
\label{fg2}
\end{figure}

\begin{figure}[!ht]
\centering
\subfigure[]{
\includegraphics[scale=0.51]{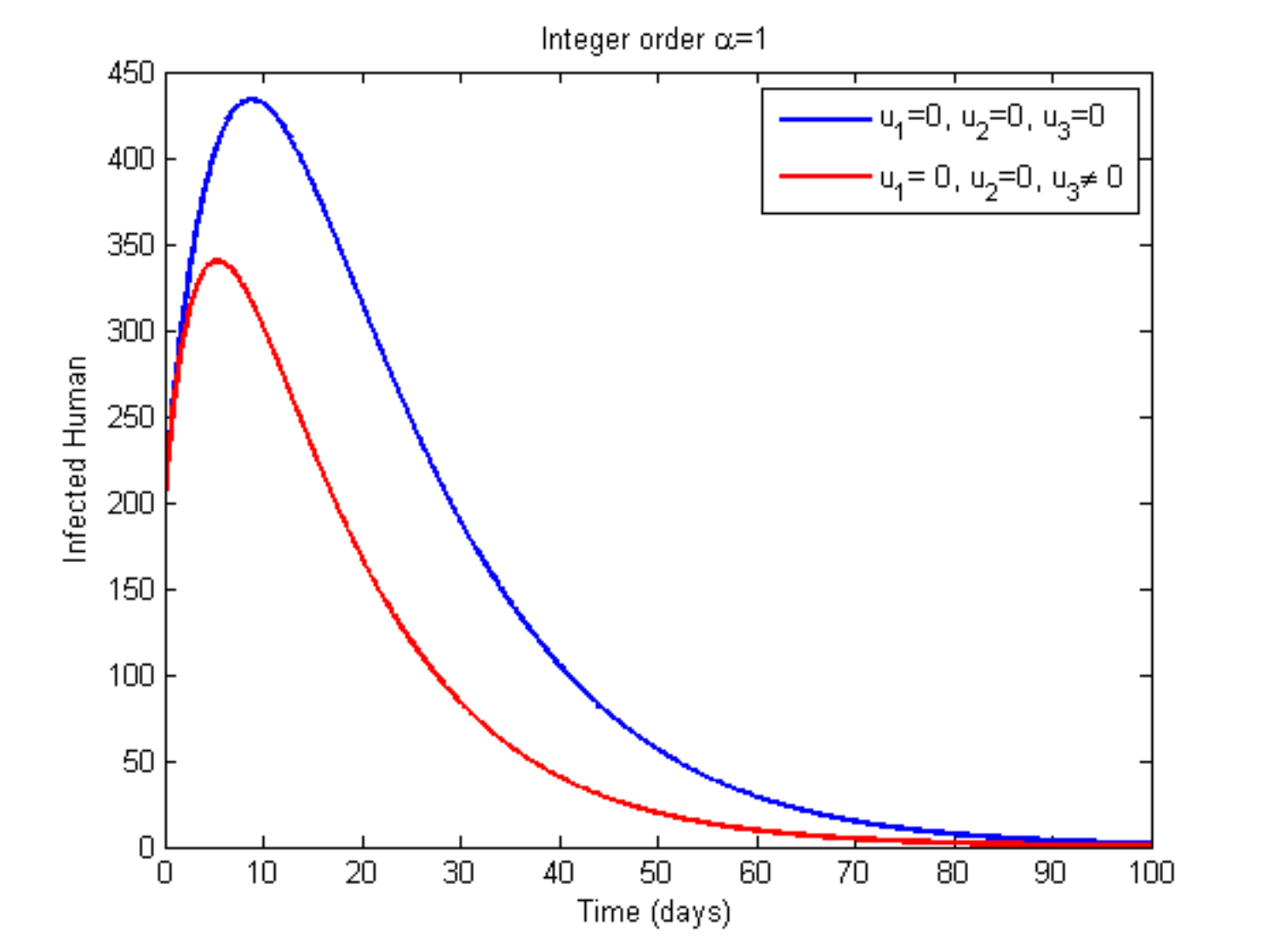}}\hfil
\subfigure[]{
\includegraphics[scale=0.51]{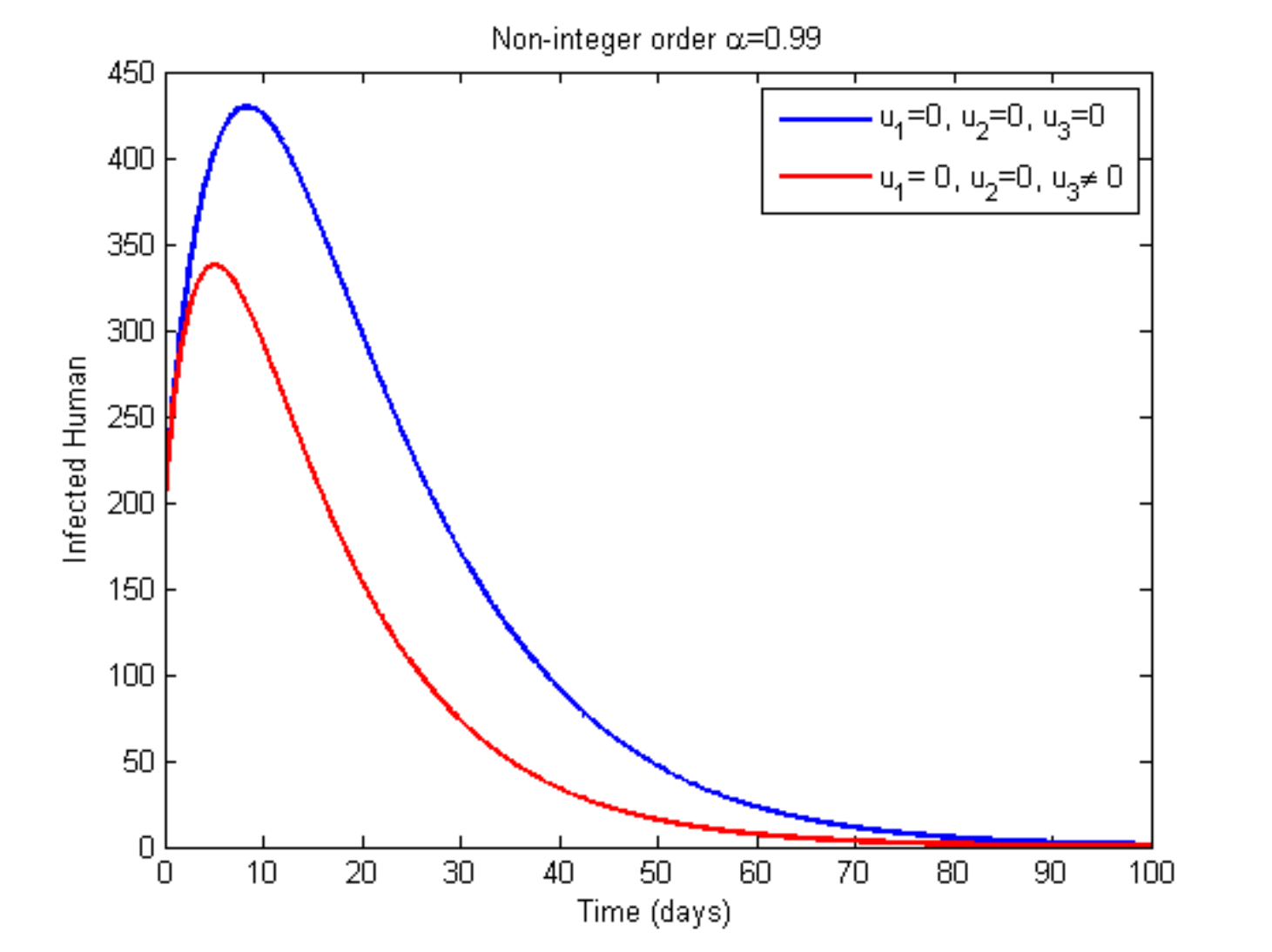}}\hfil
\subfigure[]{
\includegraphics[scale=0.51]{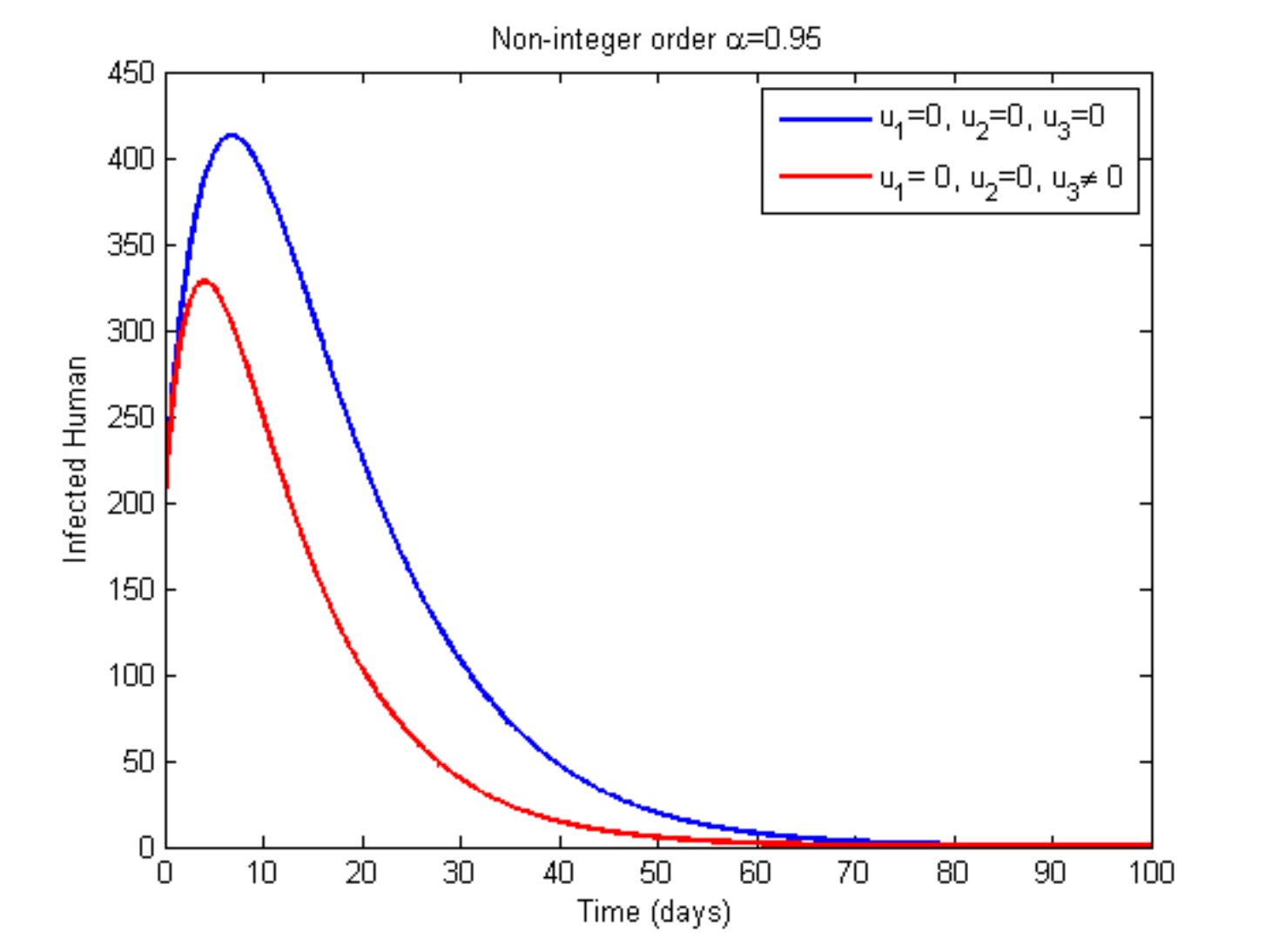}}\hfil
\subfigure[]{\includegraphics[scale=0.51]{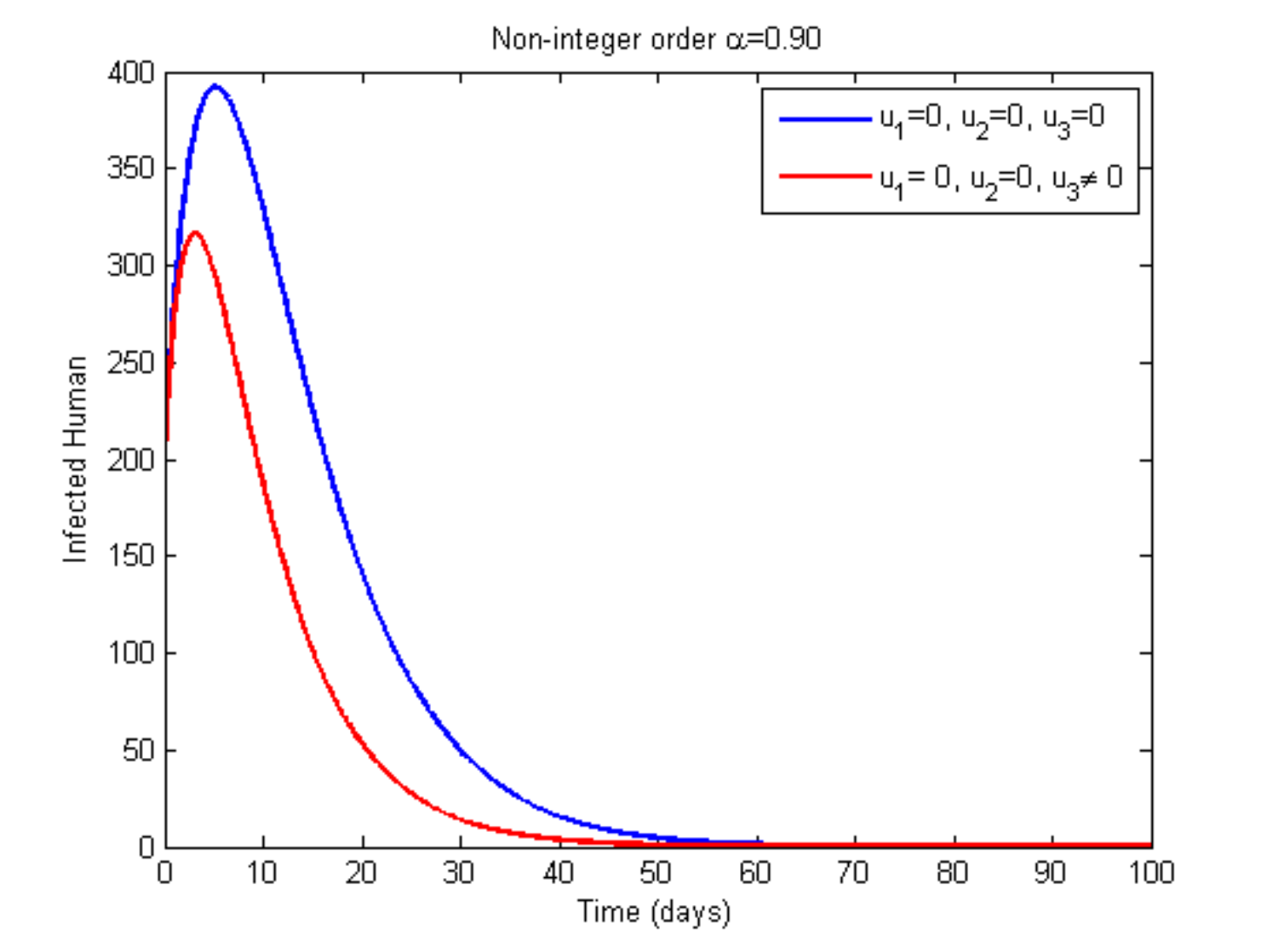}}
\caption{Numerical solutions of infected human host with insecticide spray control for $\alpha=1, 0.99, 0.95, 0.90$}
\label{fg3}
\end{figure}

\begin{figure}[!ht]
\centering
\subfigure[]{
\includegraphics[scale=0.51]{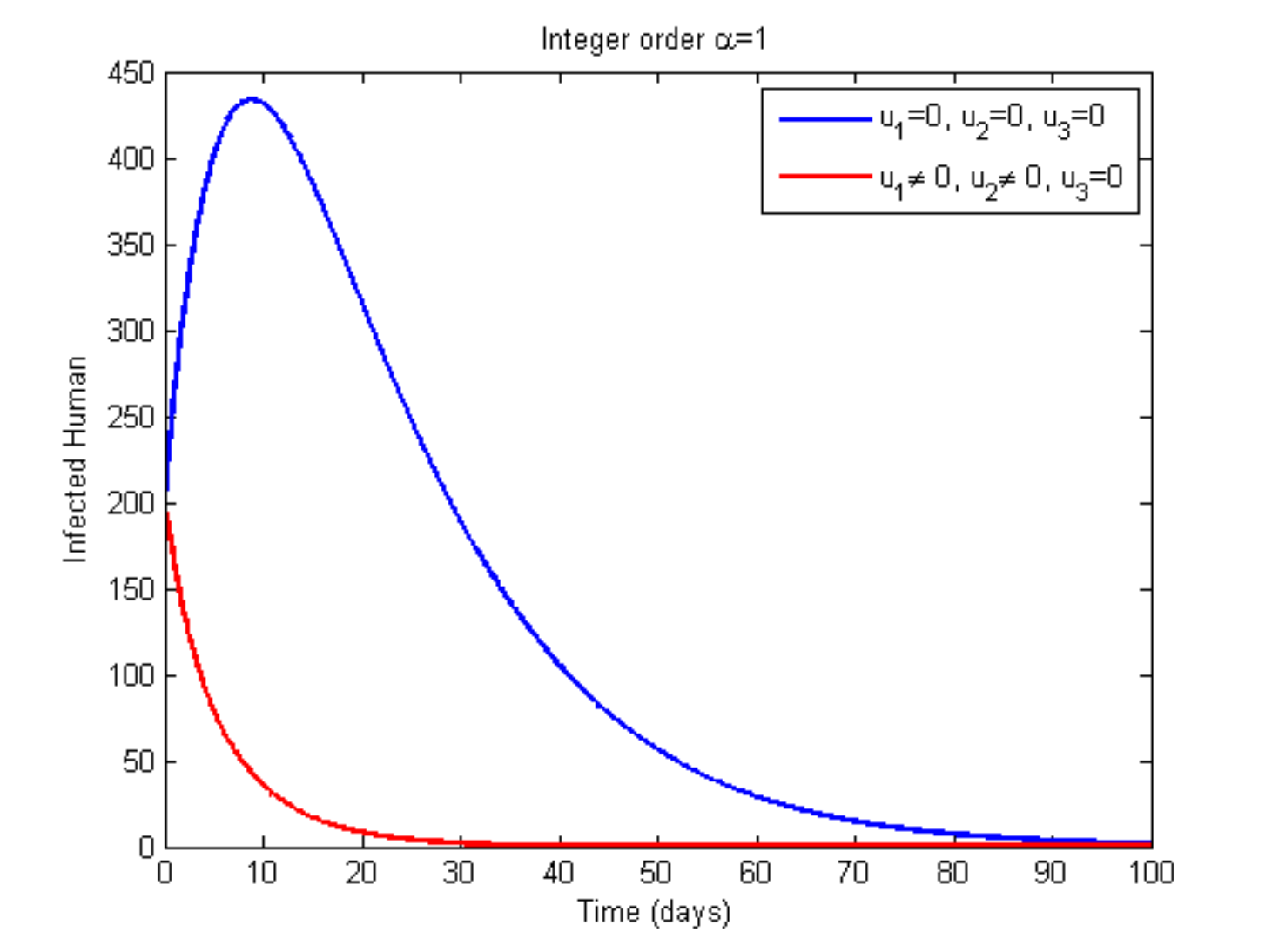}}\hfil
\subfigure[]{
\includegraphics[scale=0.51]{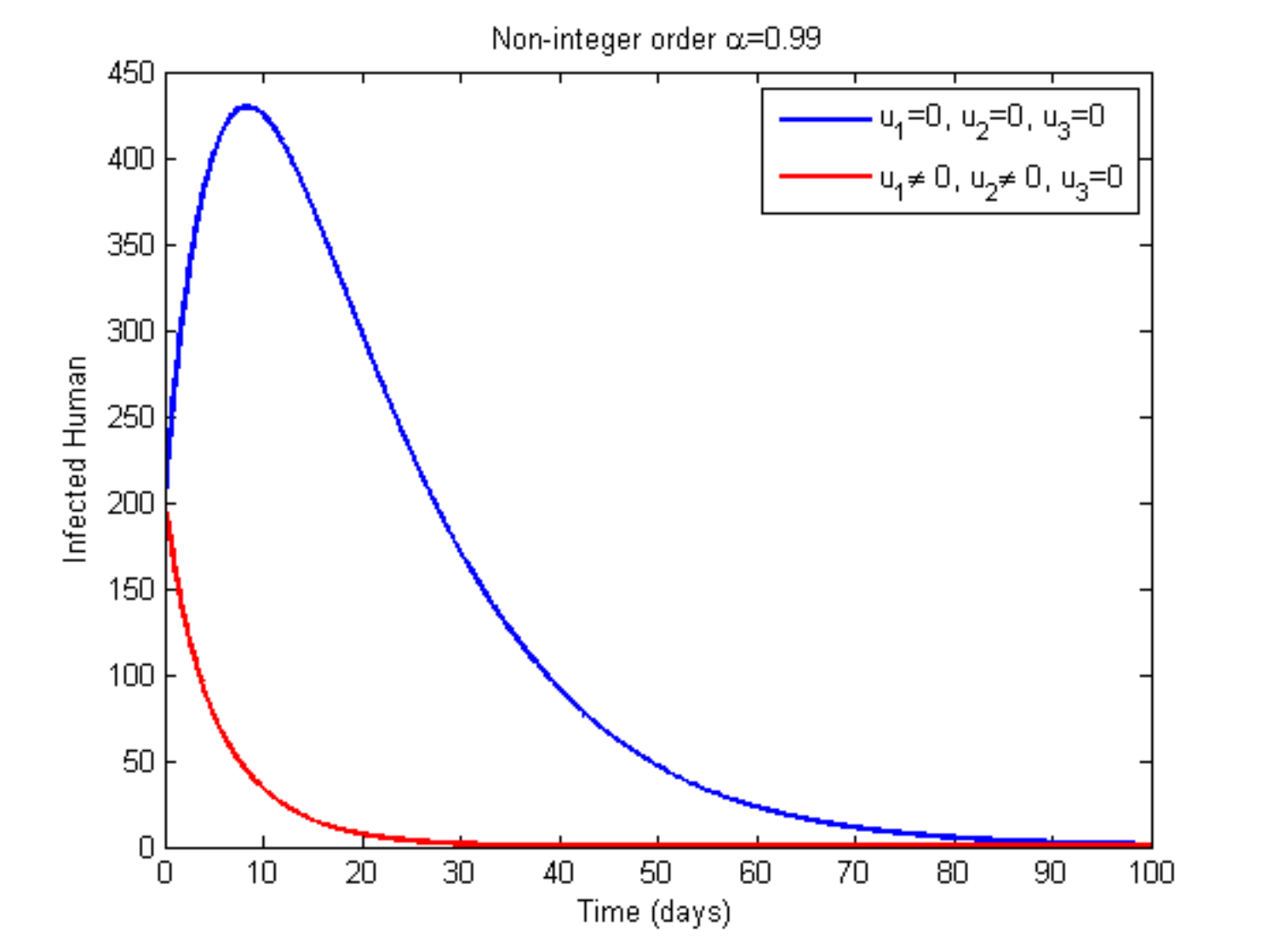}}\hfil
\subfigure[]{
\includegraphics[scale=0.51]{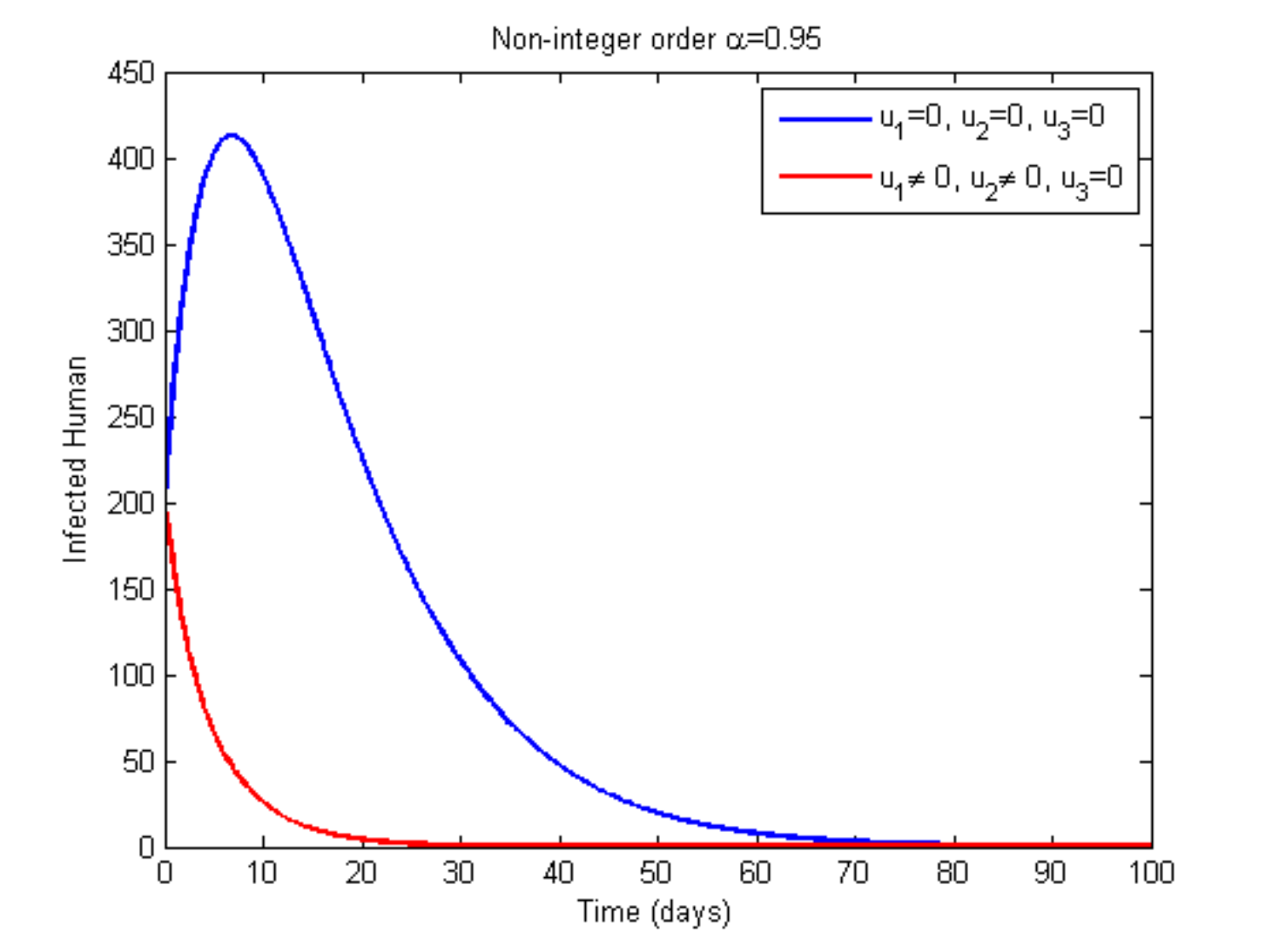}}\hfil
\subfigure[]{\includegraphics[scale=0.51]{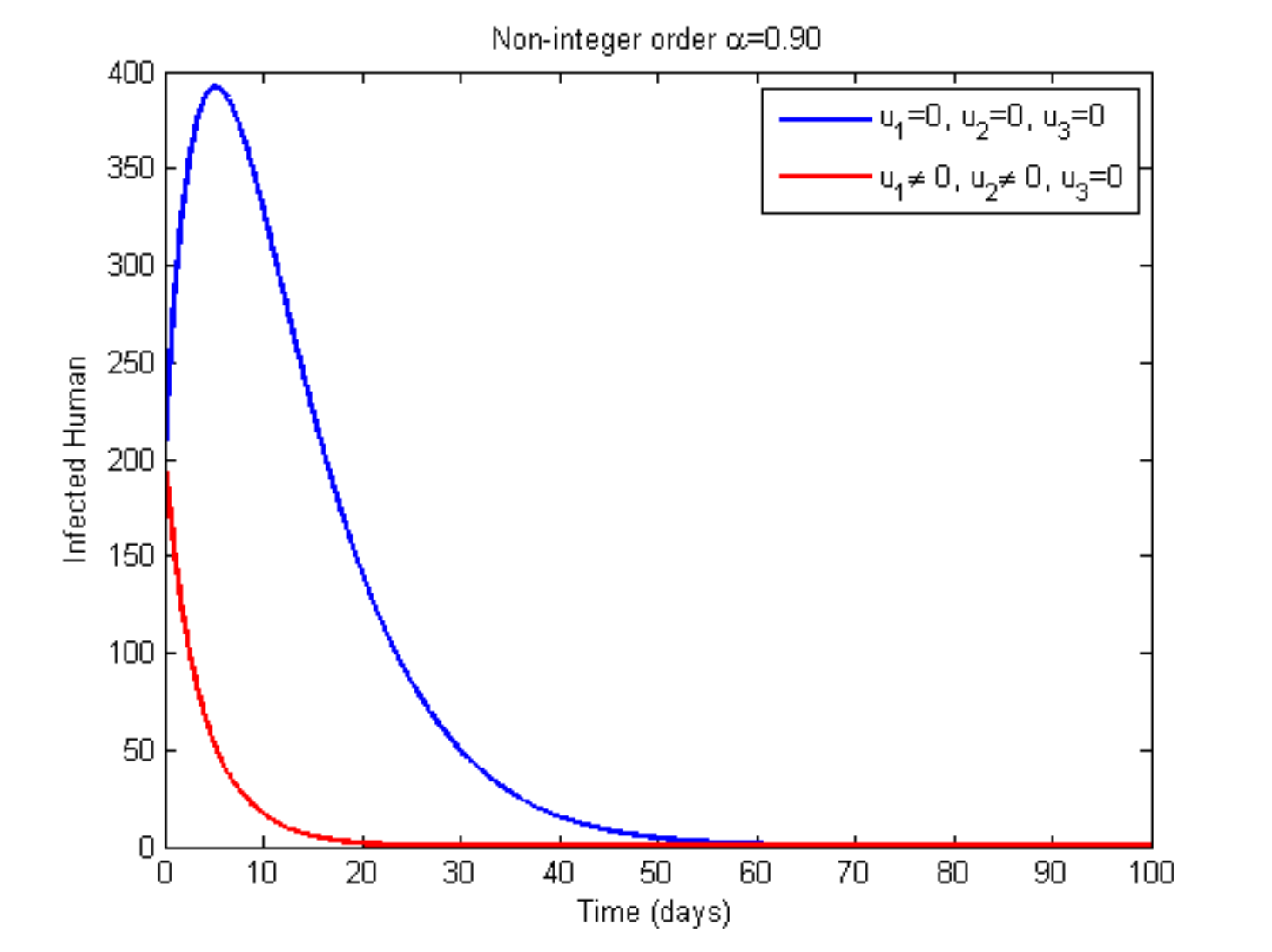}}
\caption{Numerical solutions of infected human host with treated bednets and treatment controls for $\alpha=1, 0.99, 0.95, 0.90$}
\label{fg4}
\end{figure}

\begin{figure}[!ht]
\centering
\subfigure[]{
\includegraphics[scale=0.5]{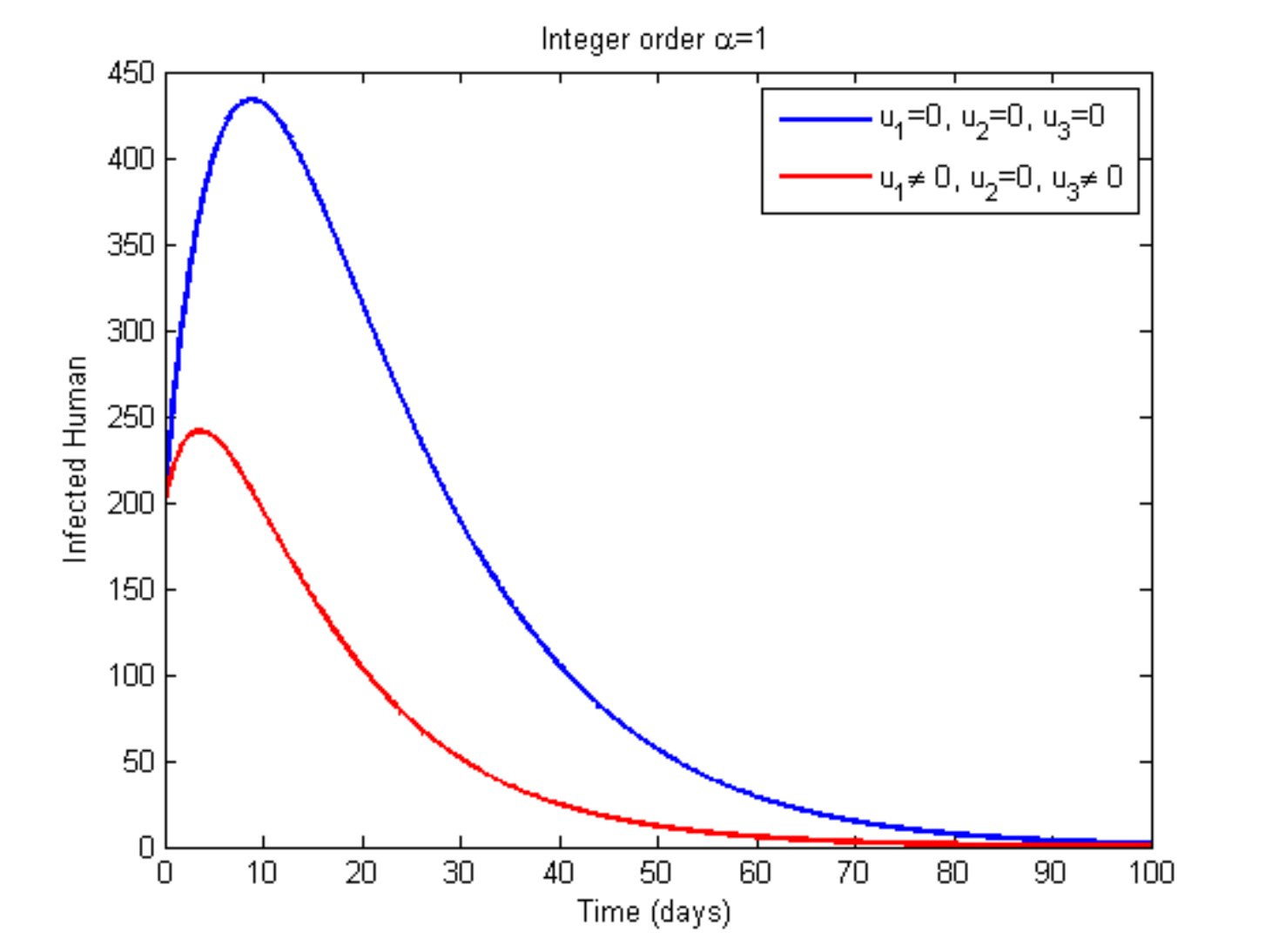}}\hfil
\subfigure[]{
\includegraphics[scale=0.5]{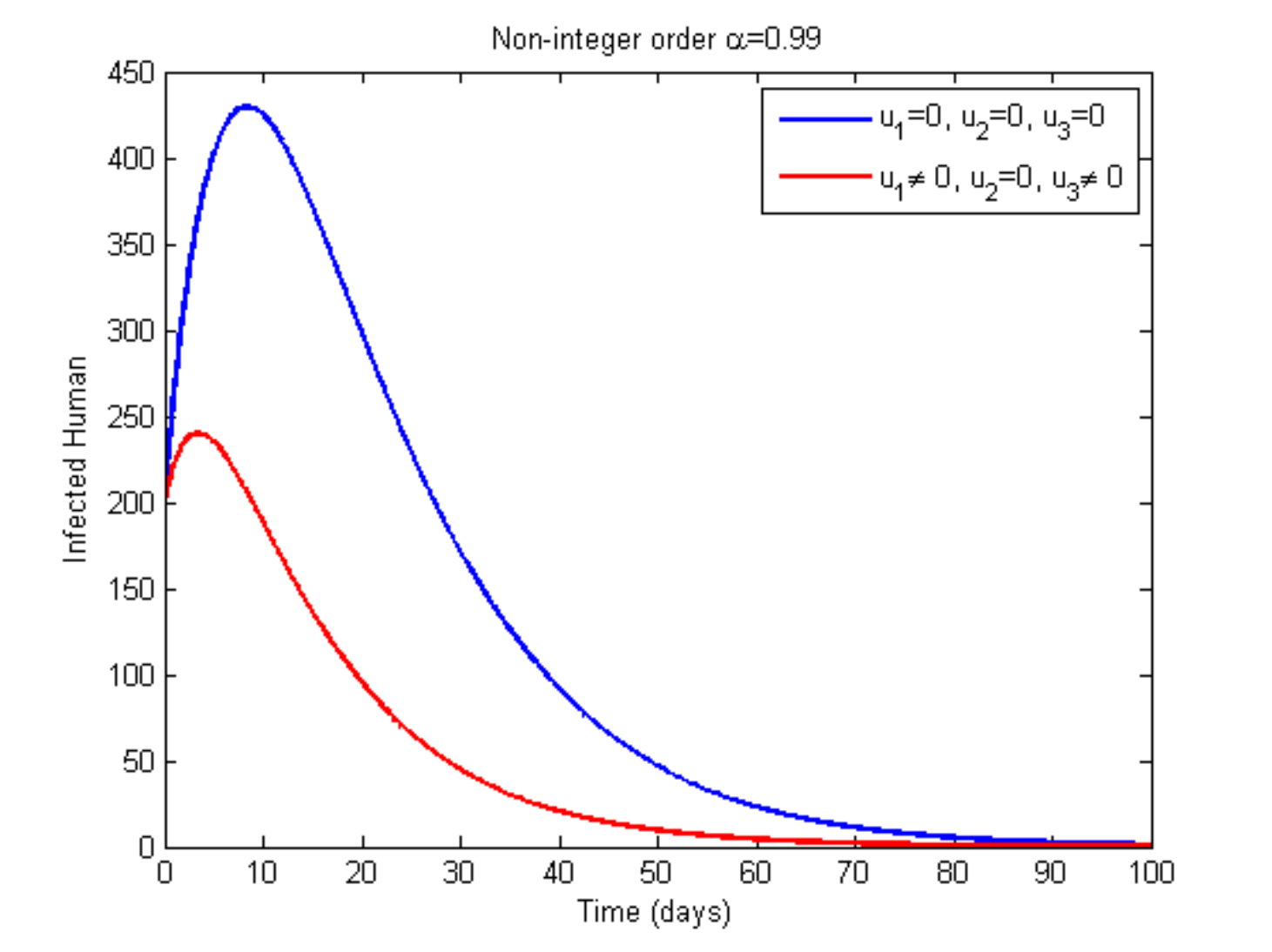}}\hfil
\subfigure[]{
\includegraphics[scale=0.5]{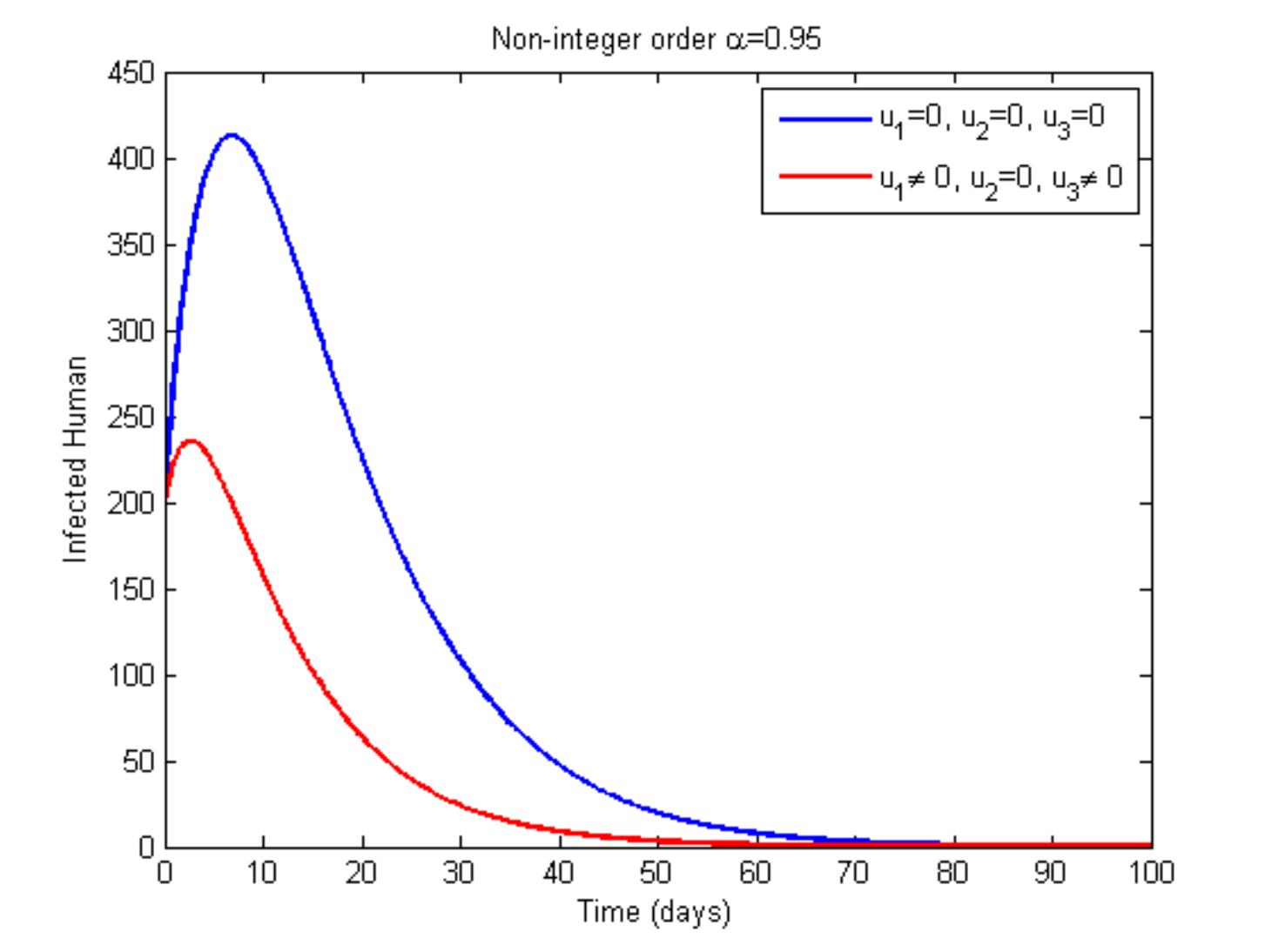}}\hfil
\subfigure[]{\includegraphics[scale=0.5]{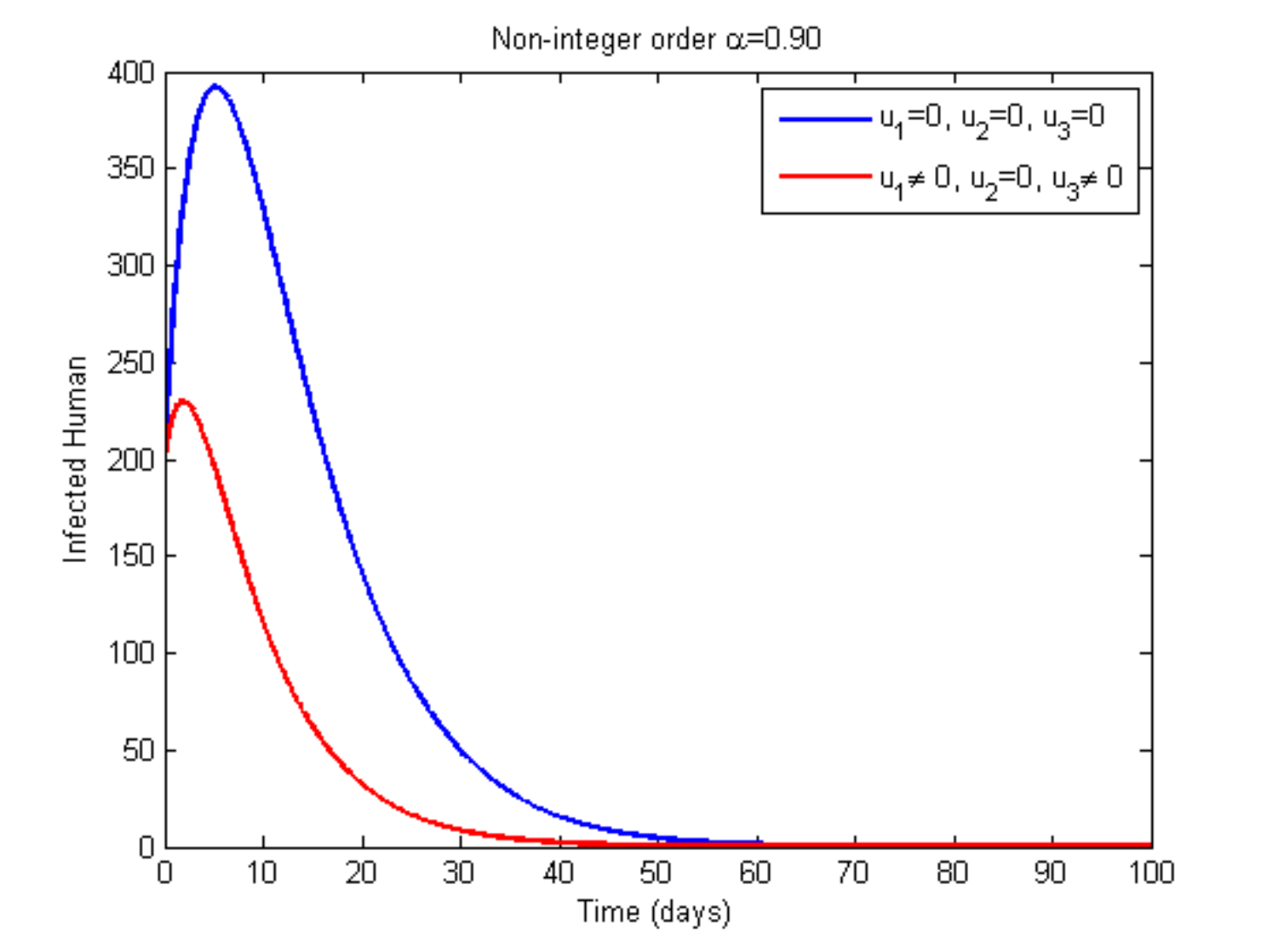}}
\caption{Numerical solutions of infected human host with treated bednets and insecticide spray controls for $\alpha=1, 0.99, 0.95, 0.90$}
\label{fg5}
\end{figure}

\begin{figure}[!ht]
\centering
\subfigure[]{
\includegraphics[scale=0.5]{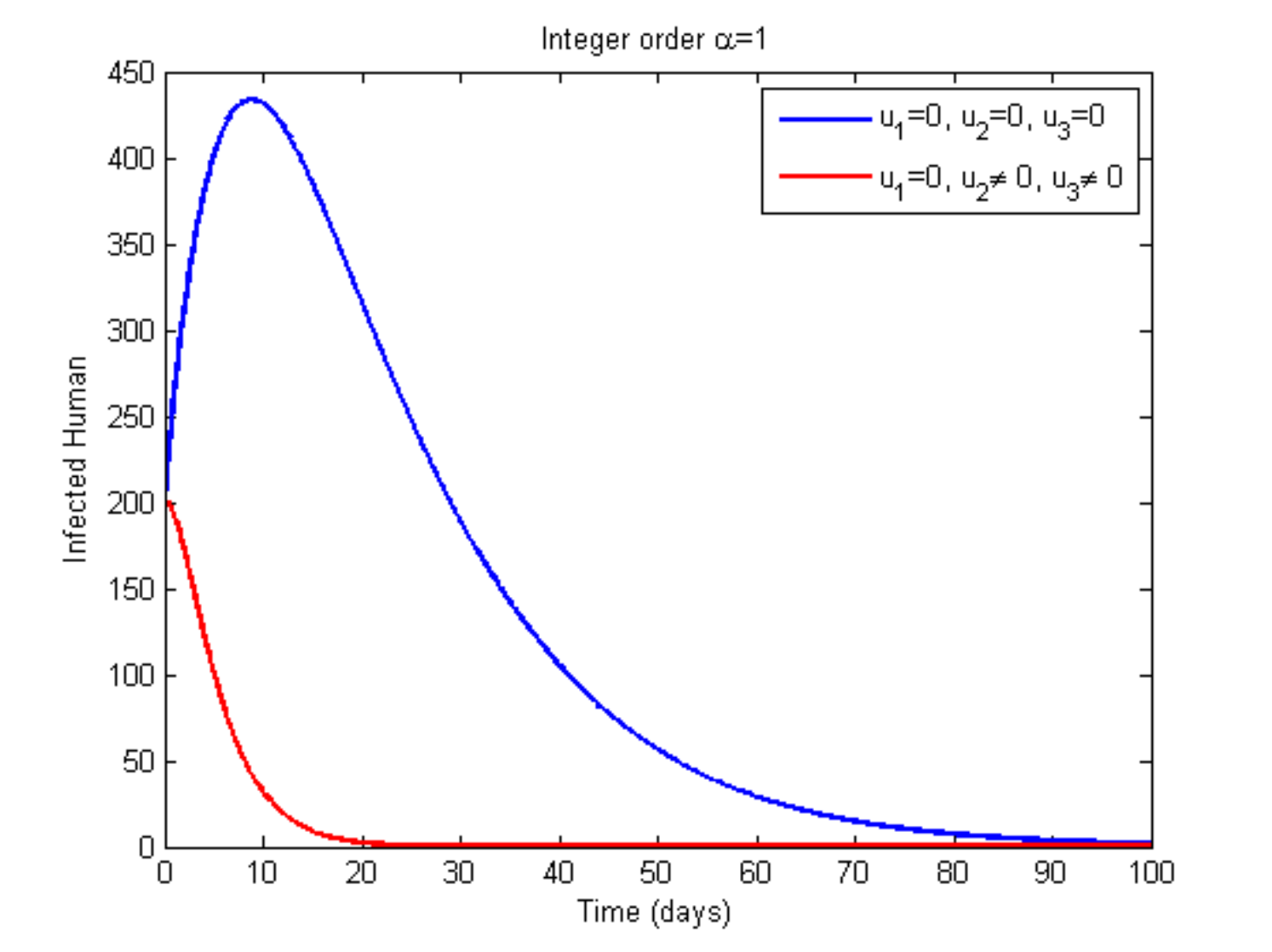}}\hfil
\subfigure[]{
\includegraphics[scale=0.5]{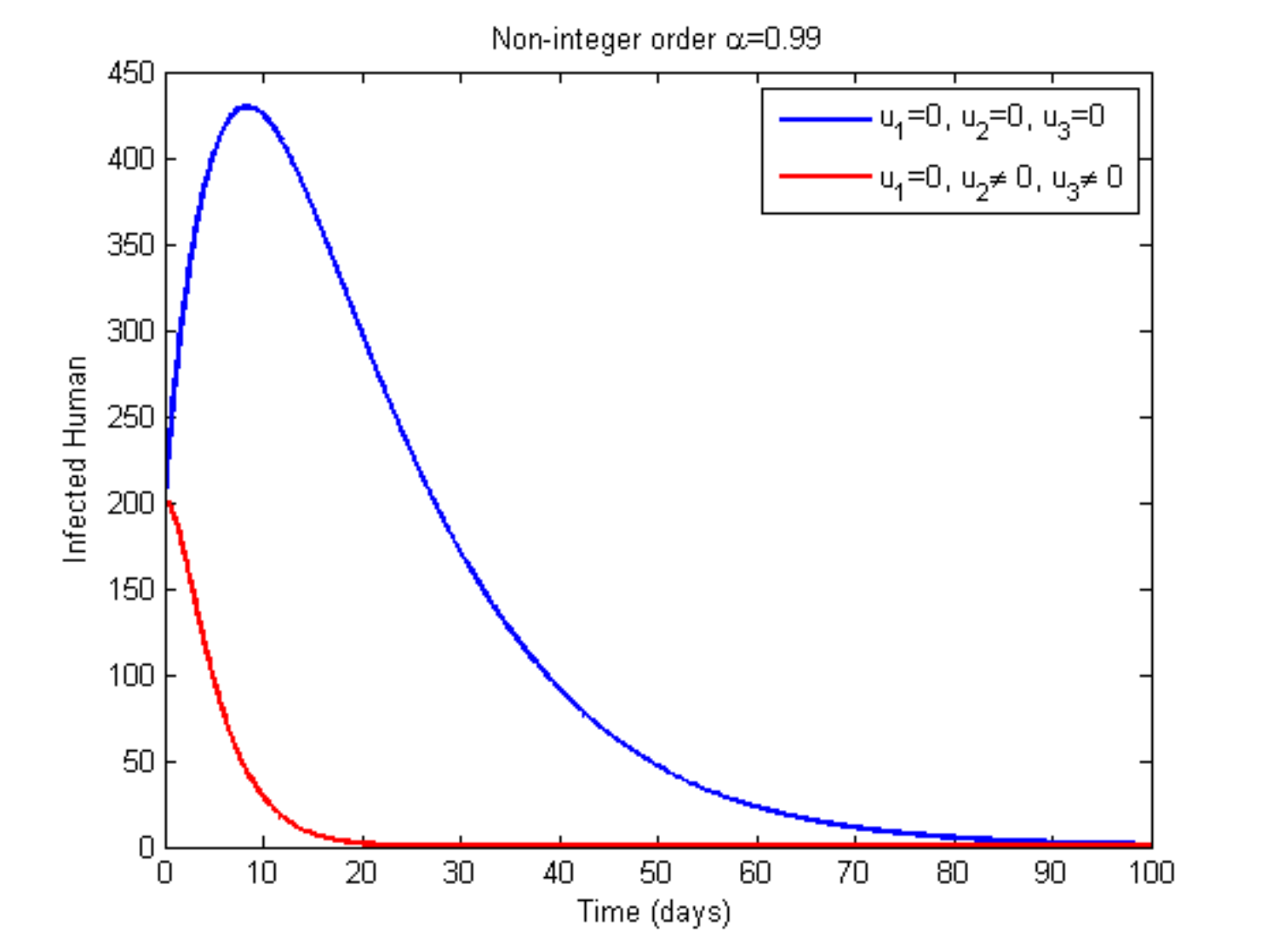}}\hfil
\subfigure[]{
\includegraphics[scale=0.5]{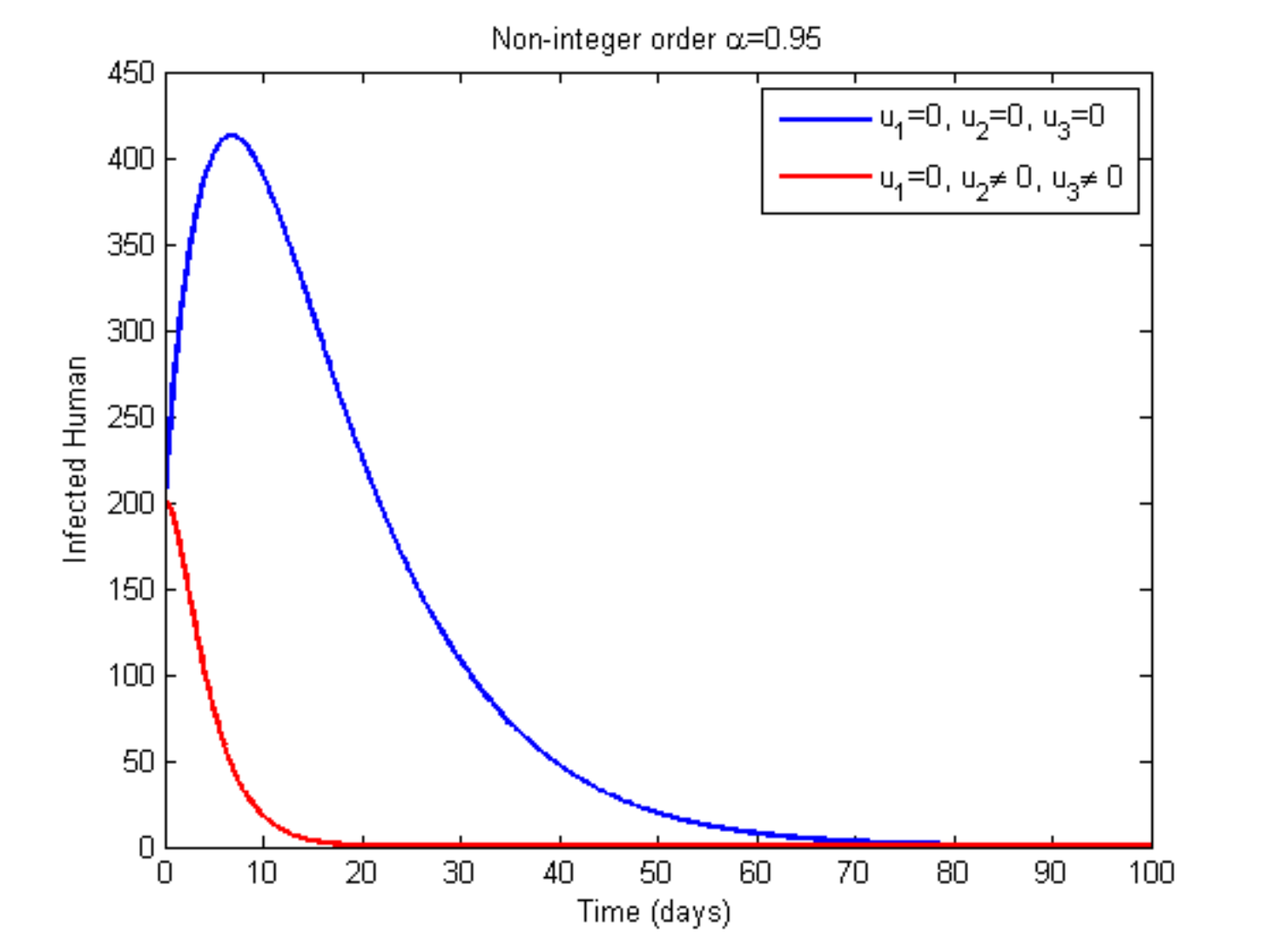}}\hfil
\subfigure[]{\includegraphics[scale=0.5]{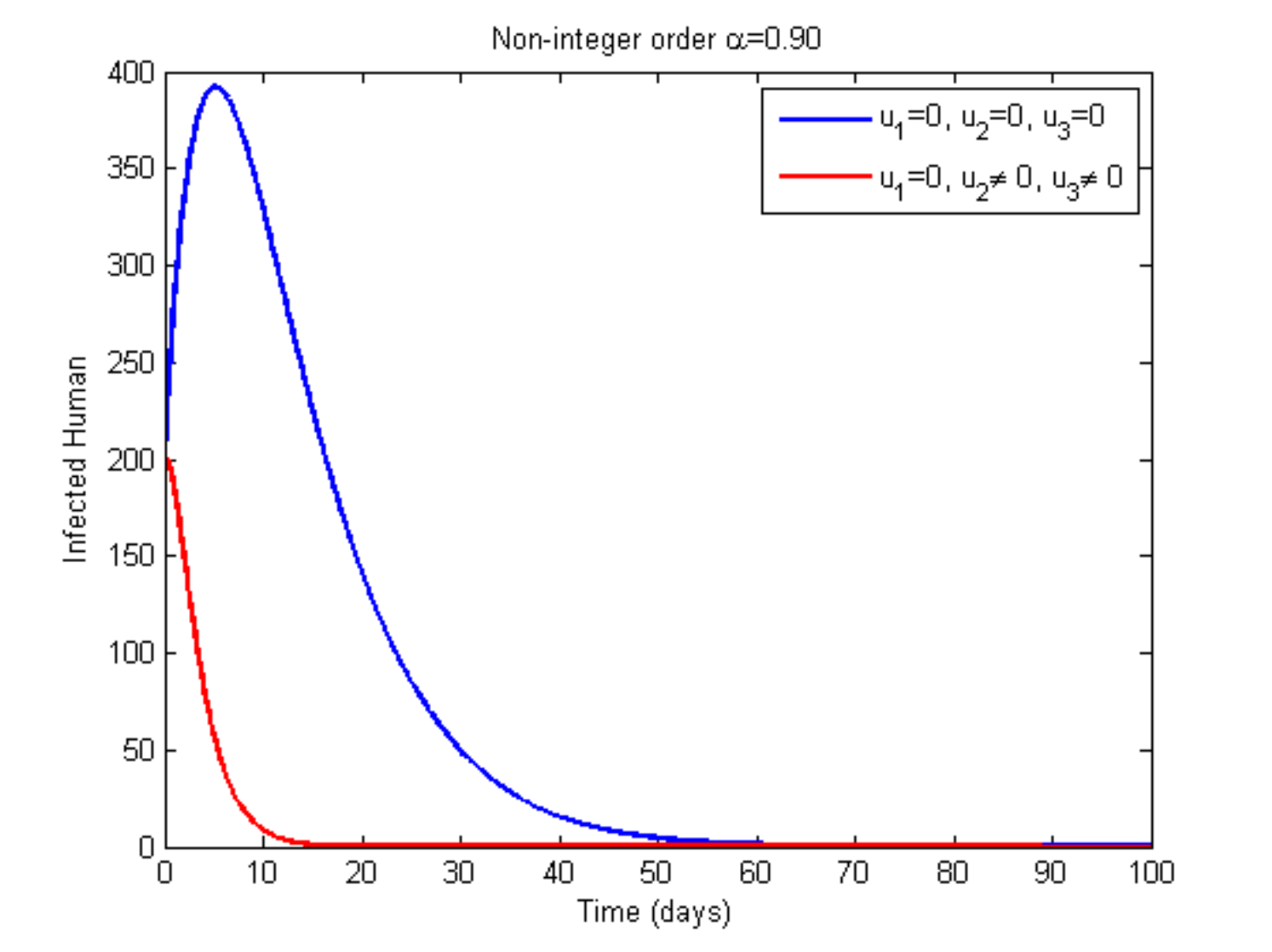}}
\caption{Numerical solutions of infected human host with treatment and insecticide spray controls for $\alpha=1, 0.99, 0.95, 0.90$}
\label{fg6}
\end{figure}

\begin{figure}[!ht]
\centering
\subfigure[]{
\includegraphics[scale=0.51]{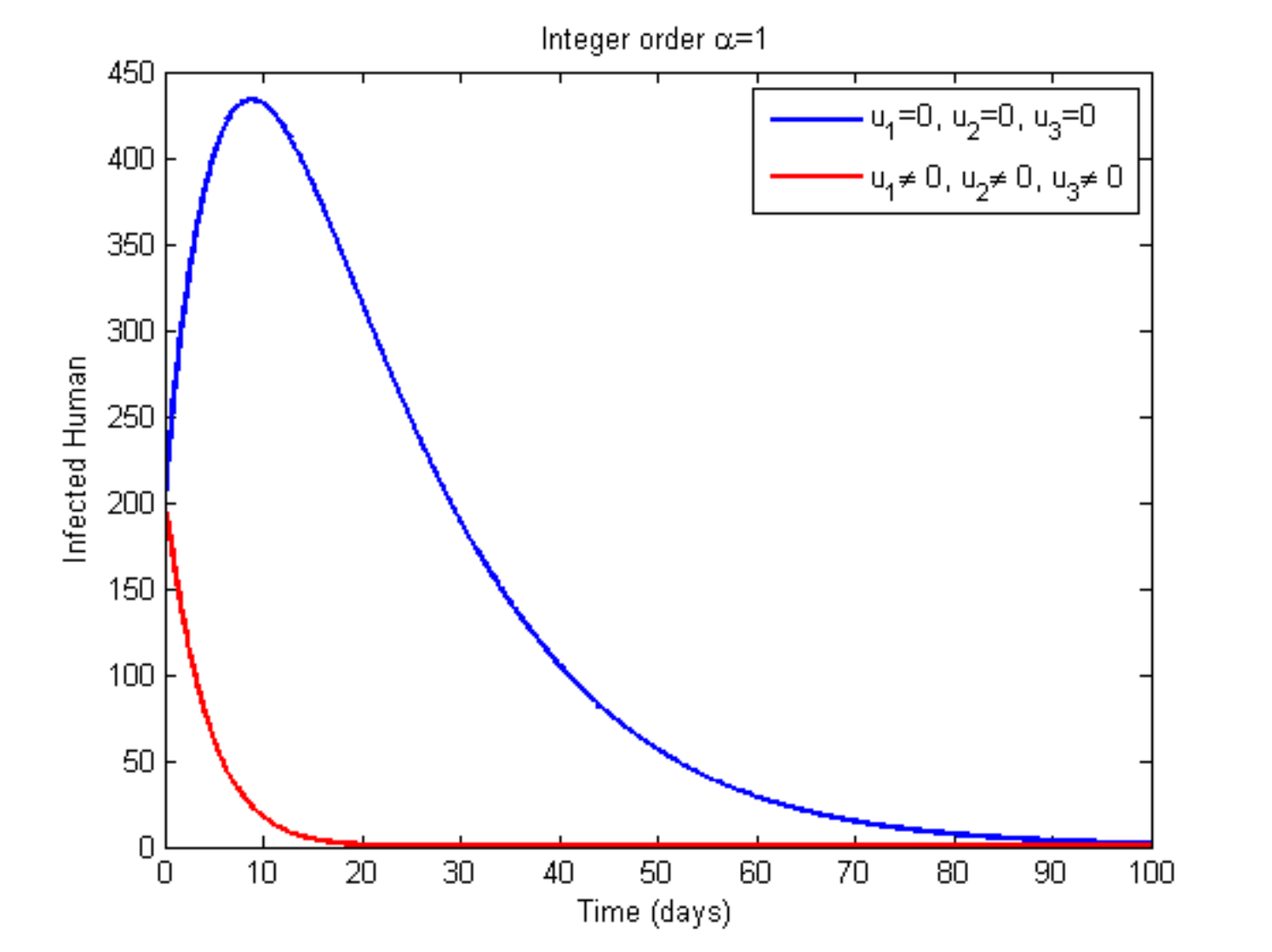}}\hfil
\subfigure[]{
\includegraphics[scale=0.51]{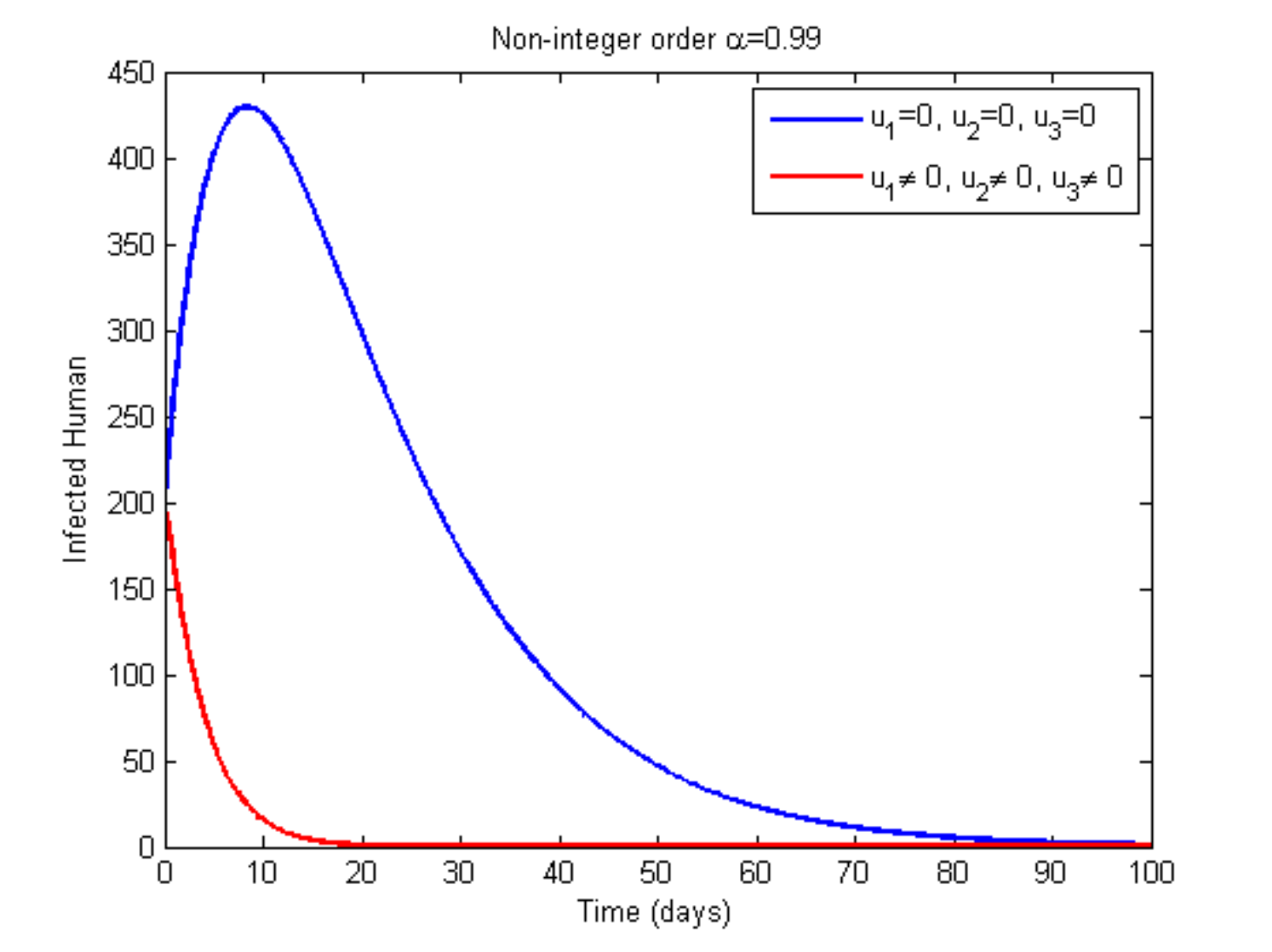}}\hfil
\subfigure[]{
\includegraphics[scale=0.51]{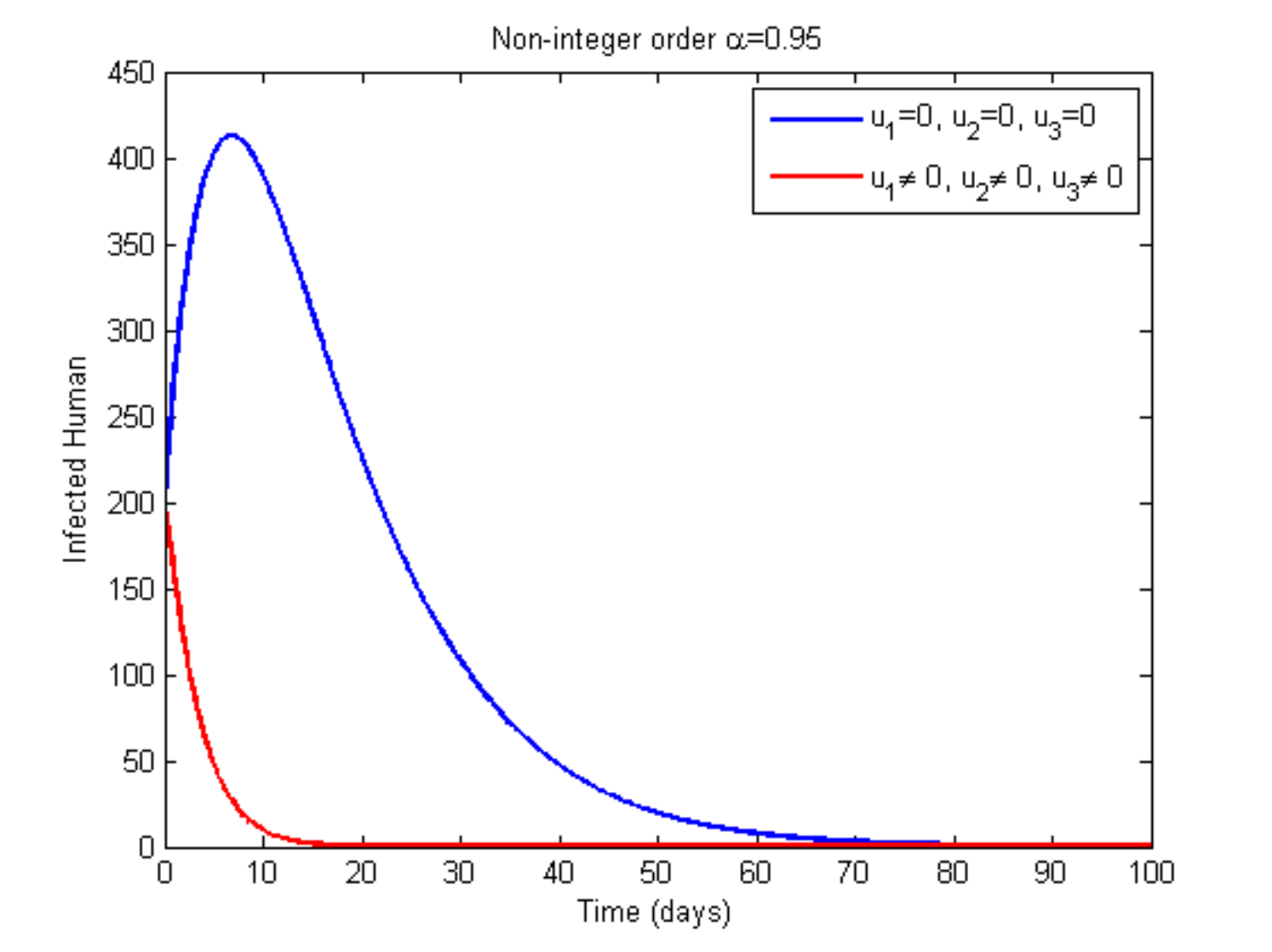}}\hfil
\subfigure[]{\includegraphics[scale=0.51]{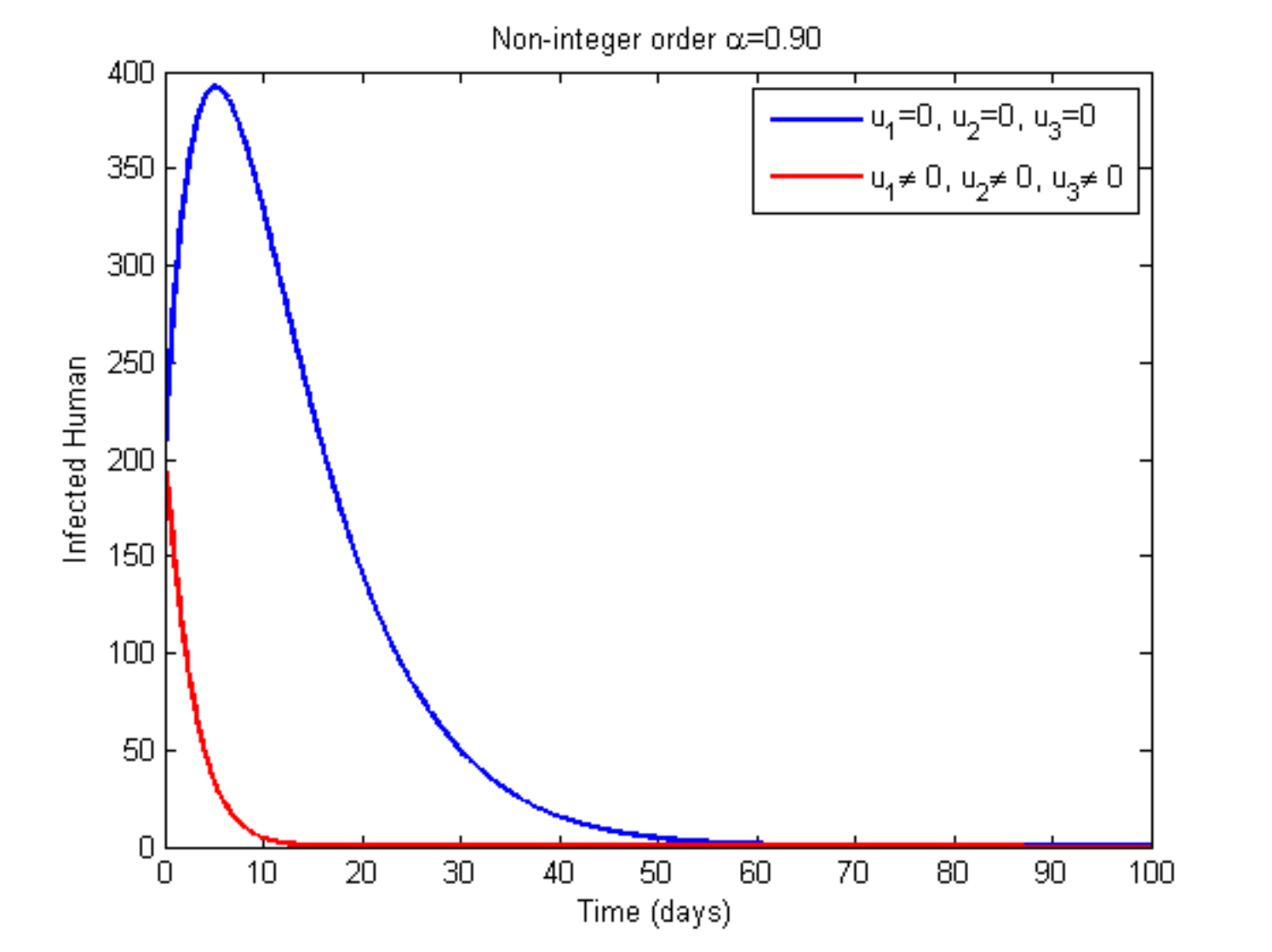}}
\caption{Numerical solutions of infected human host with all the three time dependent optimal controls (treated bednets, treatment and insecticide spray) for $\alpha=1, 0.99, 0.95, 0.90$}
\label{fg7}
\end{figure}


\begin{figure}[!ht]
\centering
\subfigure[]{
\includegraphics[scale=0.5]{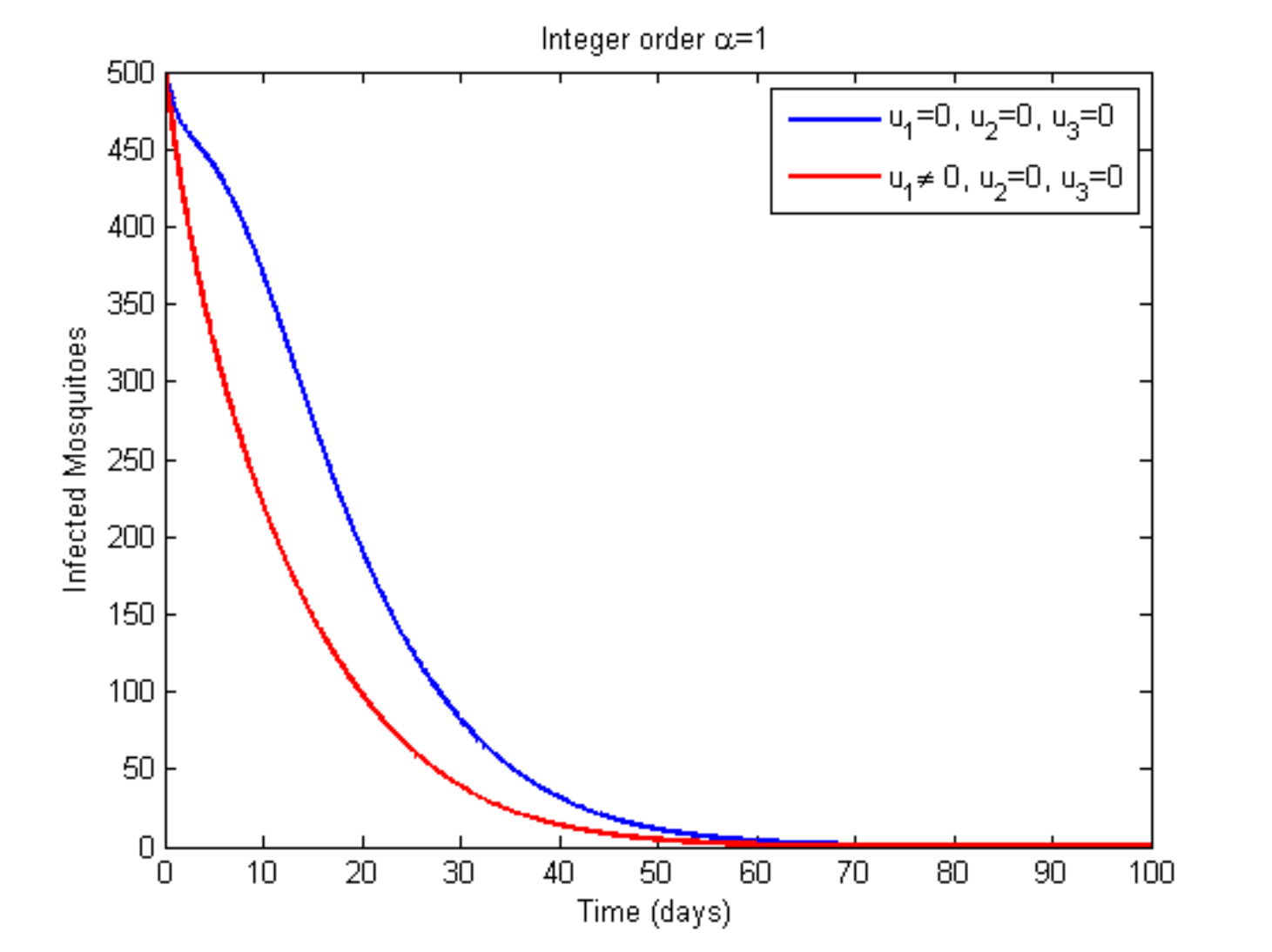}}\hfil
\subfigure[]{
\includegraphics[scale=0.5]{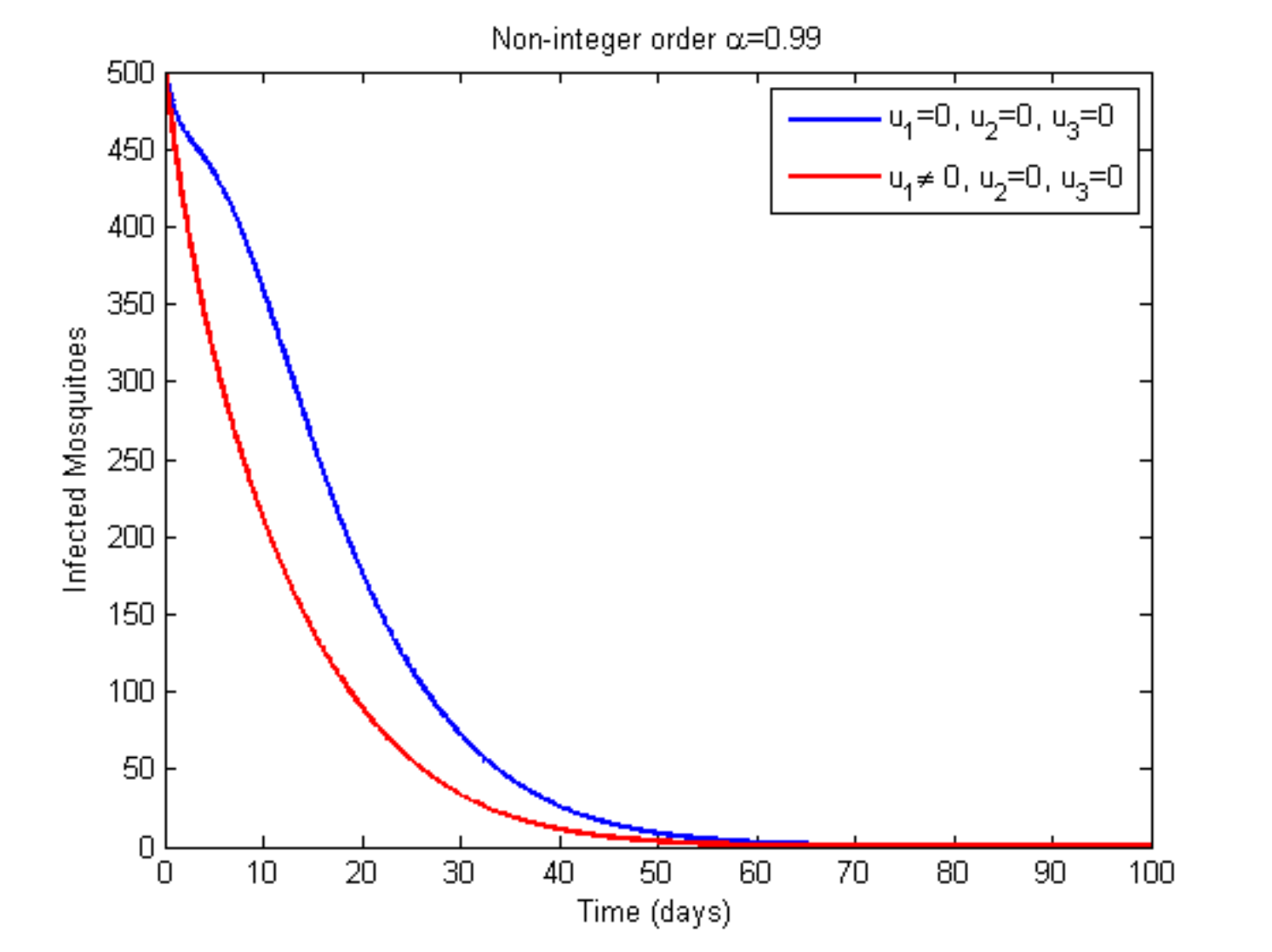}}\hfil
\subfigure[]{
\includegraphics[scale=0.5]{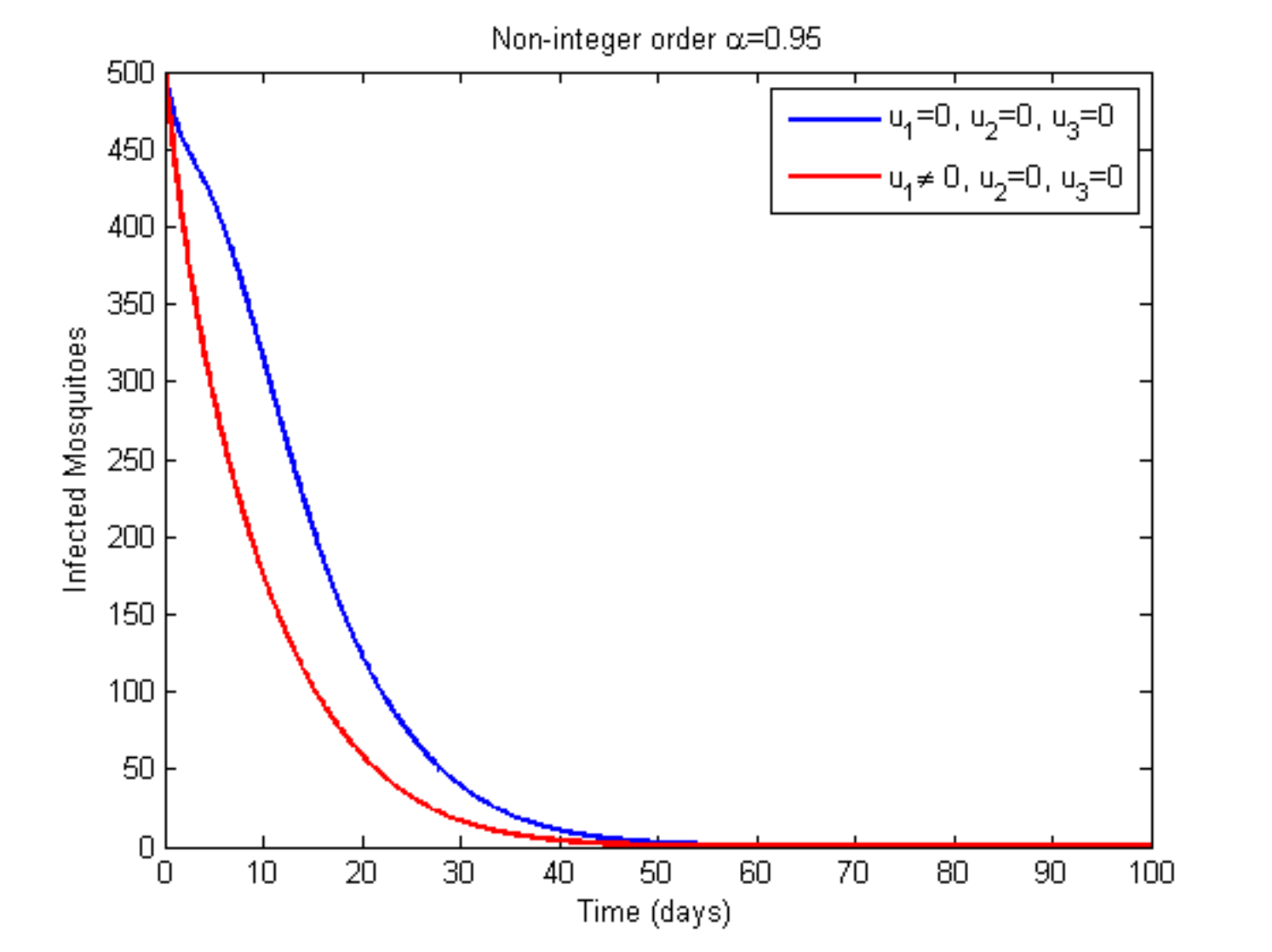}}\hfil
\subfigure[]{\includegraphics[scale=0.5]{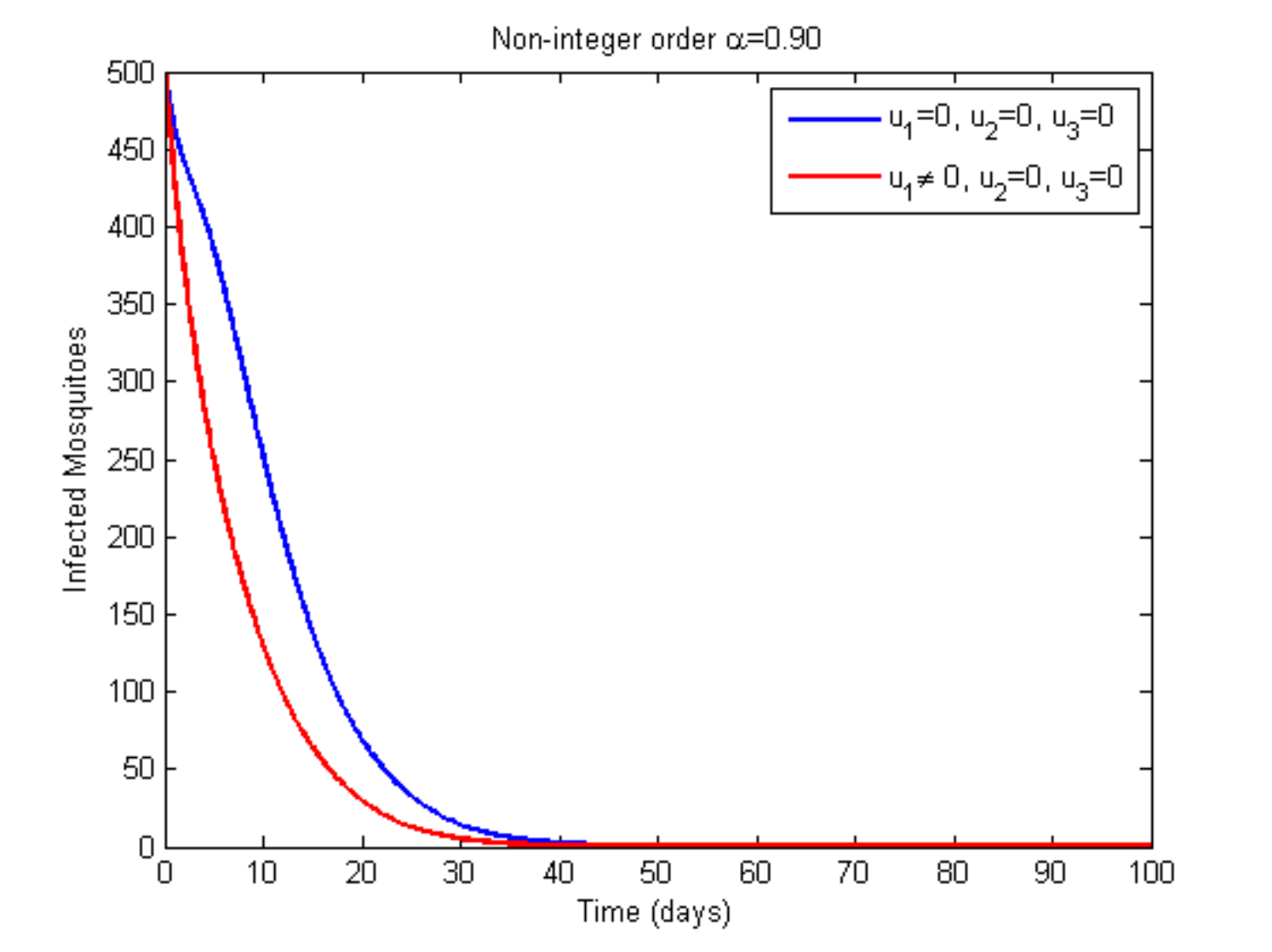}}
\caption{Numerical solutions of infected mosquitoes with treated bednets control for $\alpha=1, 0.99, 0.95, 0.90$}
\label{fg8}
\end{figure}

\begin{figure}[!ht]
\centering
\subfigure[]{
\includegraphics[scale=0.5]{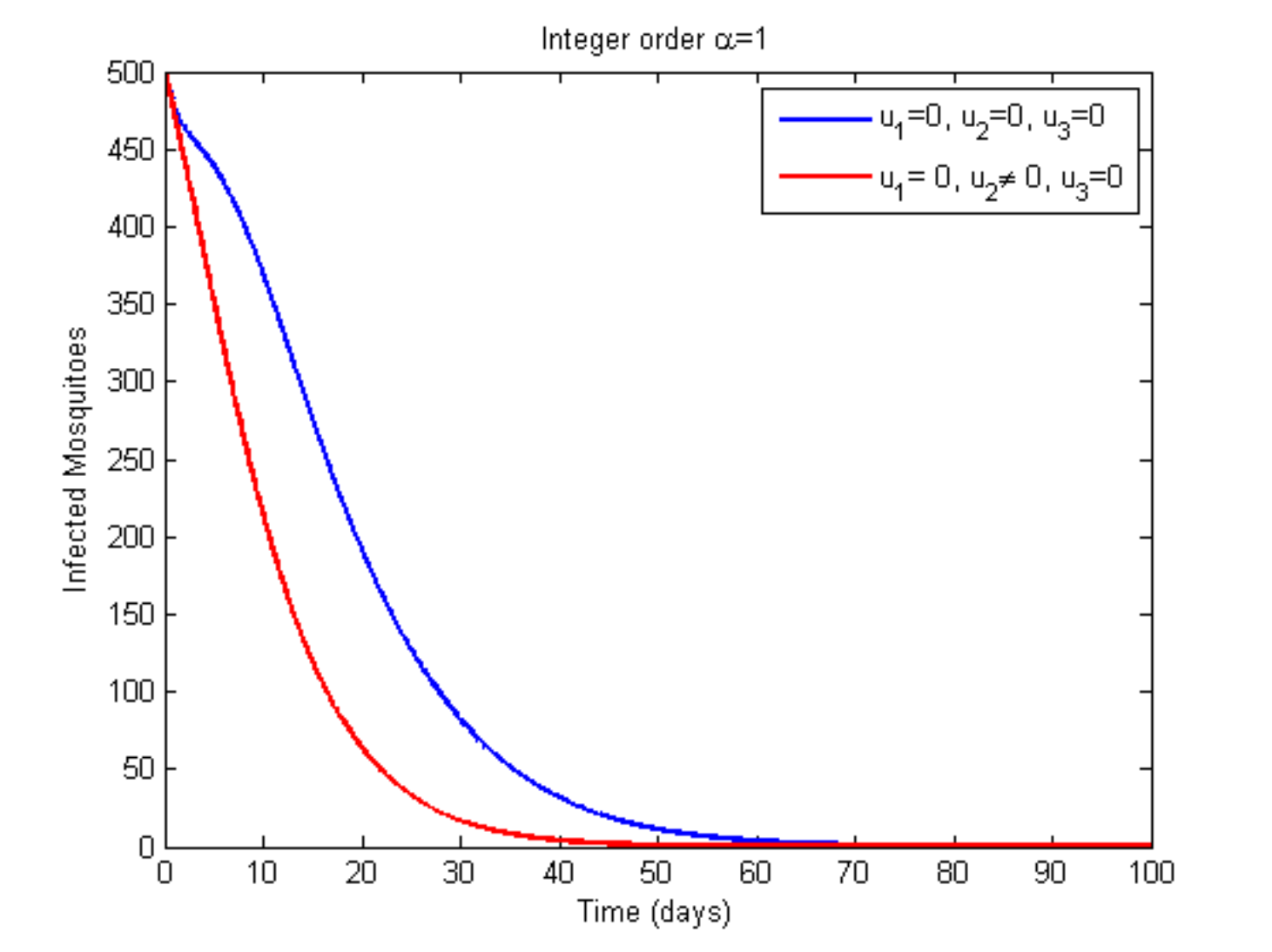}}\hfil
\subfigure[]{
\includegraphics[scale=0.5]{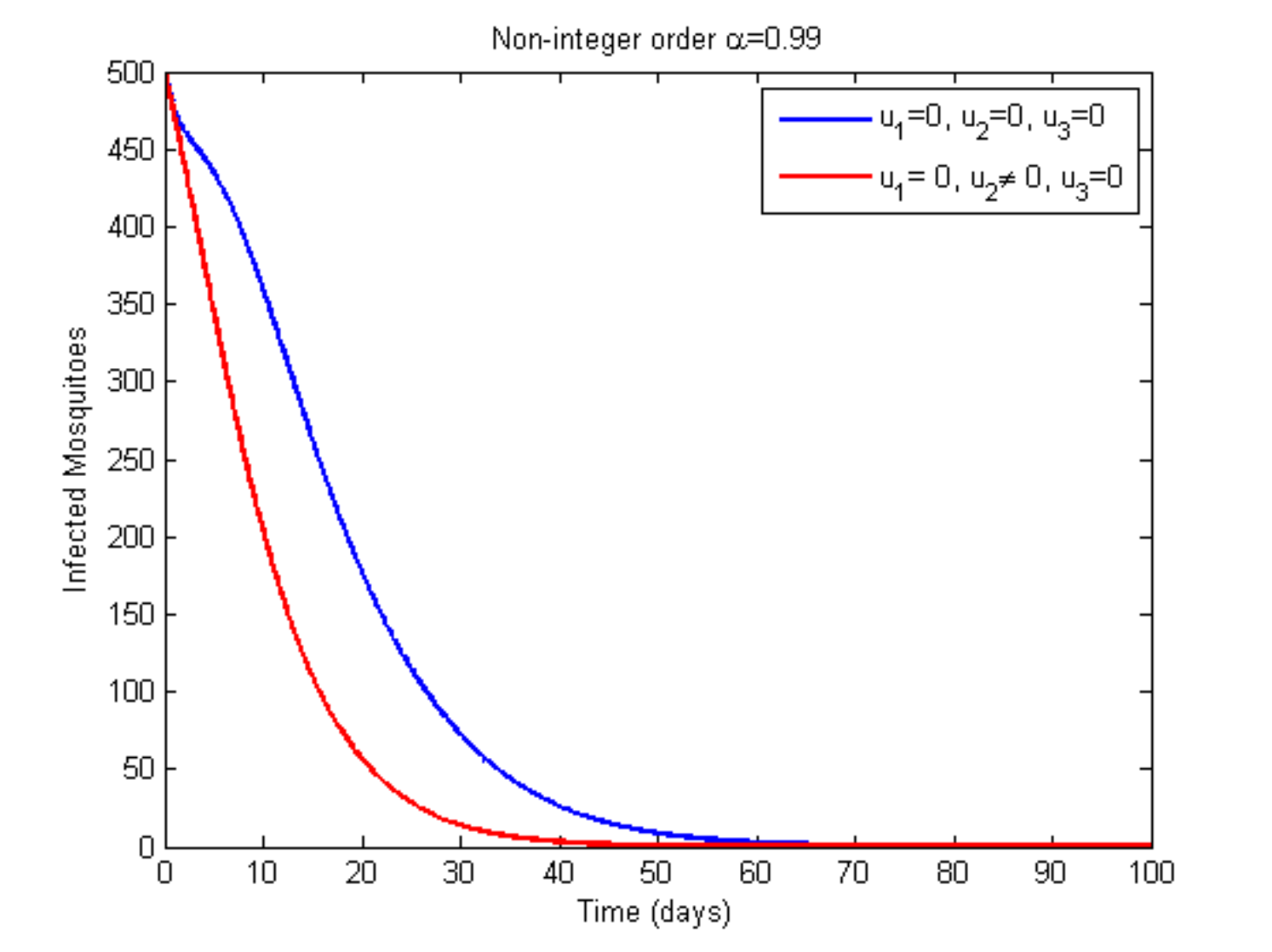}}\hfil
\subfigure[]{
\includegraphics[scale=0.5]{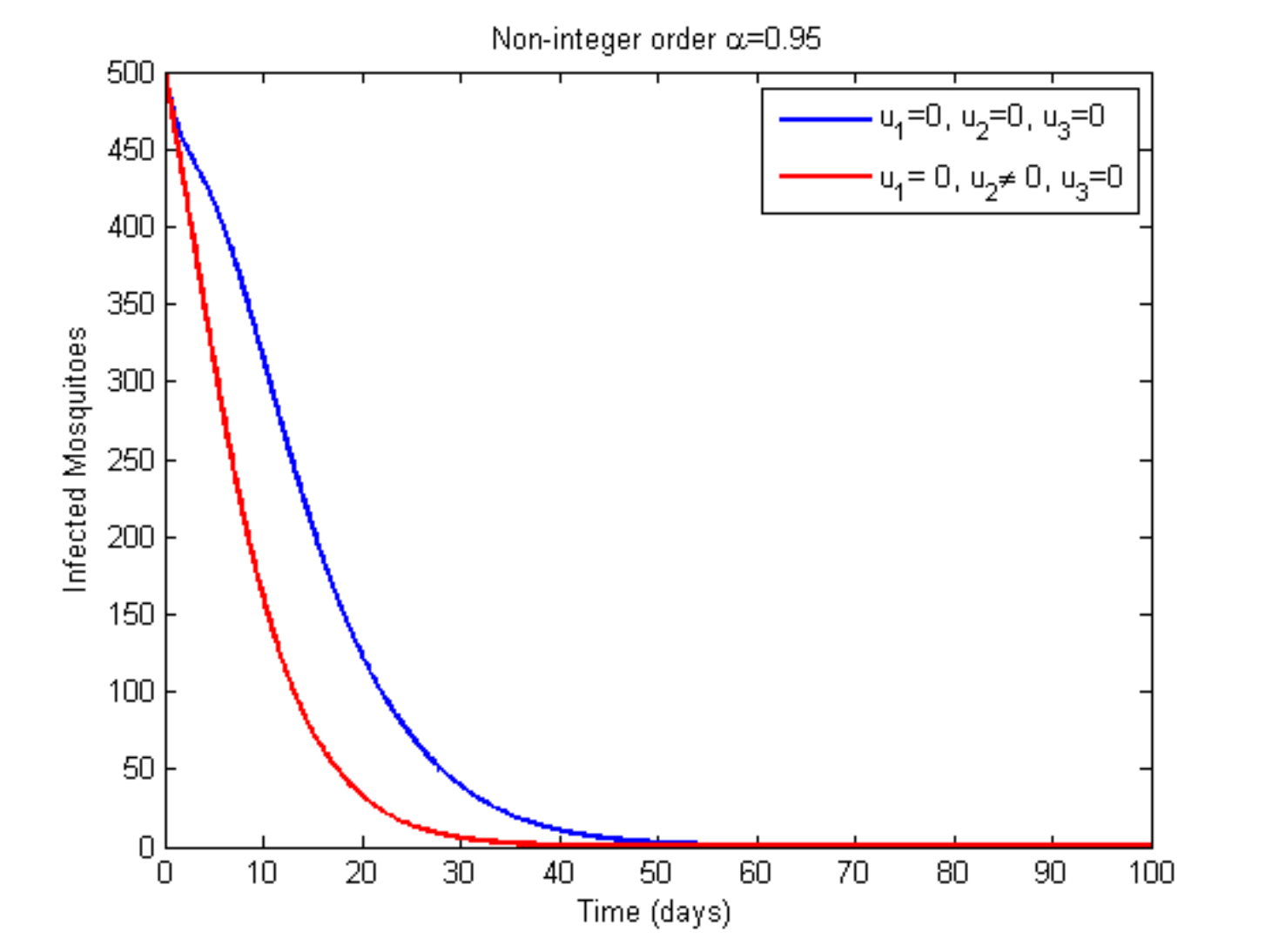}}\hfil
\subfigure[]{\includegraphics[scale=0.5]{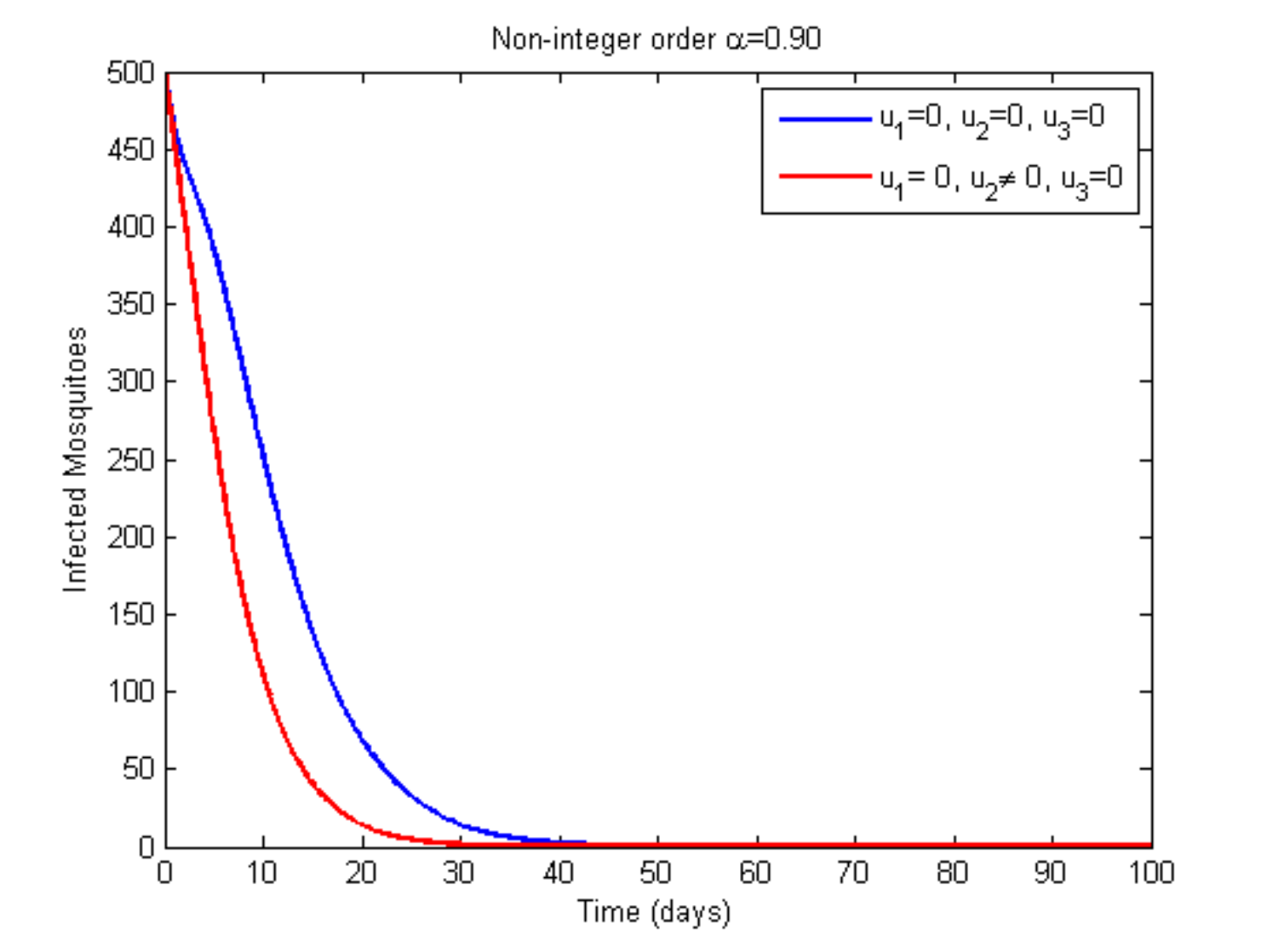}}
\caption{Numerical solutions of infected mosquitoes with treatment control for $\alpha=1, 0.99, 0.95, 0.90$}
\label{fg9}
\end{figure}

\begin{figure}[!ht]
\centering
\subfigure[]{
\includegraphics[scale=0.5]{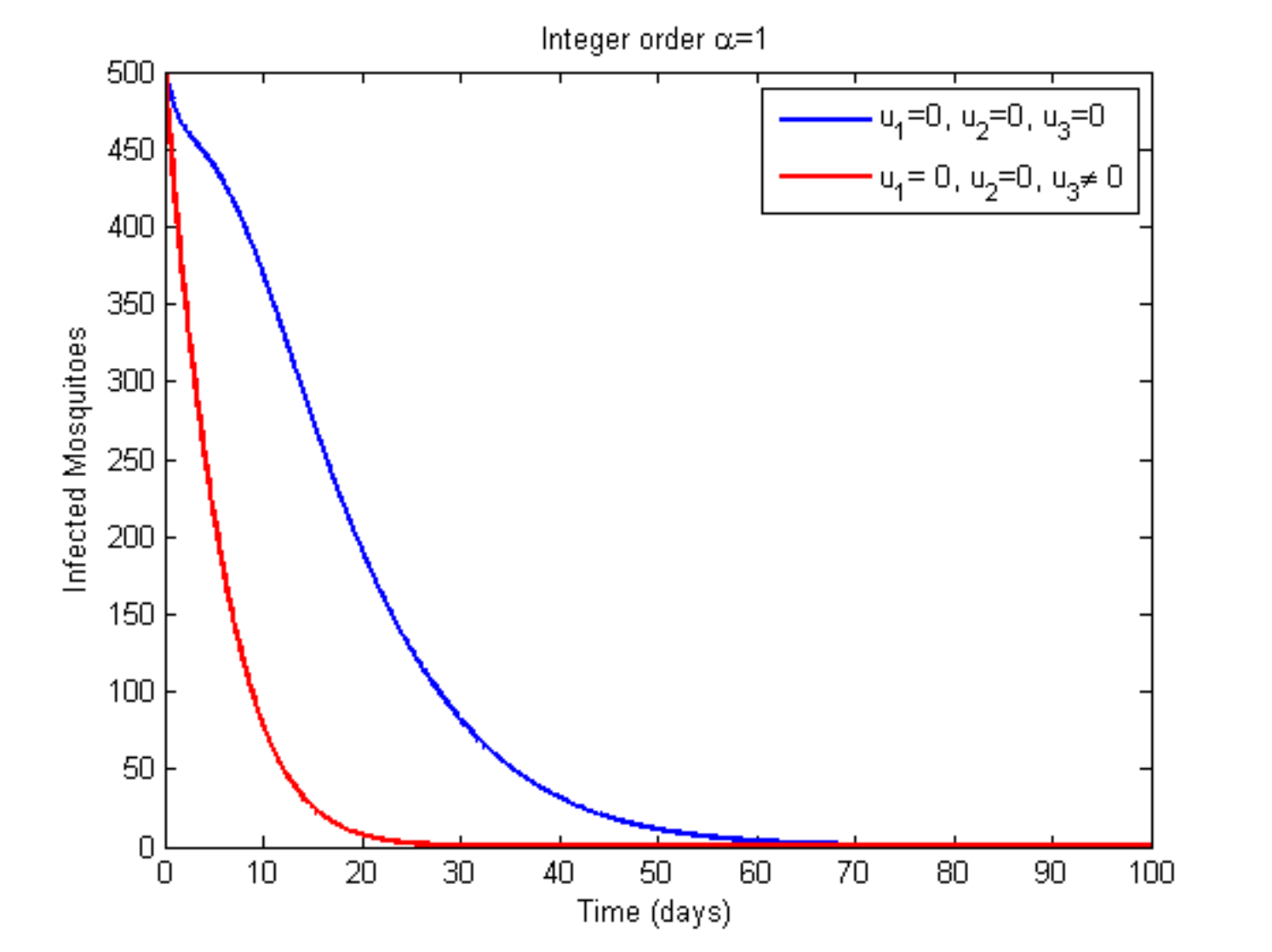}}\hfil
\subfigure[]{
\includegraphics[scale=0.5]{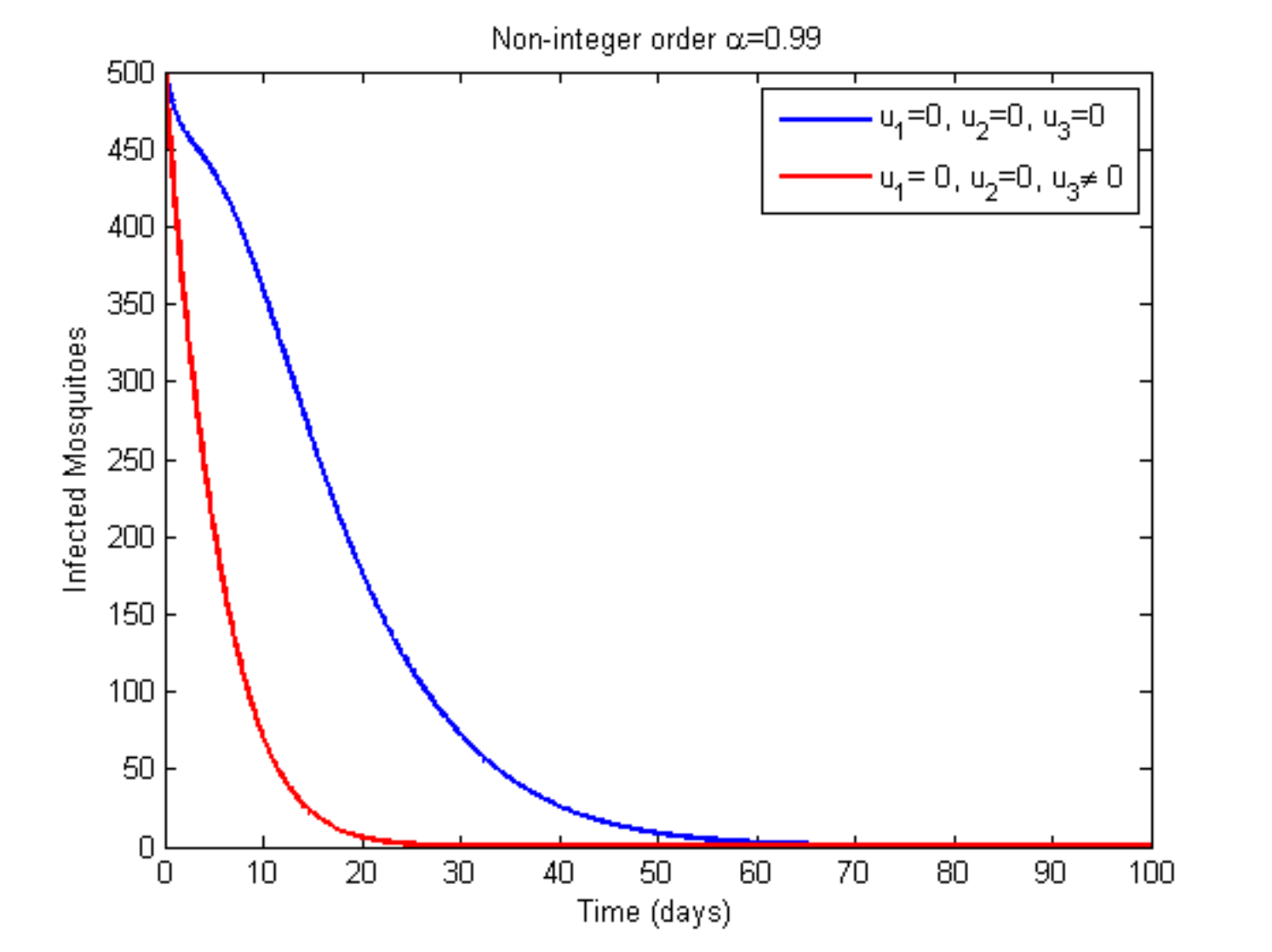}}\hfil
\subfigure[]{
\includegraphics[scale=0.5]{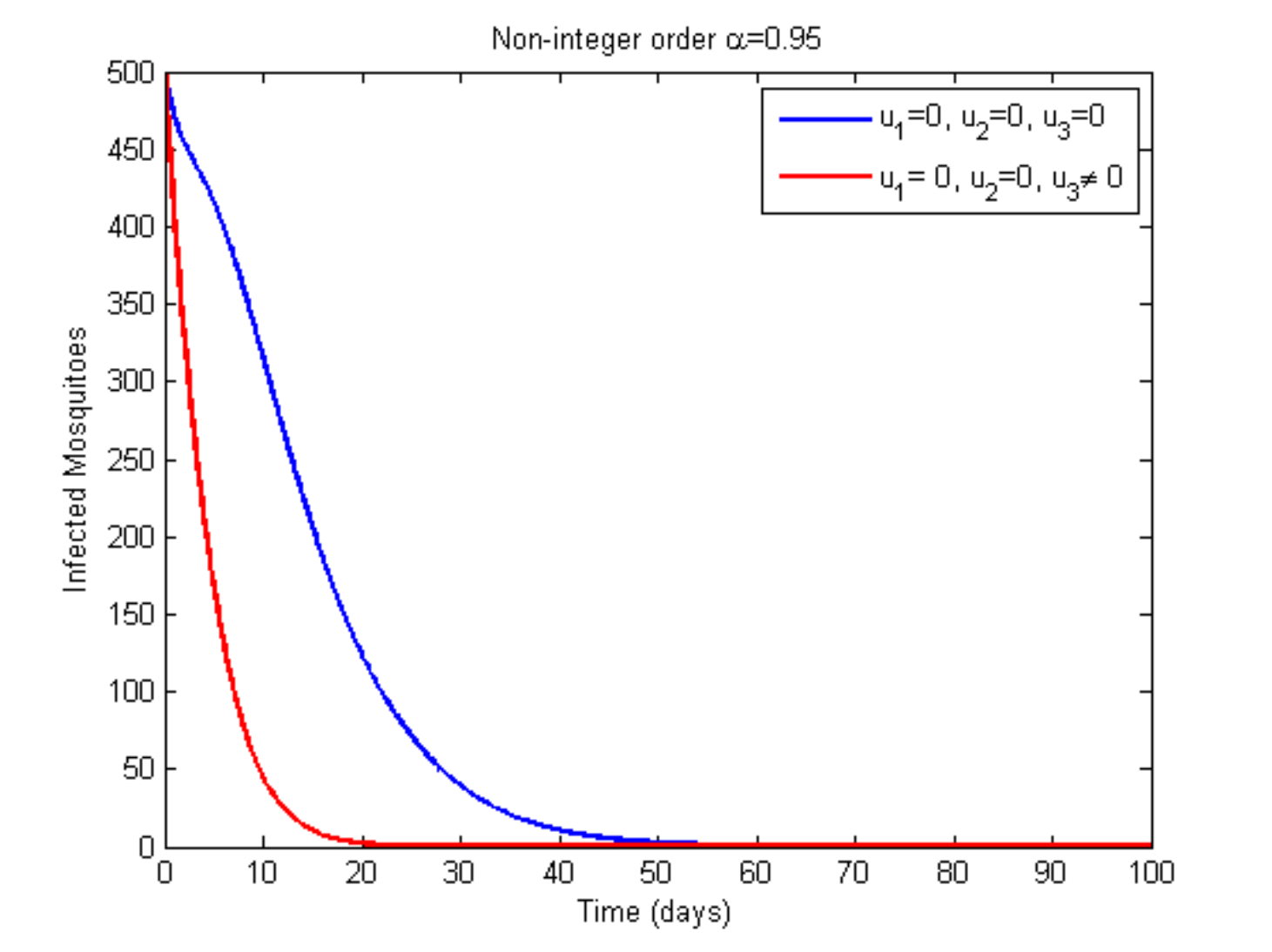}}\hfil
\subfigure[]{\includegraphics[scale=0.5]{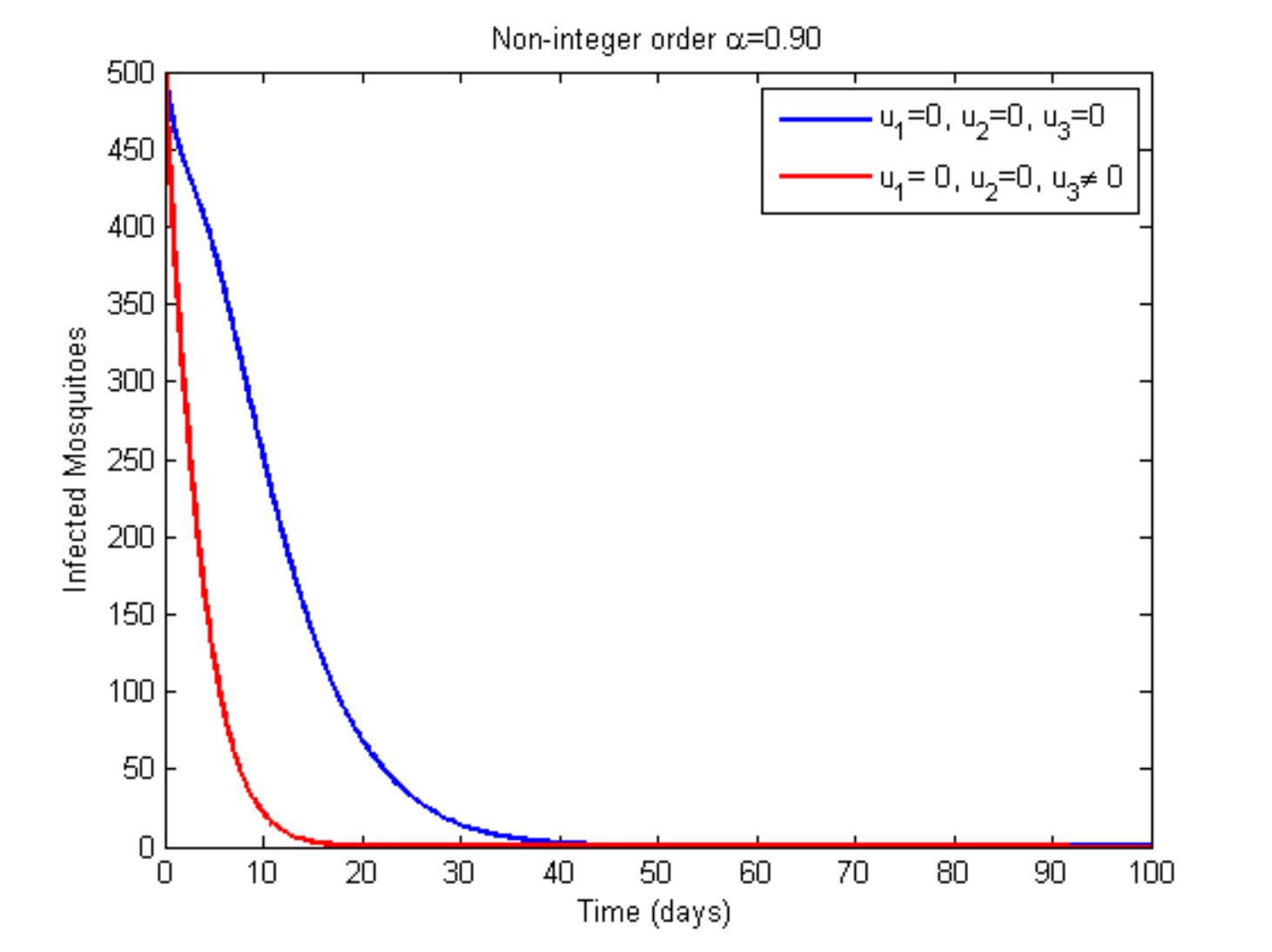}}
\caption{Numerical solutions of infected mosquitoes with insecticide spray control for $\alpha=1, 0.99, 0.95, 0.90$}
\label{fg10}
\end{figure}

\begin{figure}[!ht]
\centering
\subfigure[]{
\includegraphics[scale=0.5]{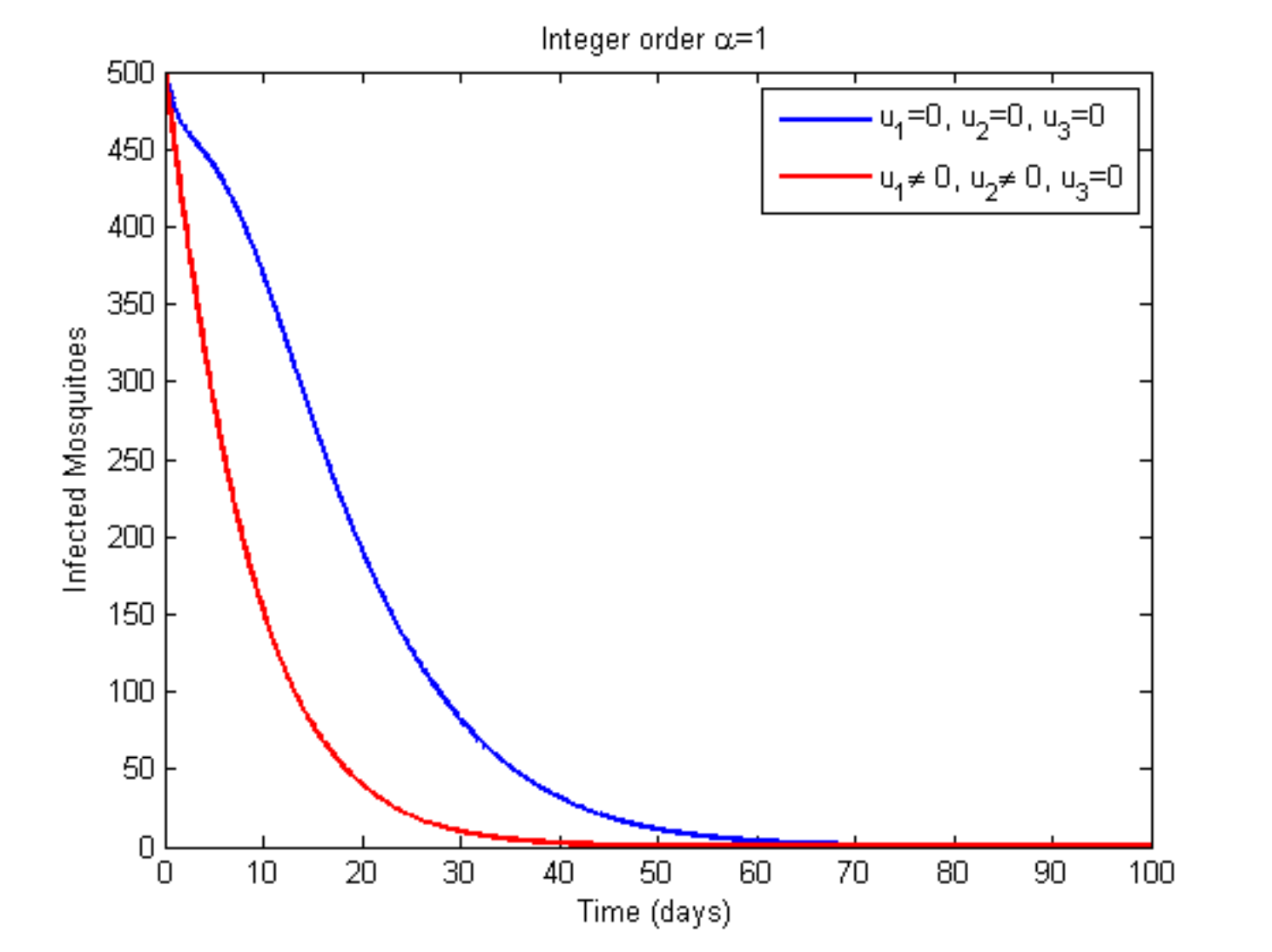}}\hfil
\subfigure[]{
\includegraphics[scale=0.5]{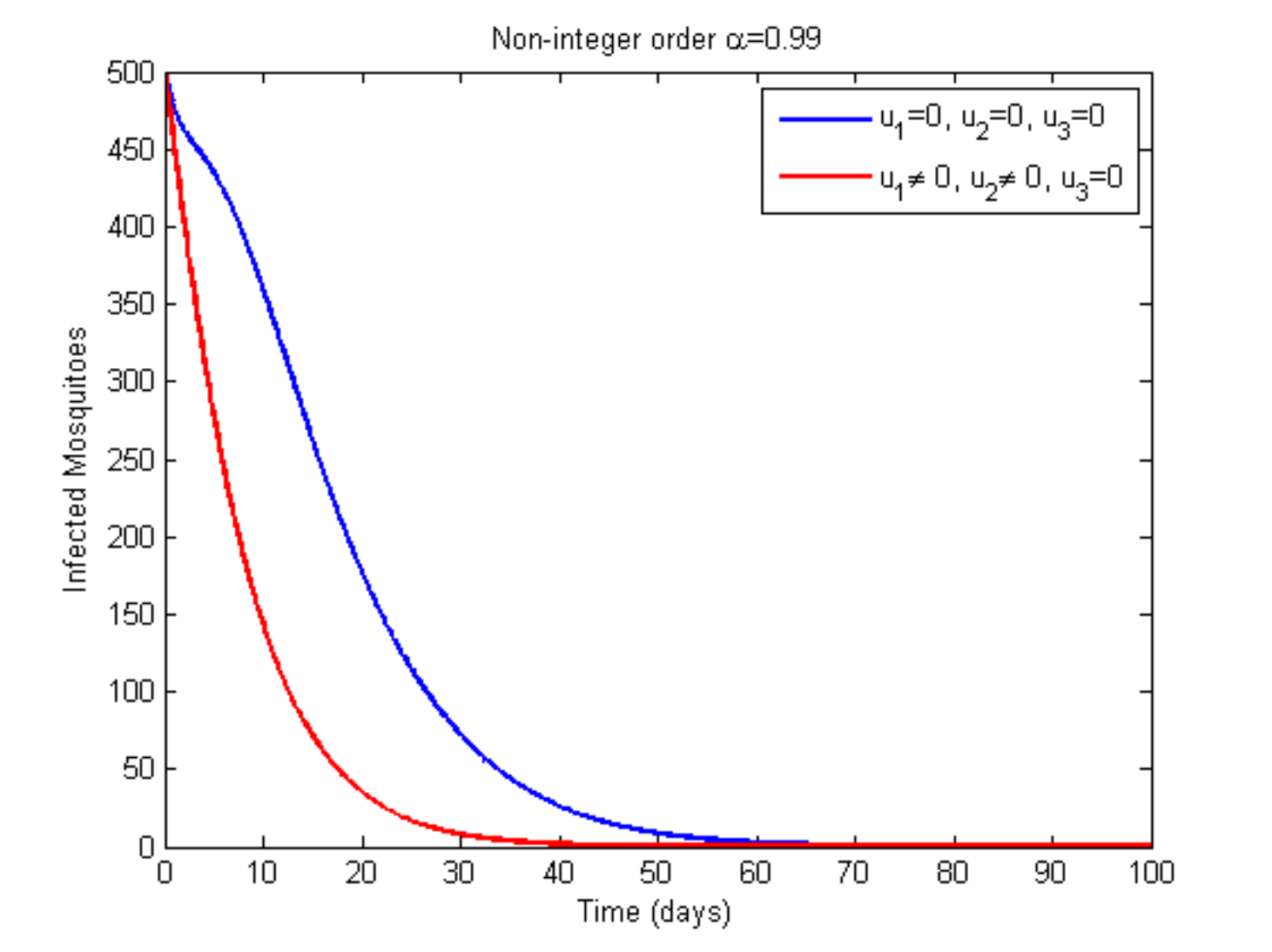}}\hfil
\subfigure[]{
\includegraphics[scale=0.5]{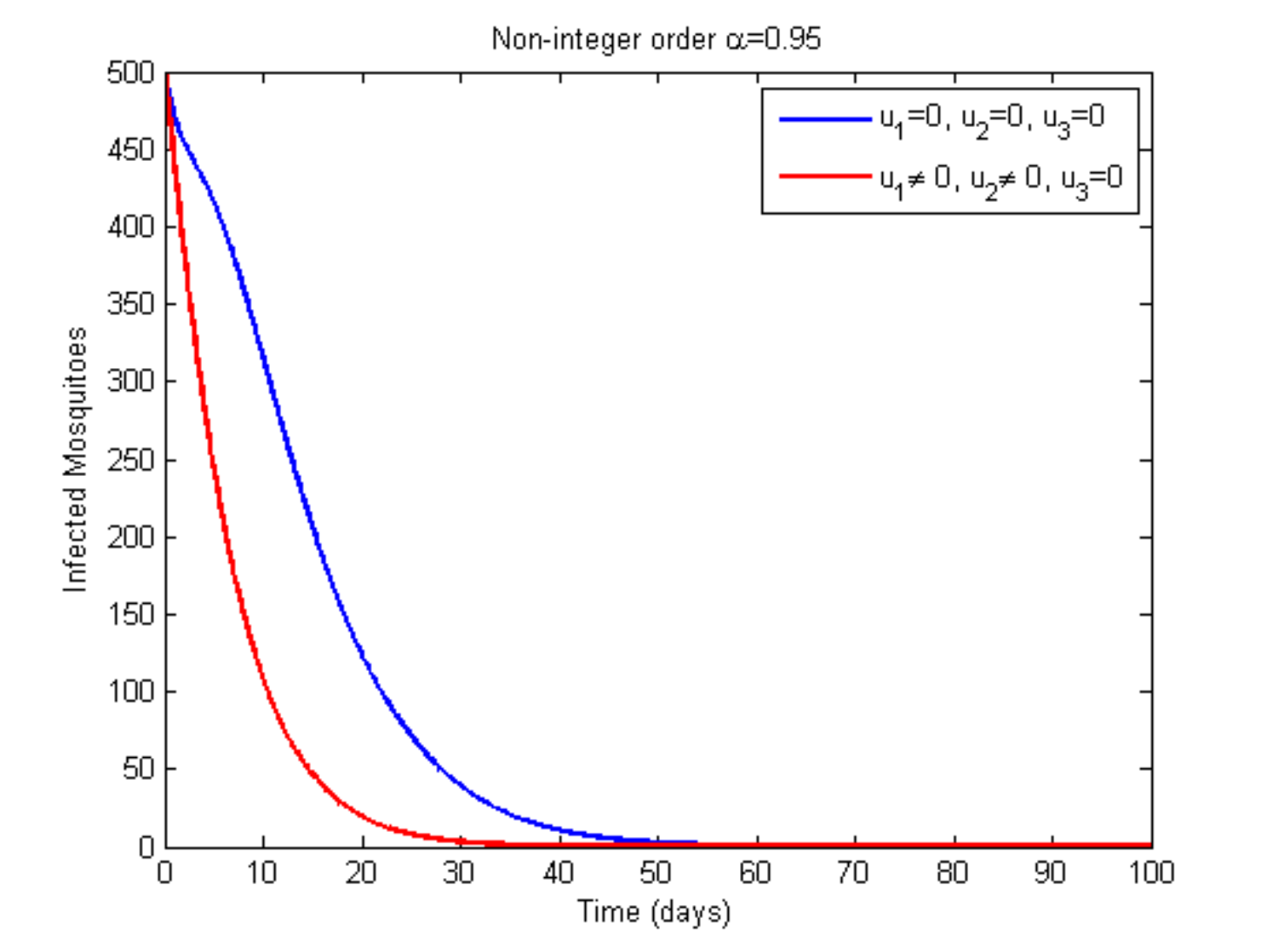}}\hfil
\subfigure[]{\includegraphics[scale=0.5]{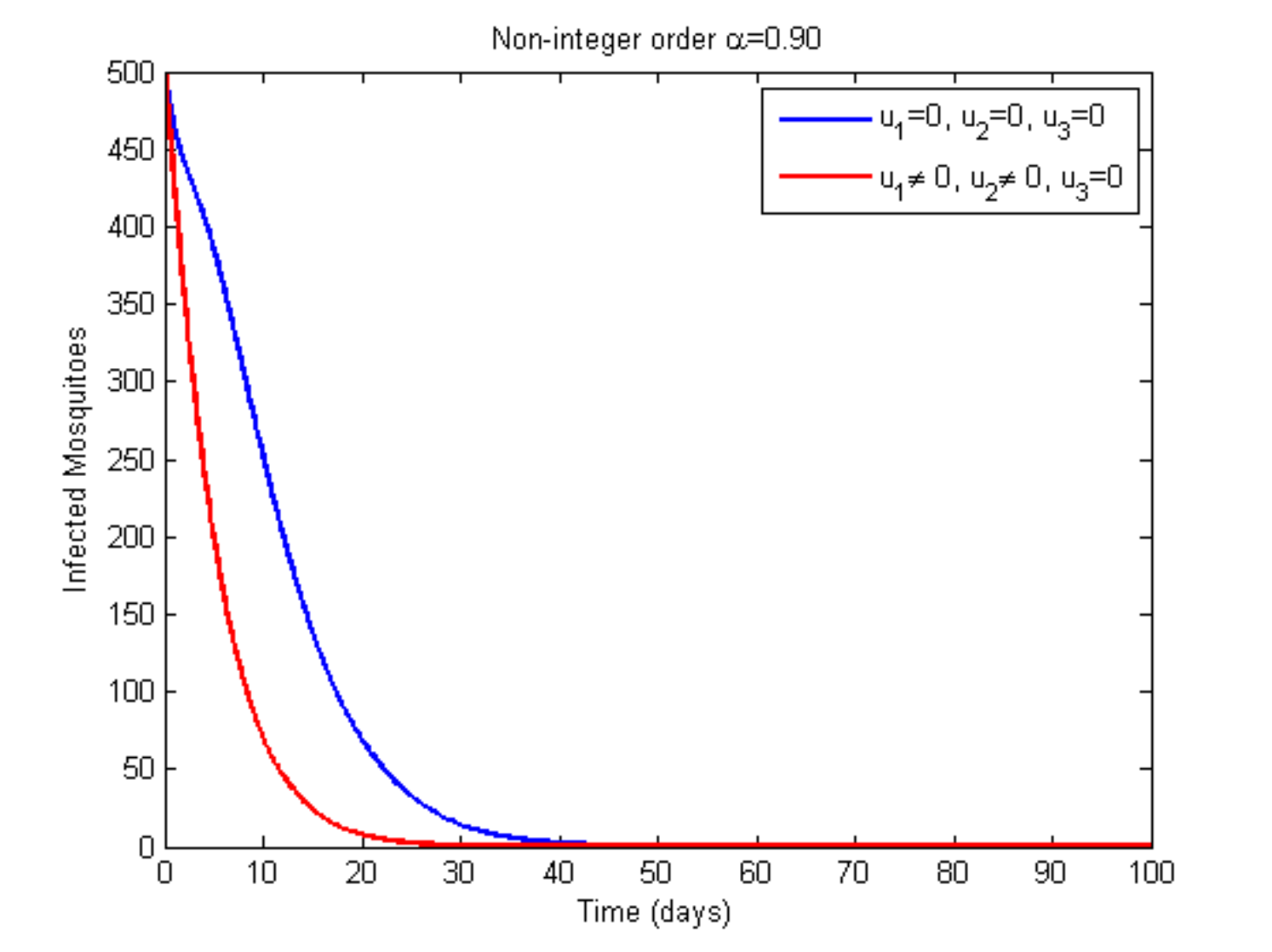}}
\caption{Numerical solutions of infected mosquitoes with treated bednets and treatment controls for $\alpha=1, 0.99, 0.95, 0.90$}
\label{fg11}
\end{figure}

\begin{figure}[!ht]
\centering
\subfigure[]{
\includegraphics[scale=0.5]{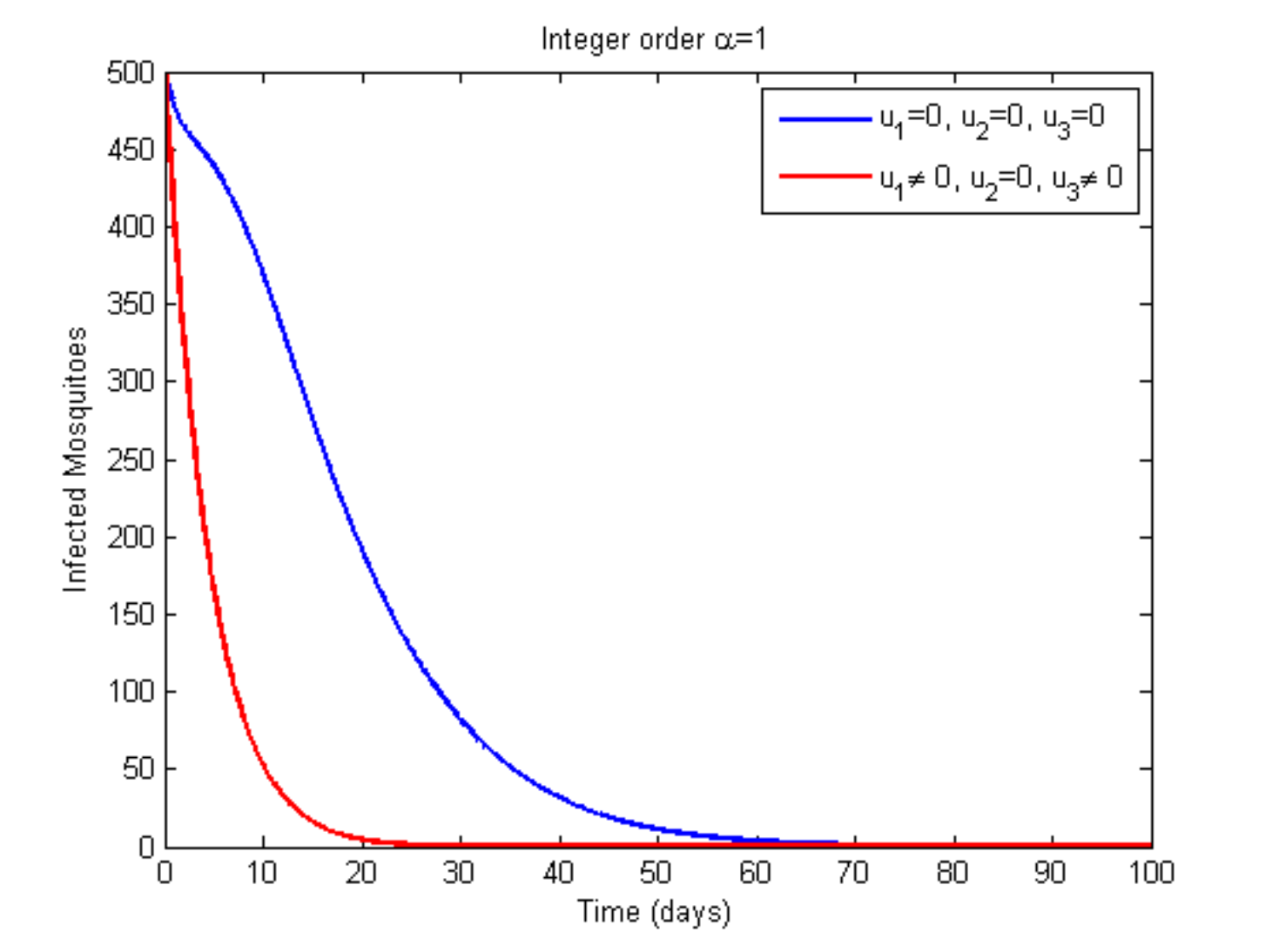}}\hfil
\subfigure[]{
\includegraphics[scale=0.5]{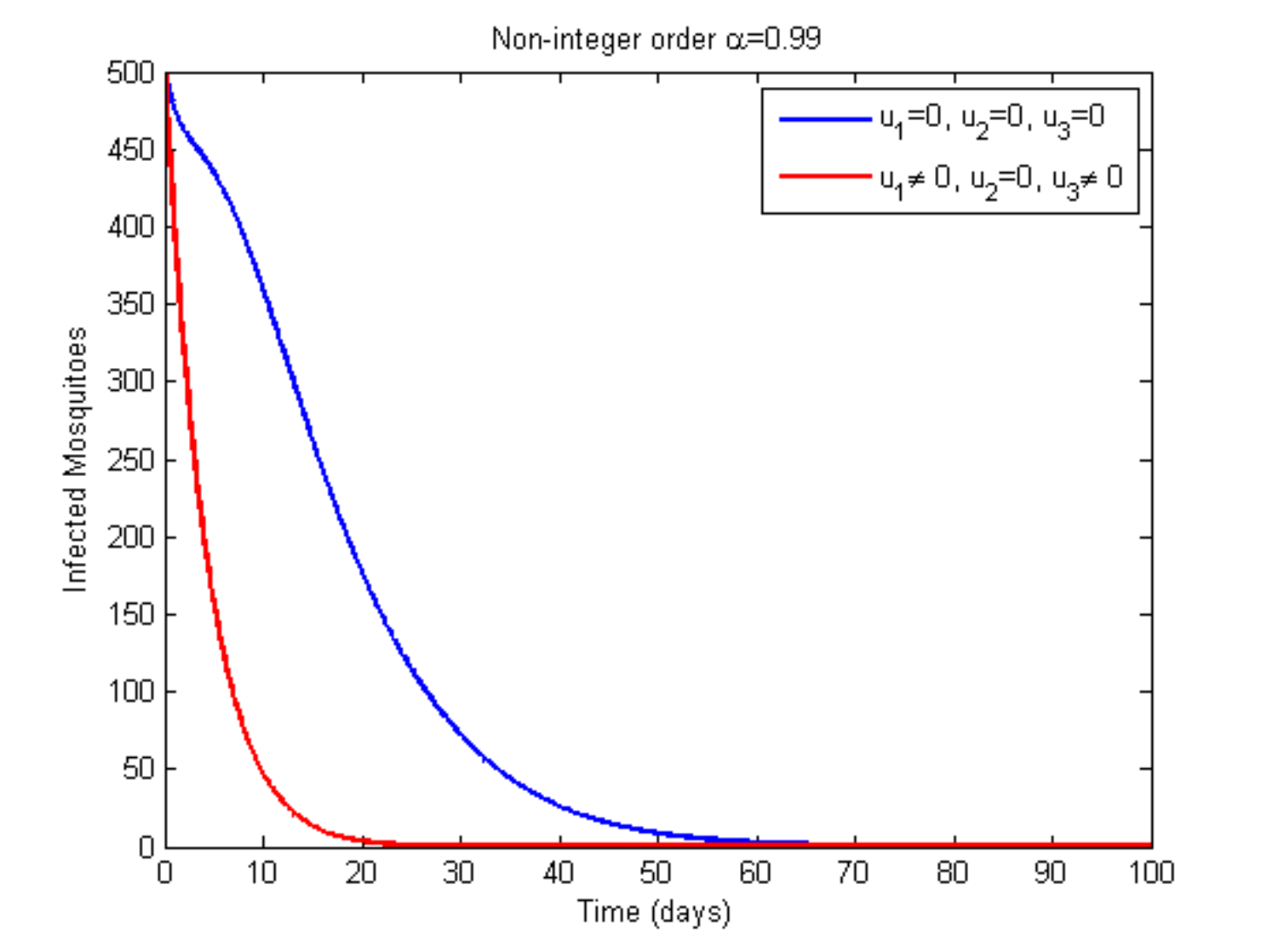}}\hfil
\subfigure[]{
\includegraphics[scale=0.5]{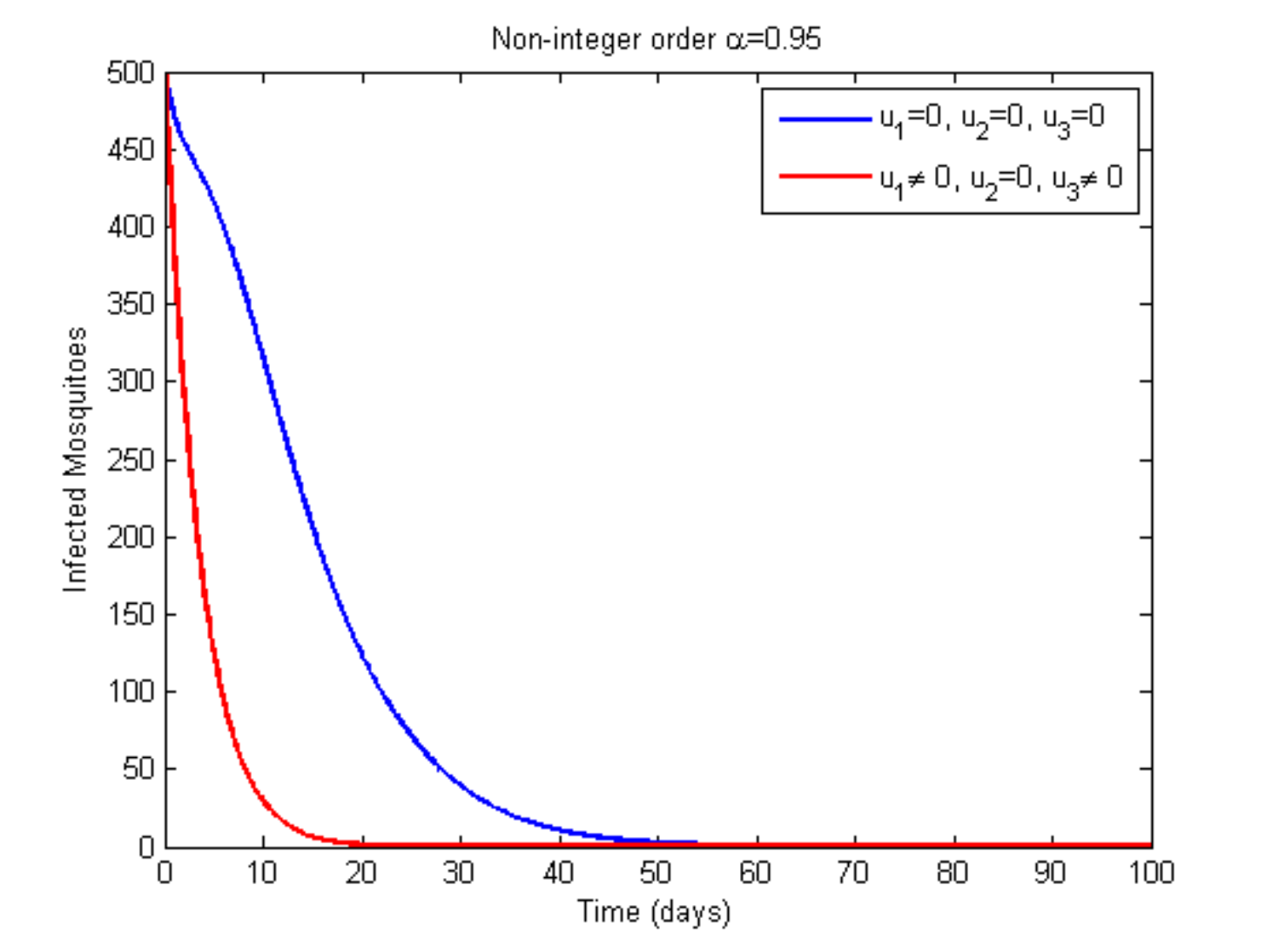}}\hfil
\subfigure[]{\includegraphics[scale=0.5]{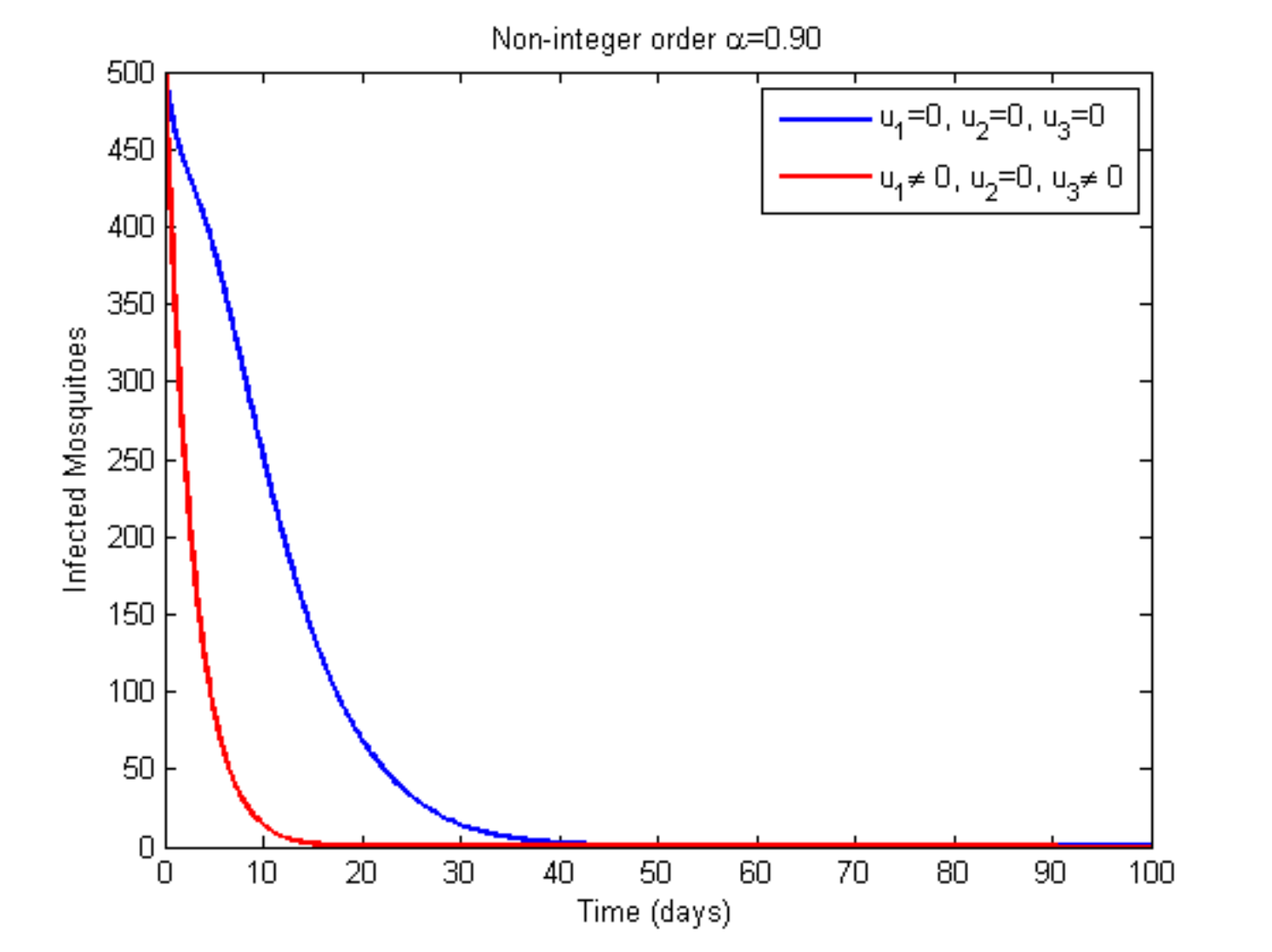}}
\caption{Numerical solutions of infected mosquitoes with treated bednets and insecticide spray controls for $\alpha=1, 0.99, 0.95, 0.90$}
\label{fg12}
\end{figure}

\begin{figure}[!ht]
\centering
\subfigure[]{
\includegraphics[scale=0.5]{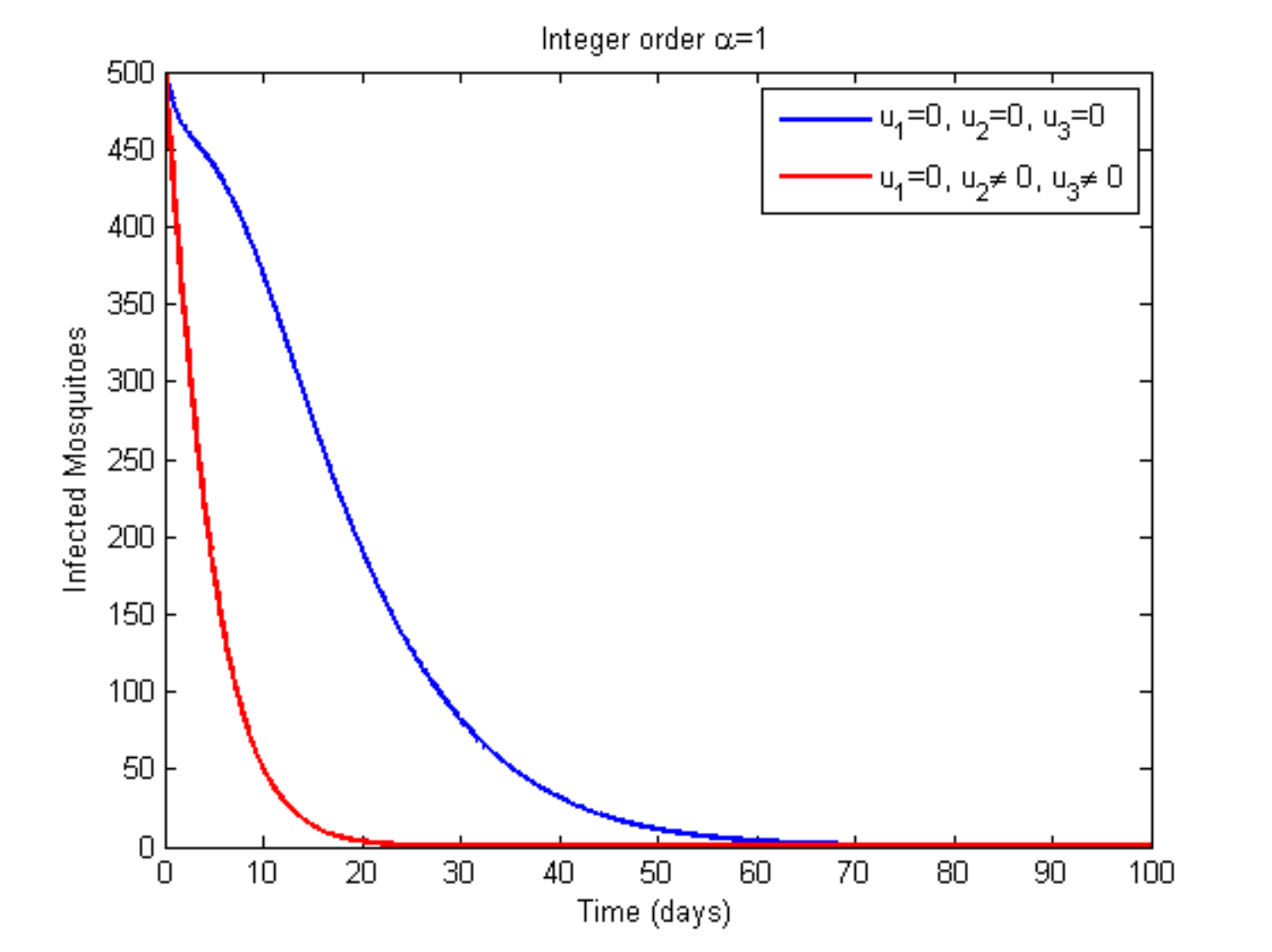}}\hfil
\subfigure[]{
\includegraphics[scale=0.5]{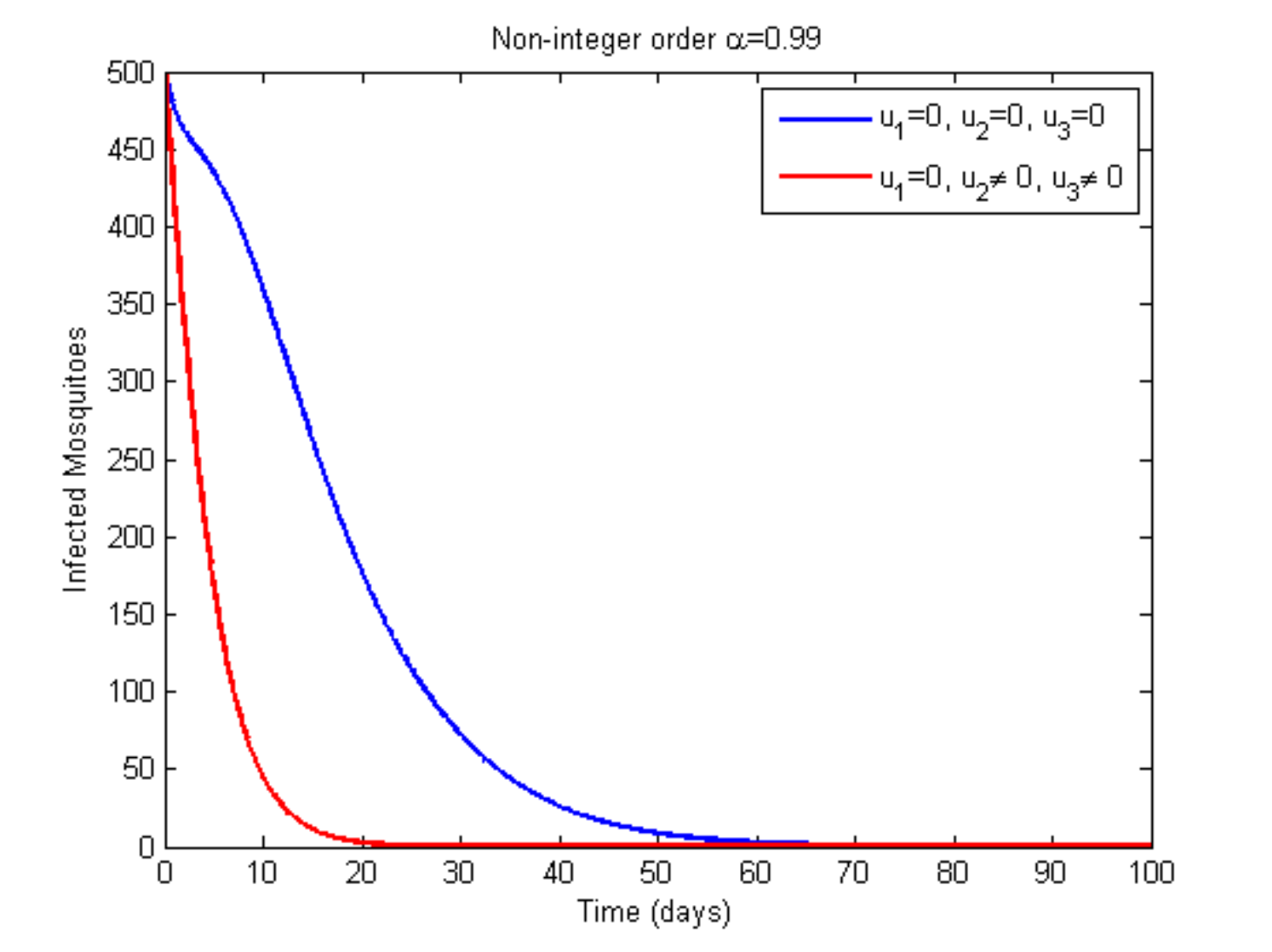}}\hfil
\subfigure[]{
\includegraphics[scale=0.5]{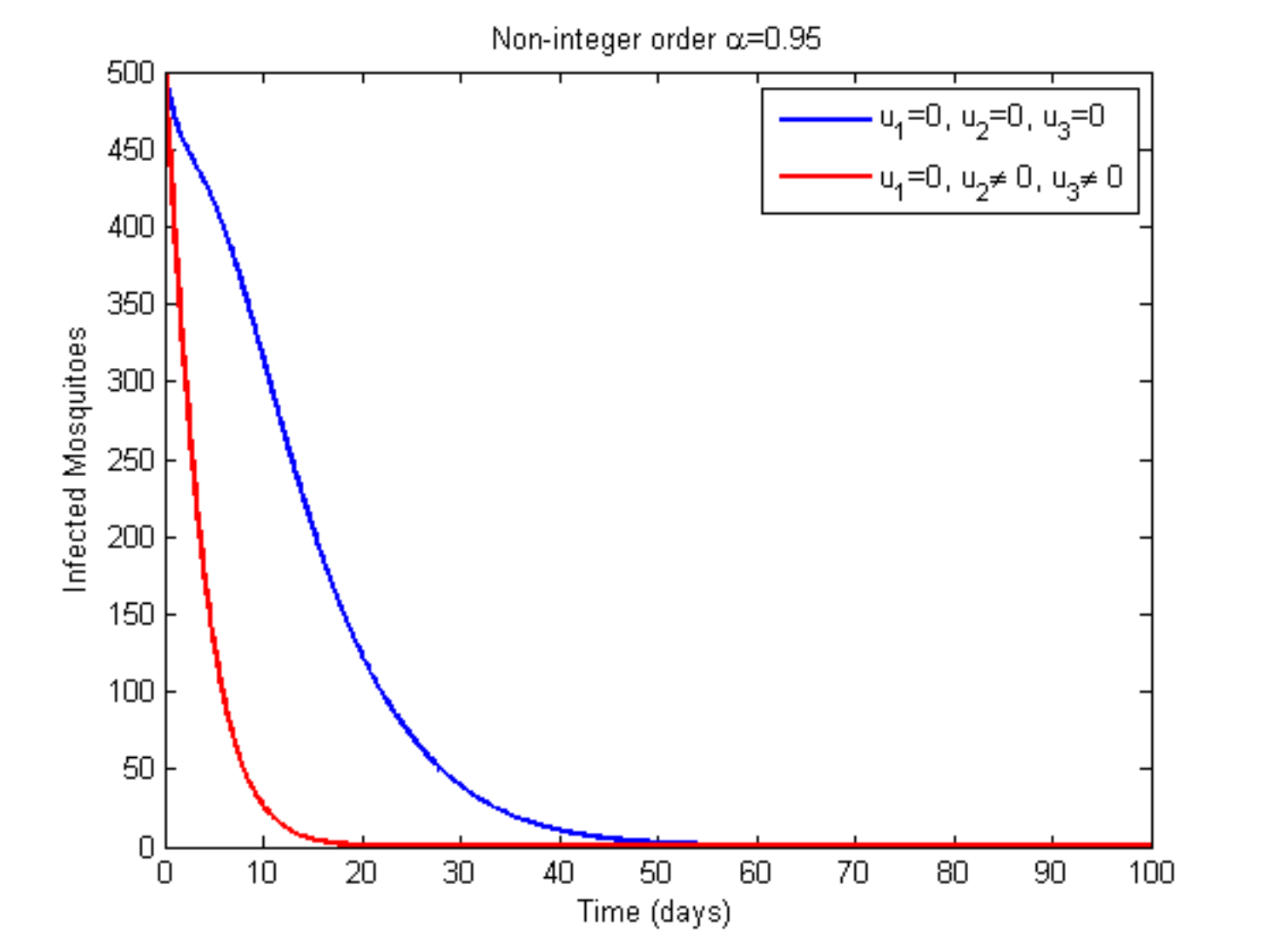}}\hfil
\subfigure[]{\includegraphics[scale=0.5]{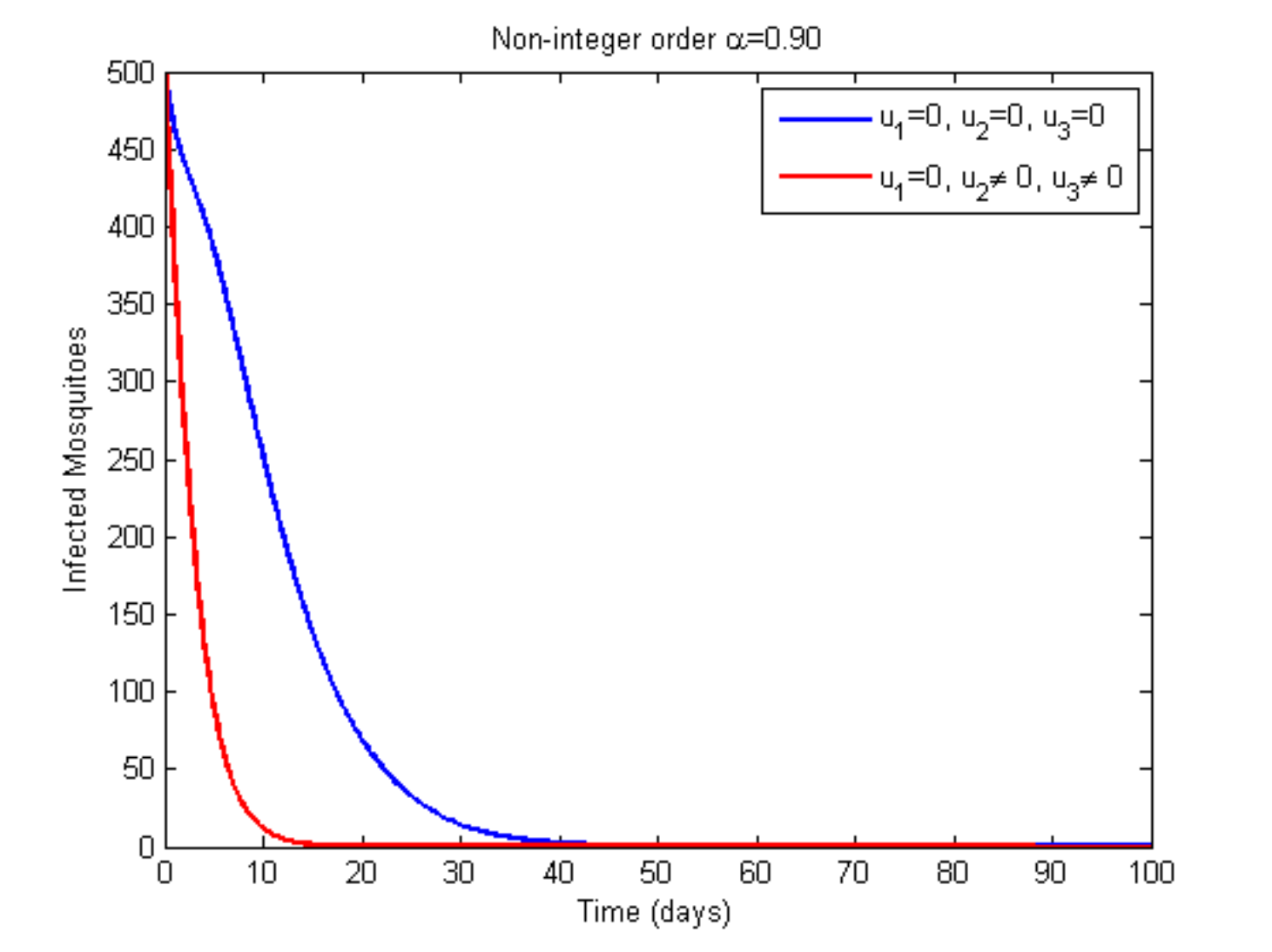}}
\caption{Numerical solutions of infected mosquitoes with treatment and insecticide spray controls for $\alpha=1, 0.99, 0.95, 0.90$}
\label{fg13}
\end{figure}

\begin{figure}[!ht]
\centering
\subfigure[]{
\includegraphics[scale=0.5]{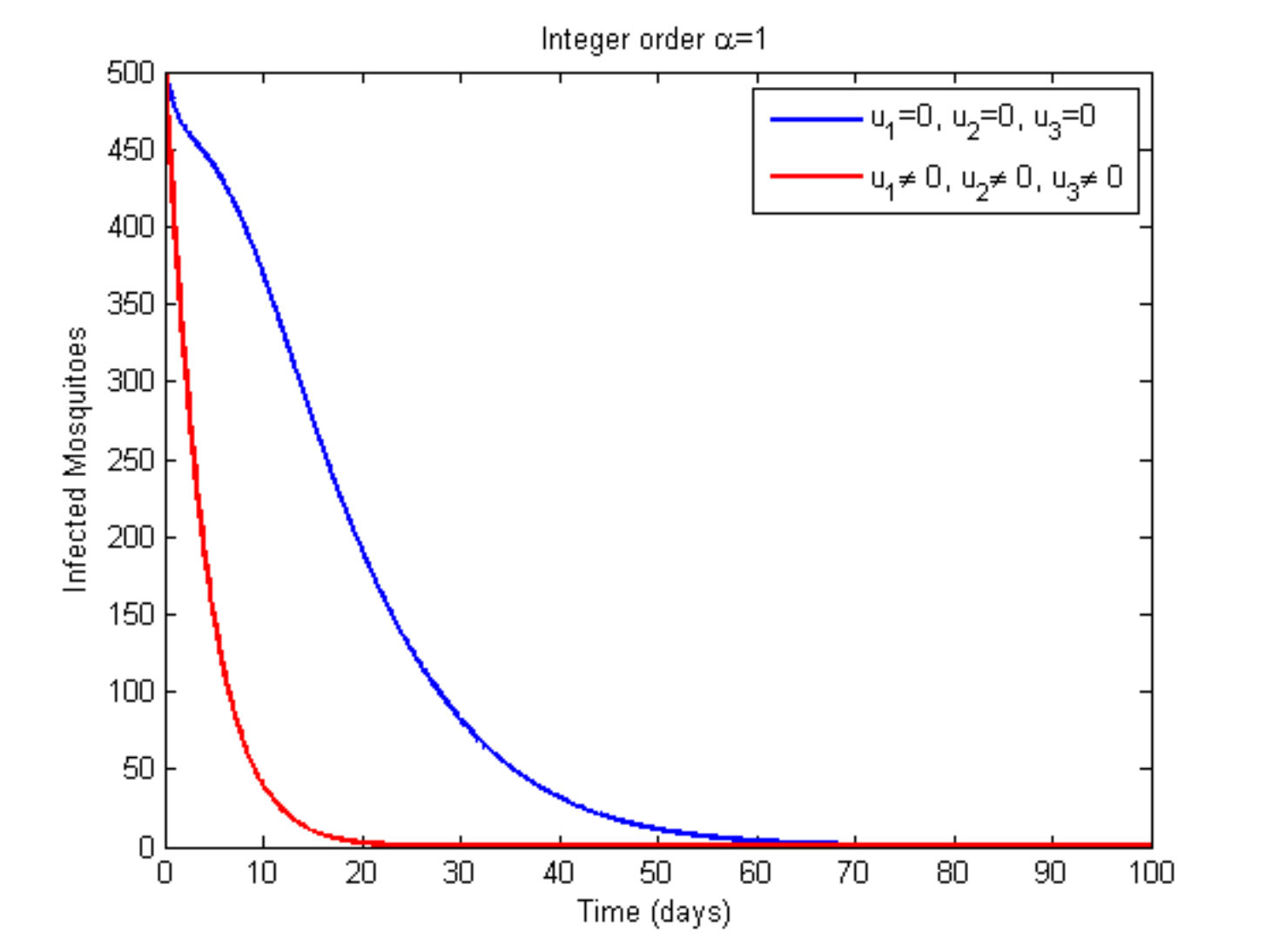}}\hfil
\subfigure[]{
\includegraphics[scale=0.5]{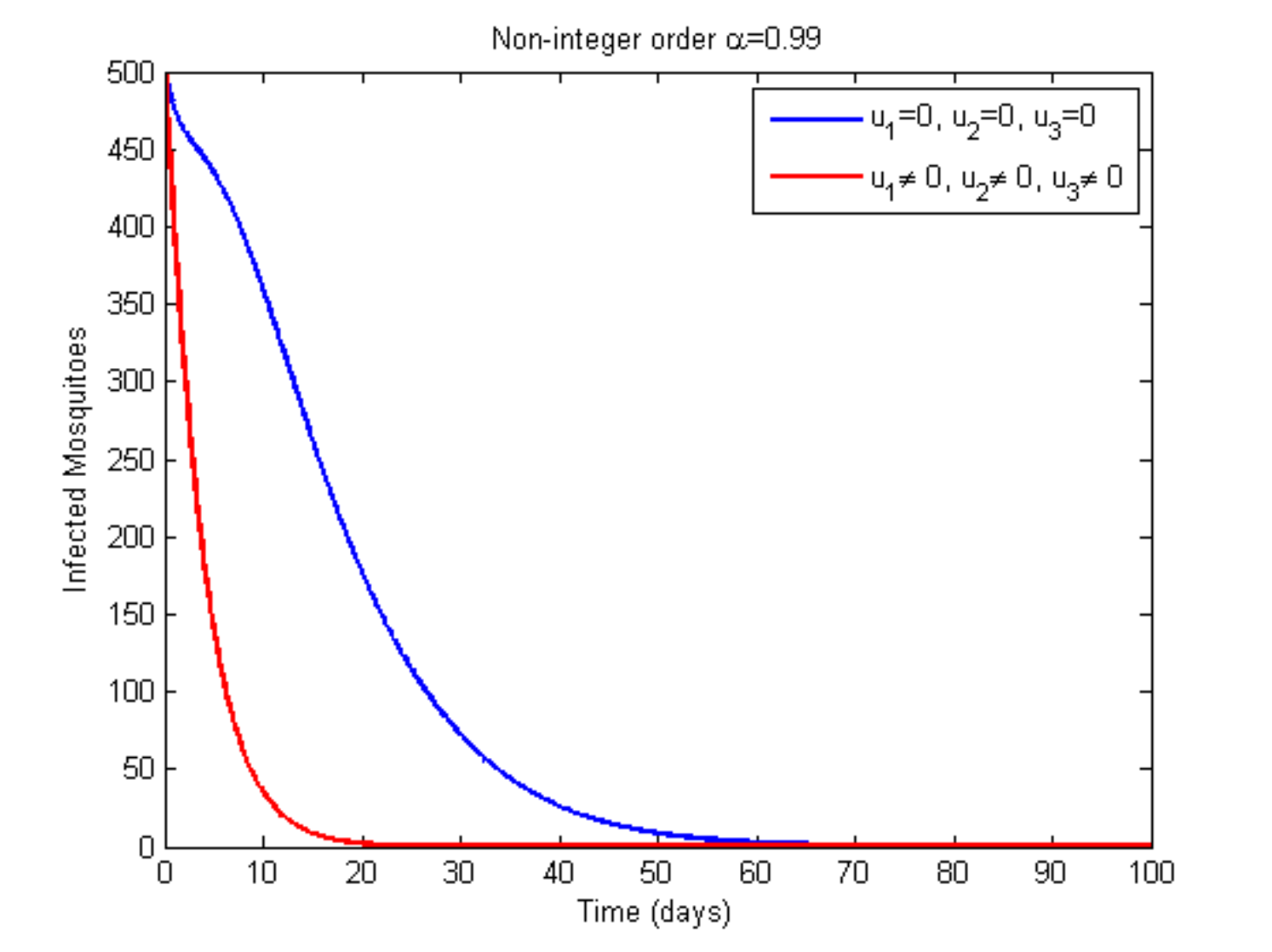}}\hfil
\subfigure[]{
\includegraphics[scale=0.5]{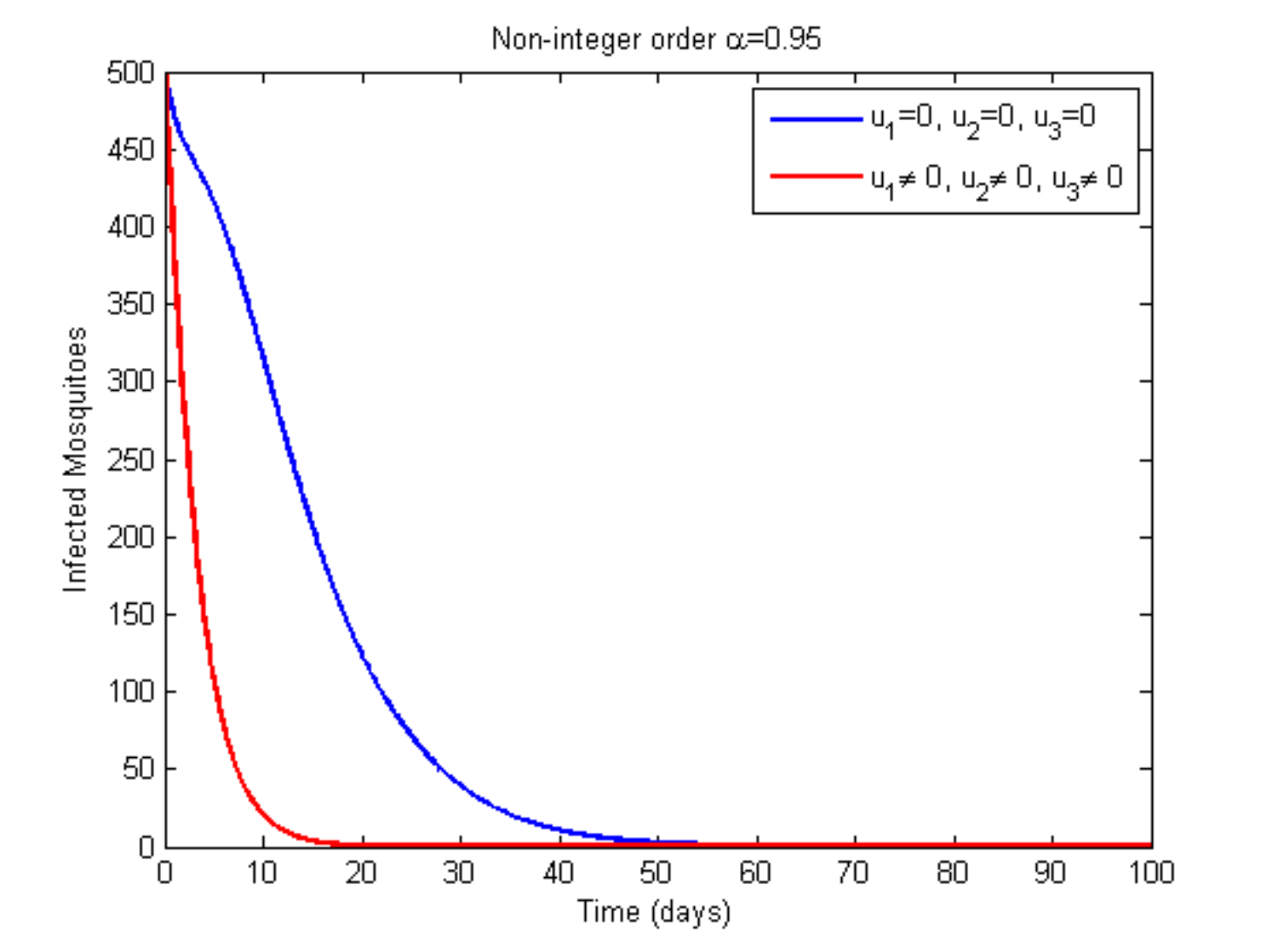}}\hfil
\subfigure[]{\includegraphics[scale=0.5]{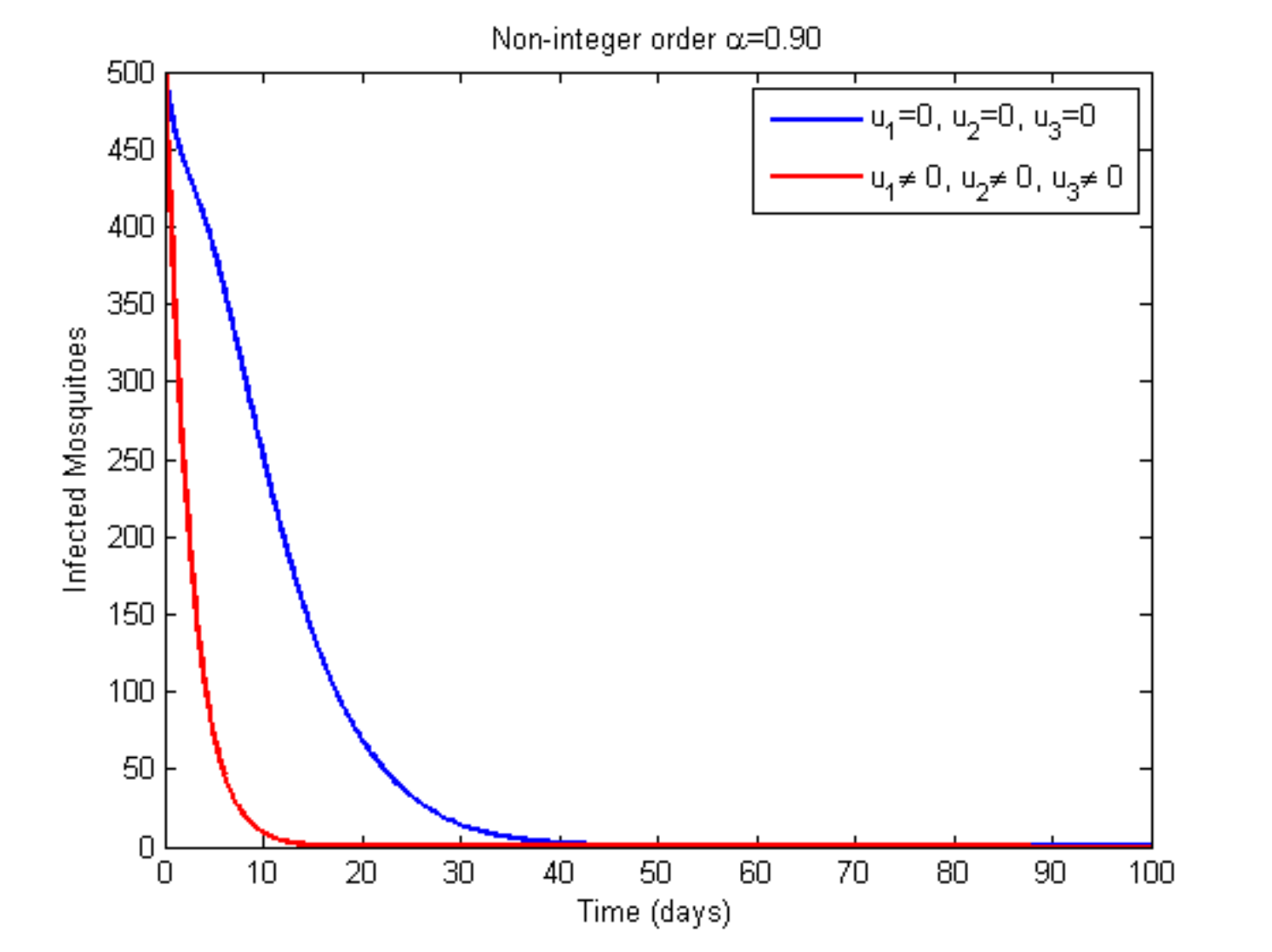}}
\caption{Numerical solutions of infected mosquitoes with all the three time dependent controls (treated bednets, treatment and insecticide spray) for $\alpha=1, 0.99, 0.95, 0.90$}
\label{fg14}
\end{figure}



\begin{figure}[!ht]
\centering
\subfigure[]{
\includegraphics[scale=0.5]{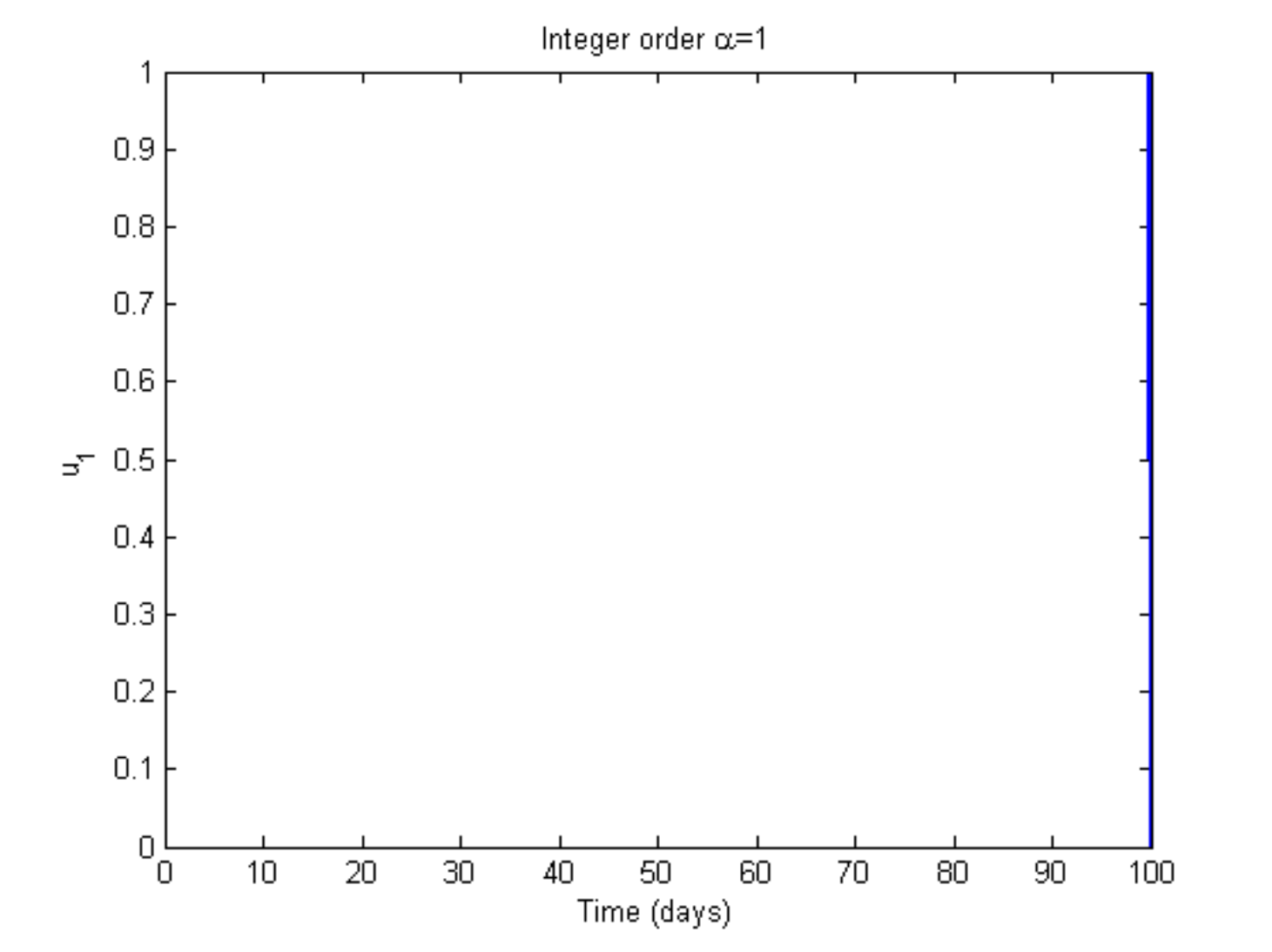}}\hfil
\subfigure[]{
\includegraphics[scale=0.5]{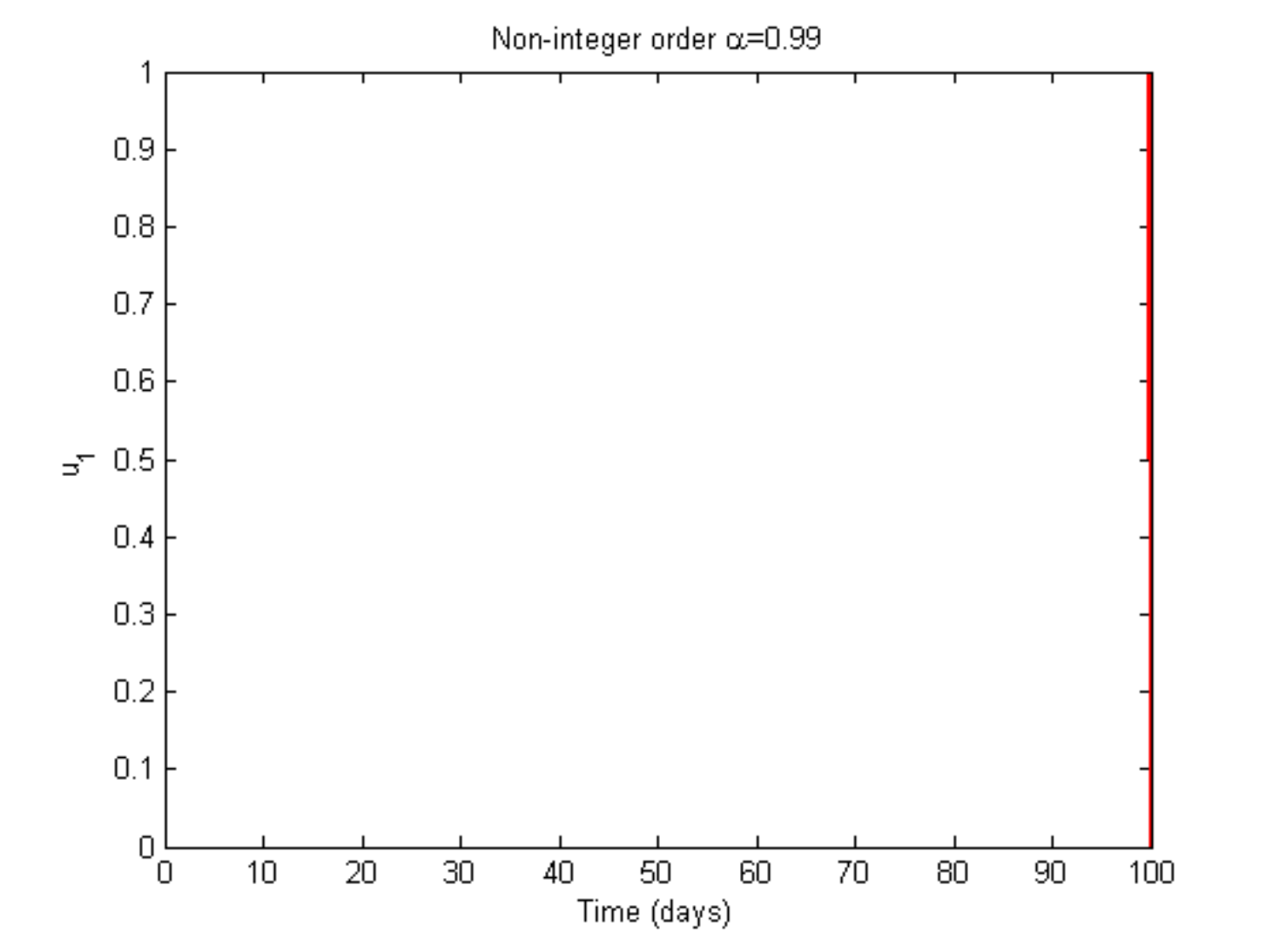}}\hfil
\subfigure[]{
\includegraphics[scale=0.5]{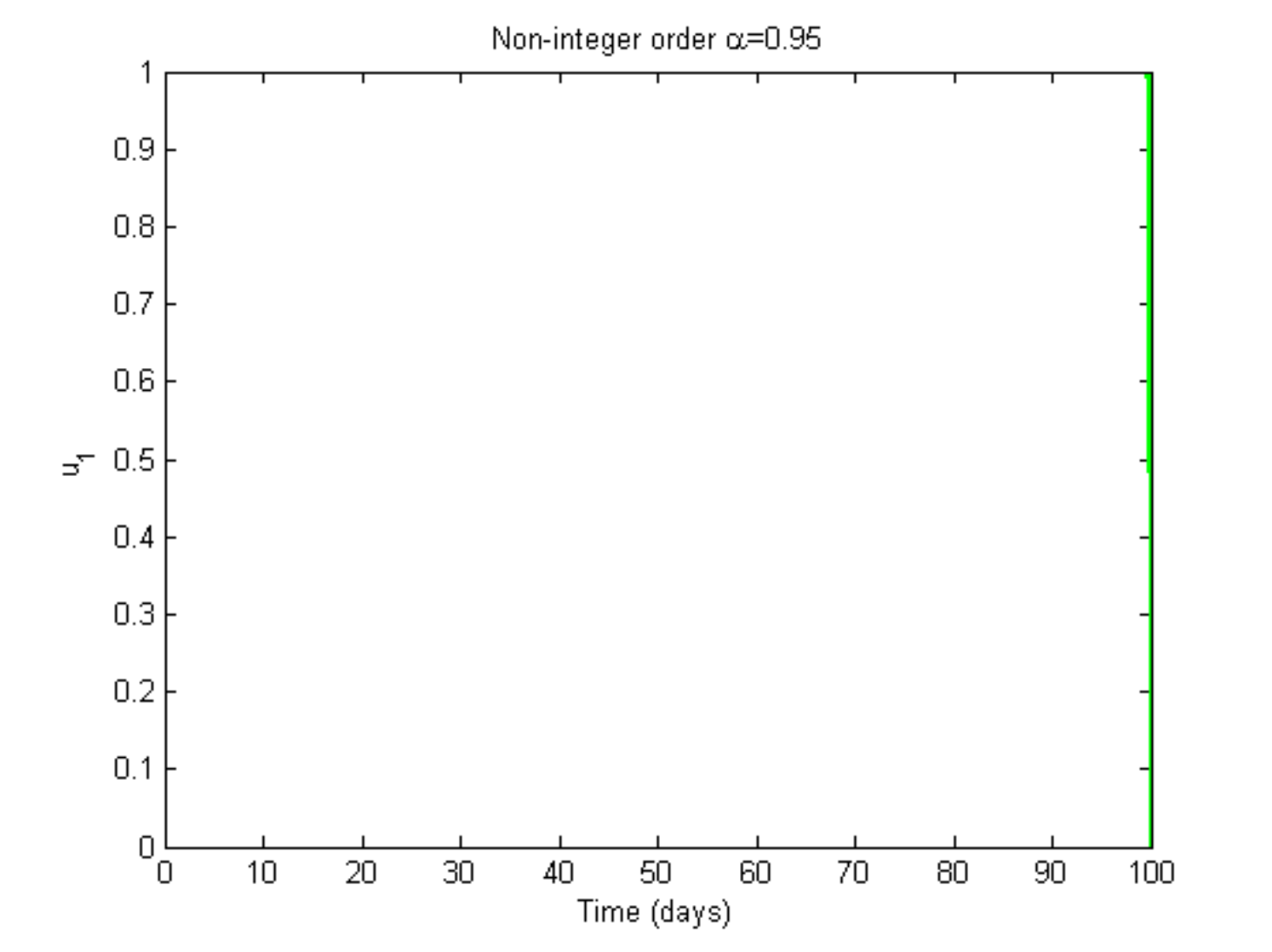}}\hfil
\subfigure[]{\includegraphics[scale=0.5]{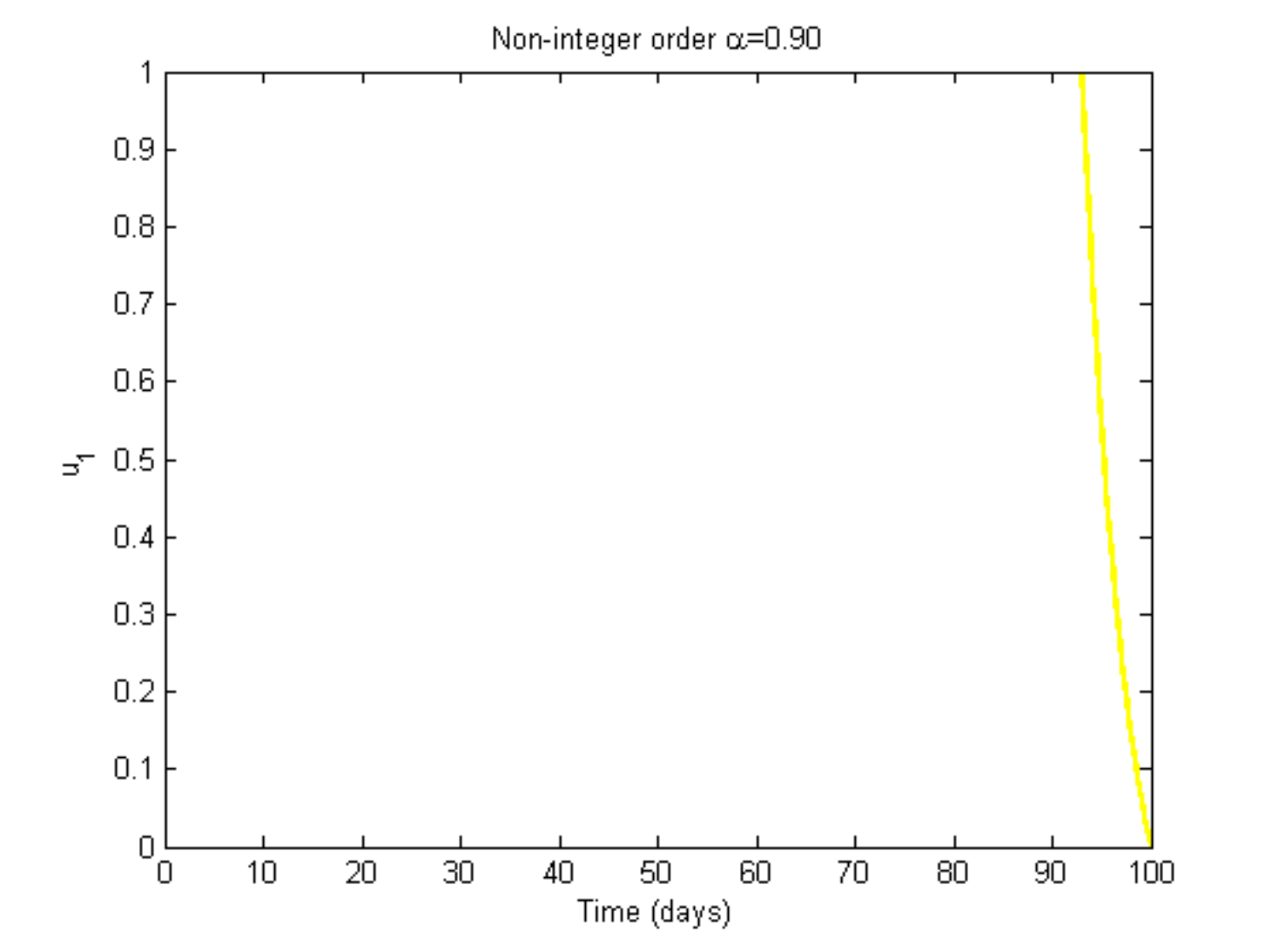}}
\caption{Optimal control function $u_1$ for $\alpha=1, 0.99, 0.95, 0.90$}
\label{fg15}
\end{figure}

\begin{figure}[!ht]
\centering
\subfigure[]{
\includegraphics[scale=0.5]{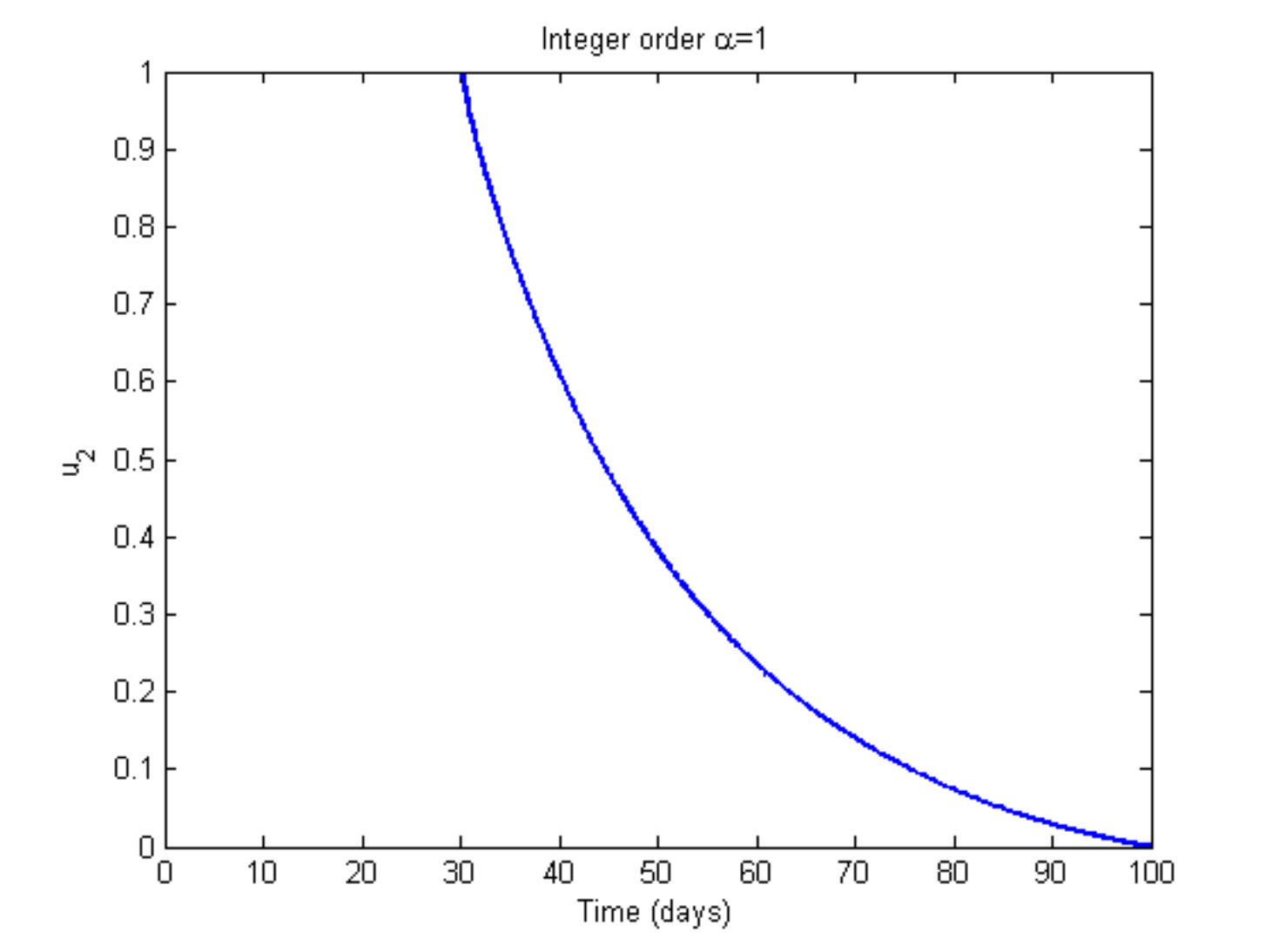}}\hfil
\subfigure[]{
\includegraphics[scale=0.5]{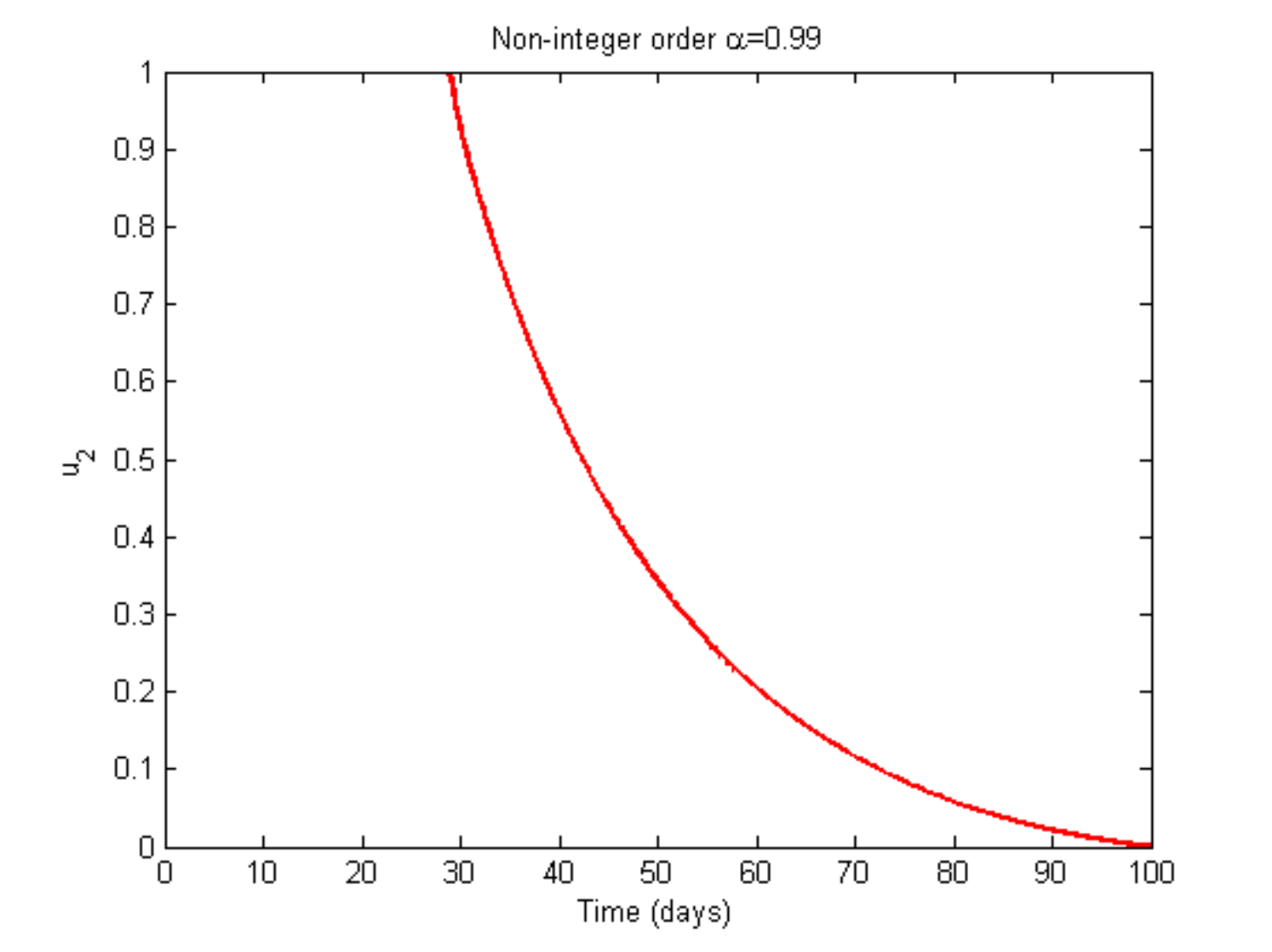}}\hfil
\subfigure[]{
\includegraphics[scale=0.5]{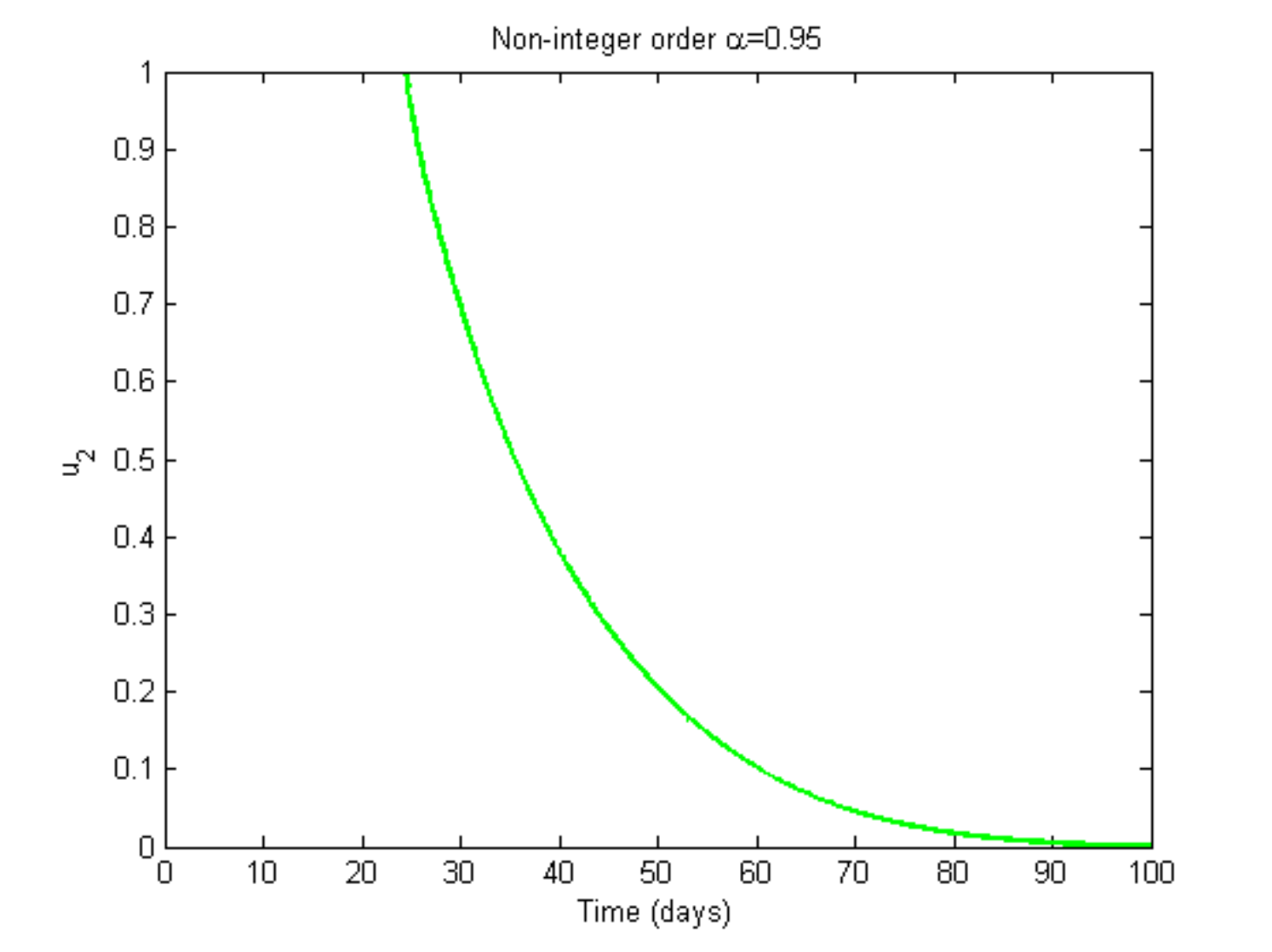}}\hfil
\subfigure[]{\includegraphics[scale=0.5]{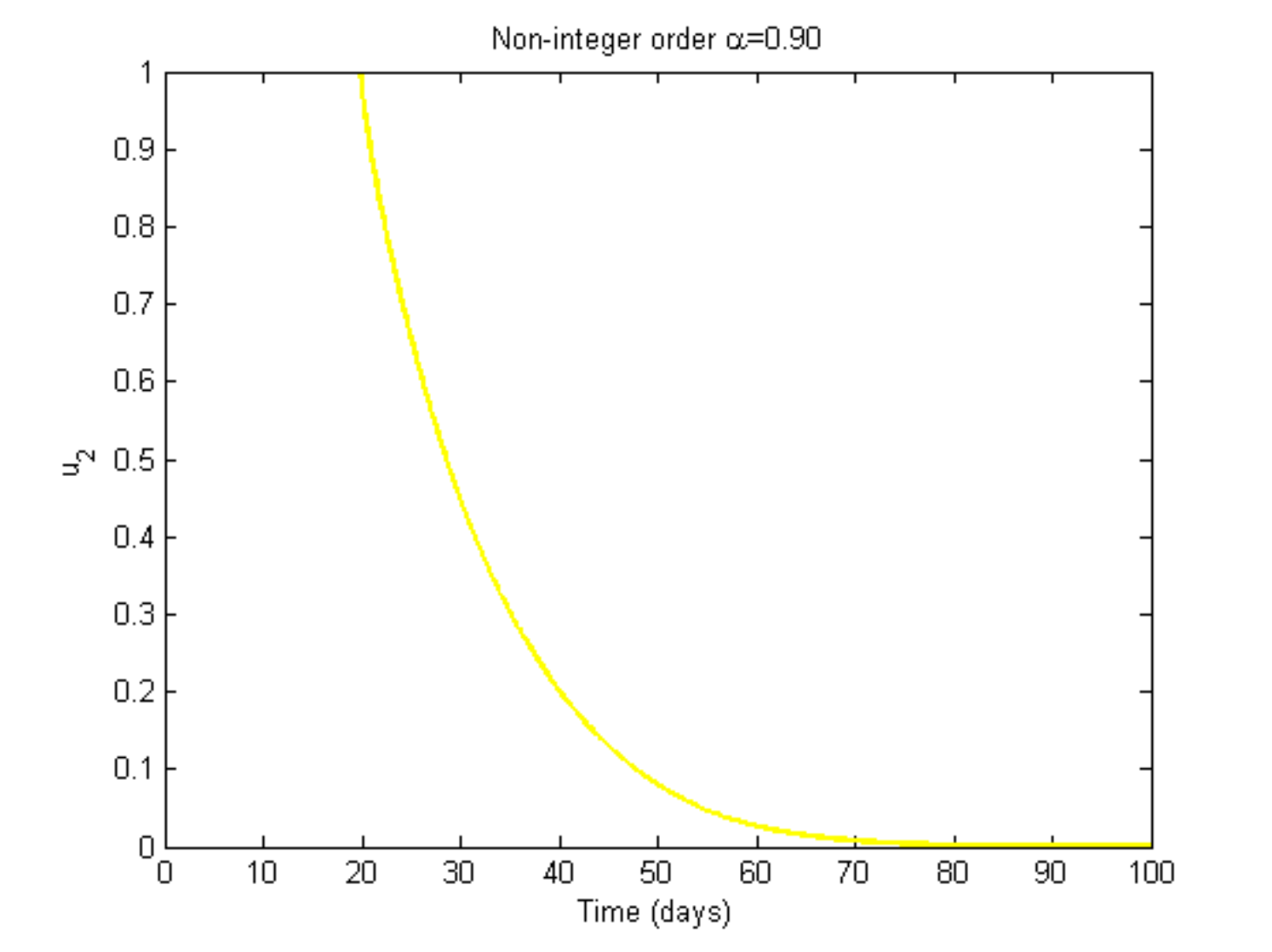}}
\caption{Optimal control function $u_2$ for $\alpha=1, 0.99, 0.95, 0.90$}
\label{fg16}
\end{figure}

\begin{figure}[!ht]
\centering
\subfigure[]{
\includegraphics[scale=0.5]{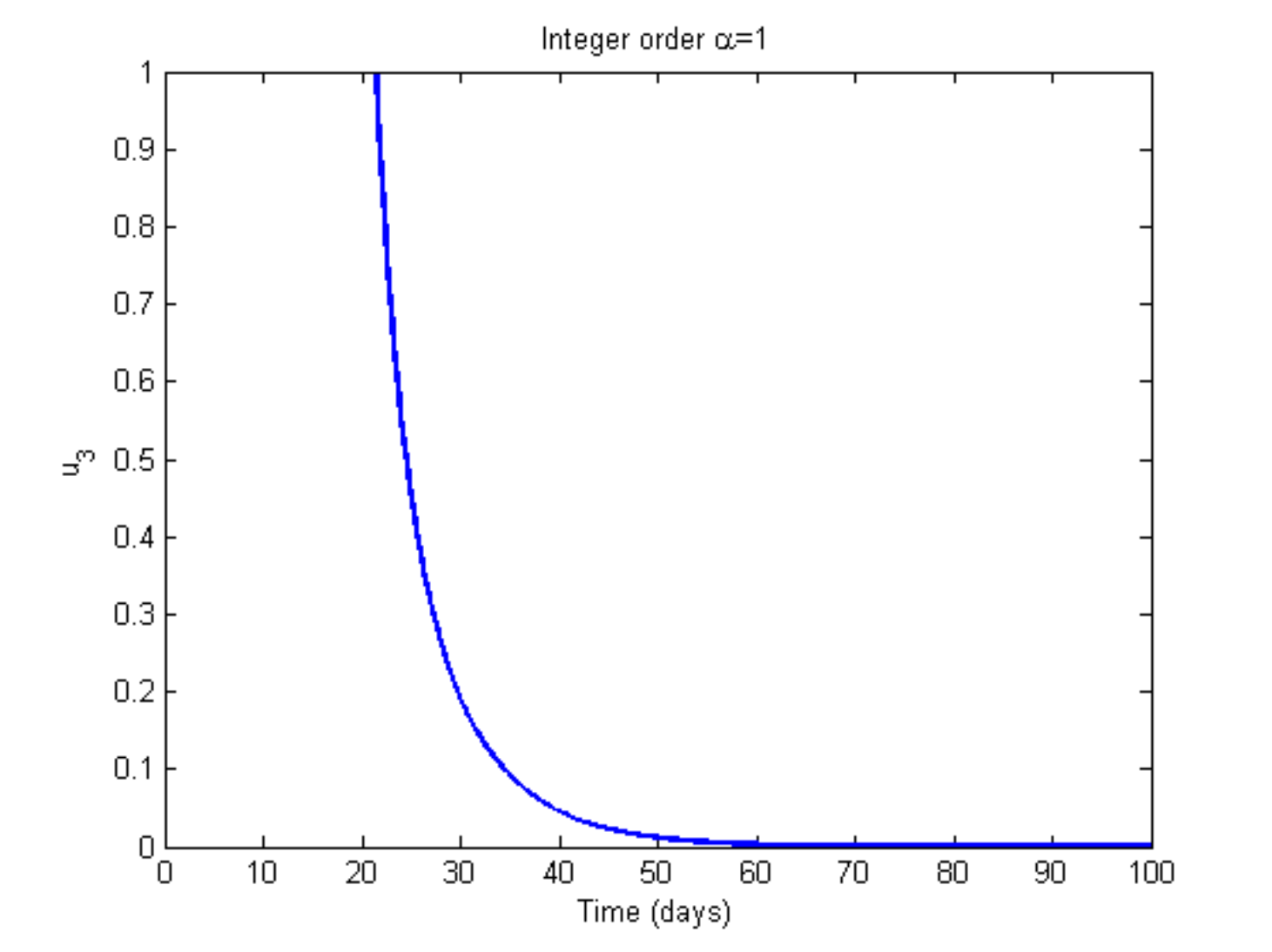}}\hfil
\subfigure[]{
\includegraphics[scale=0.5]{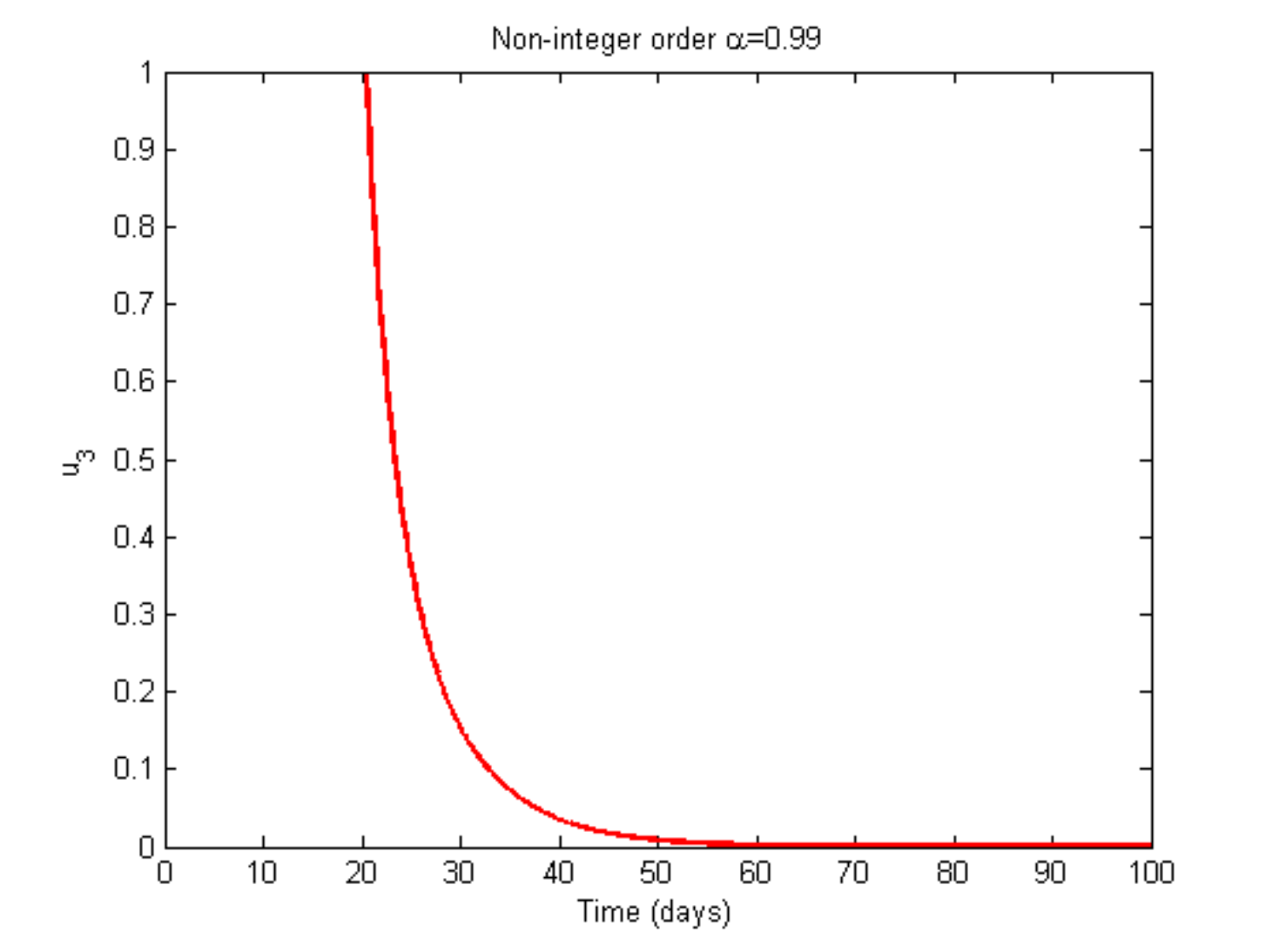}}\hfil
\subfigure[]{
\includegraphics[scale=0.5]{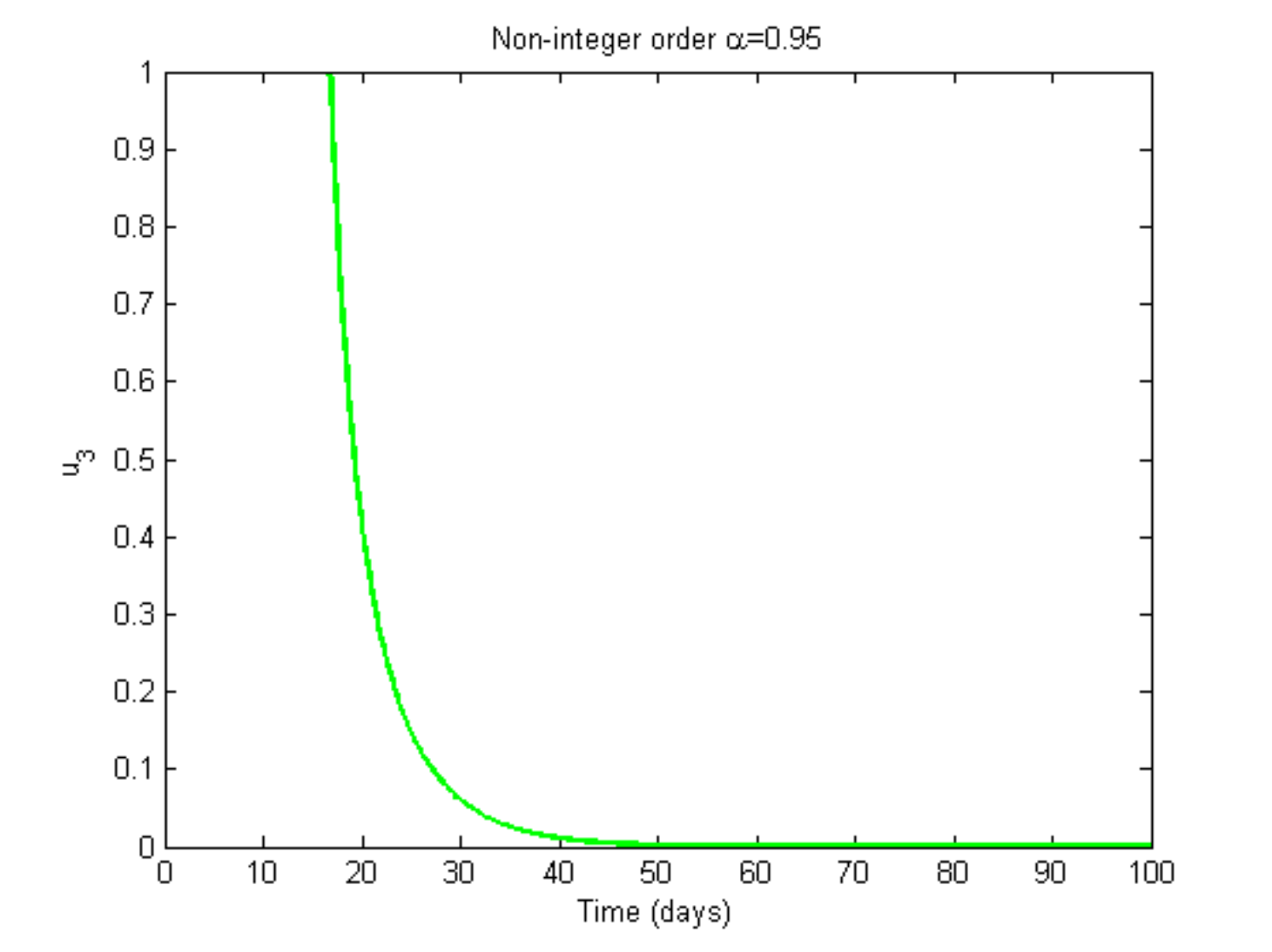}}\hfil
\subfigure[]{\includegraphics[scale=0.5]{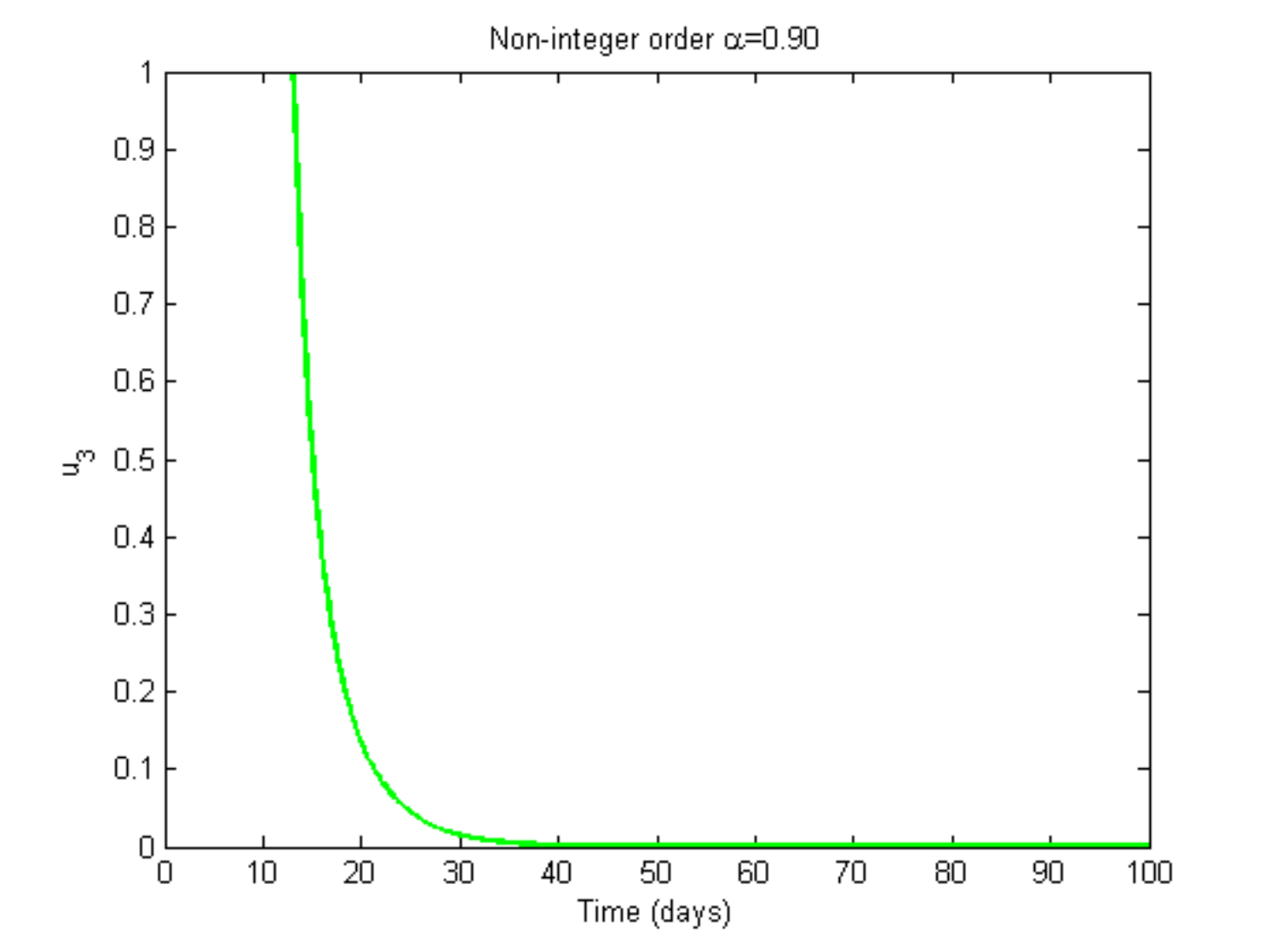}}
\caption{Optimal control function $u_3$ for $\alpha=1, 0.99, 0.95, 0.90$}
\label{fg17}
\end{figure}

\begin{figure}[!ht]
\centering
\subfigure[]{
\includegraphics[scale=0.5]{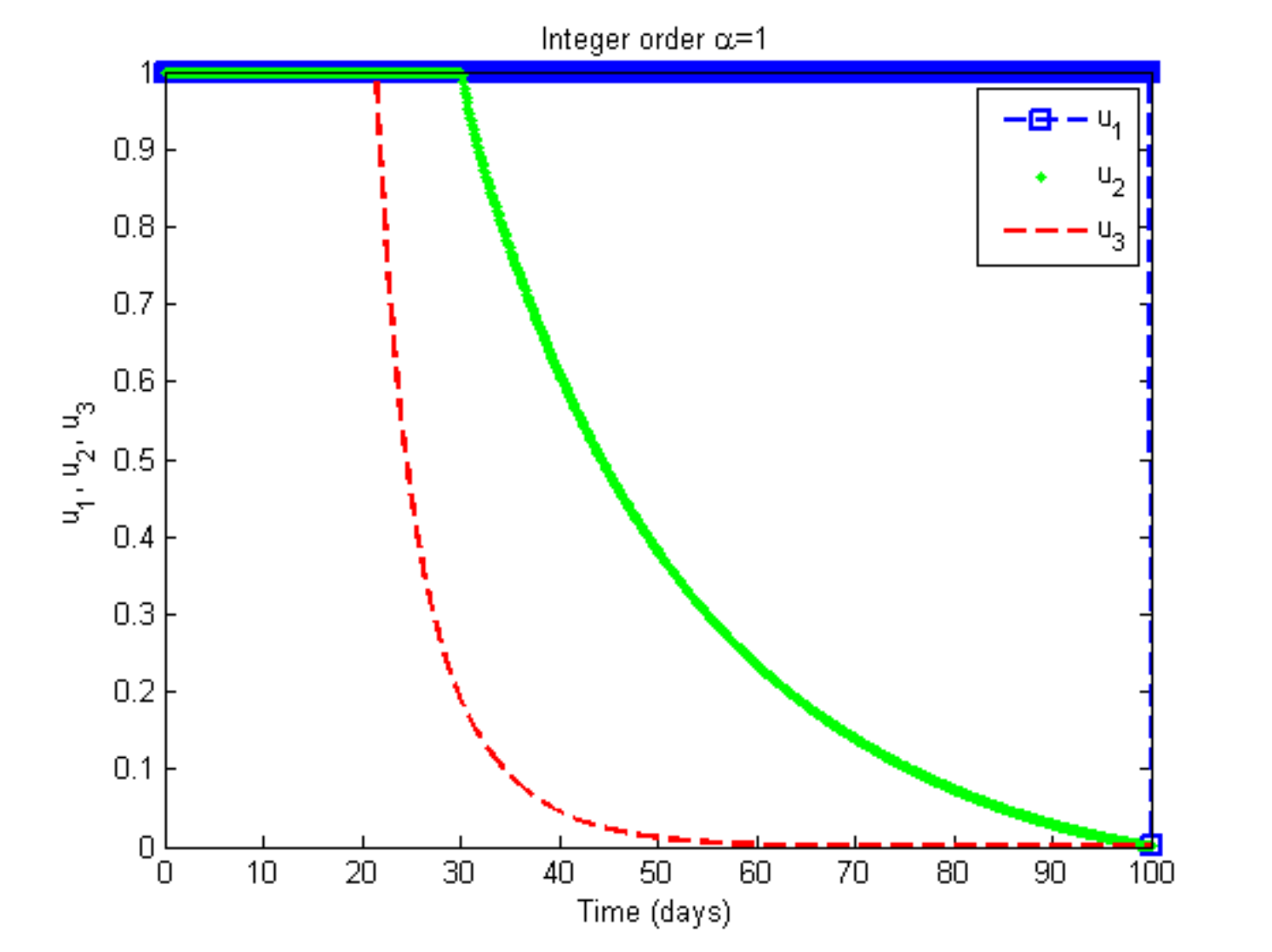}}\hfil
\subfigure[]{
\includegraphics[scale=0.5]{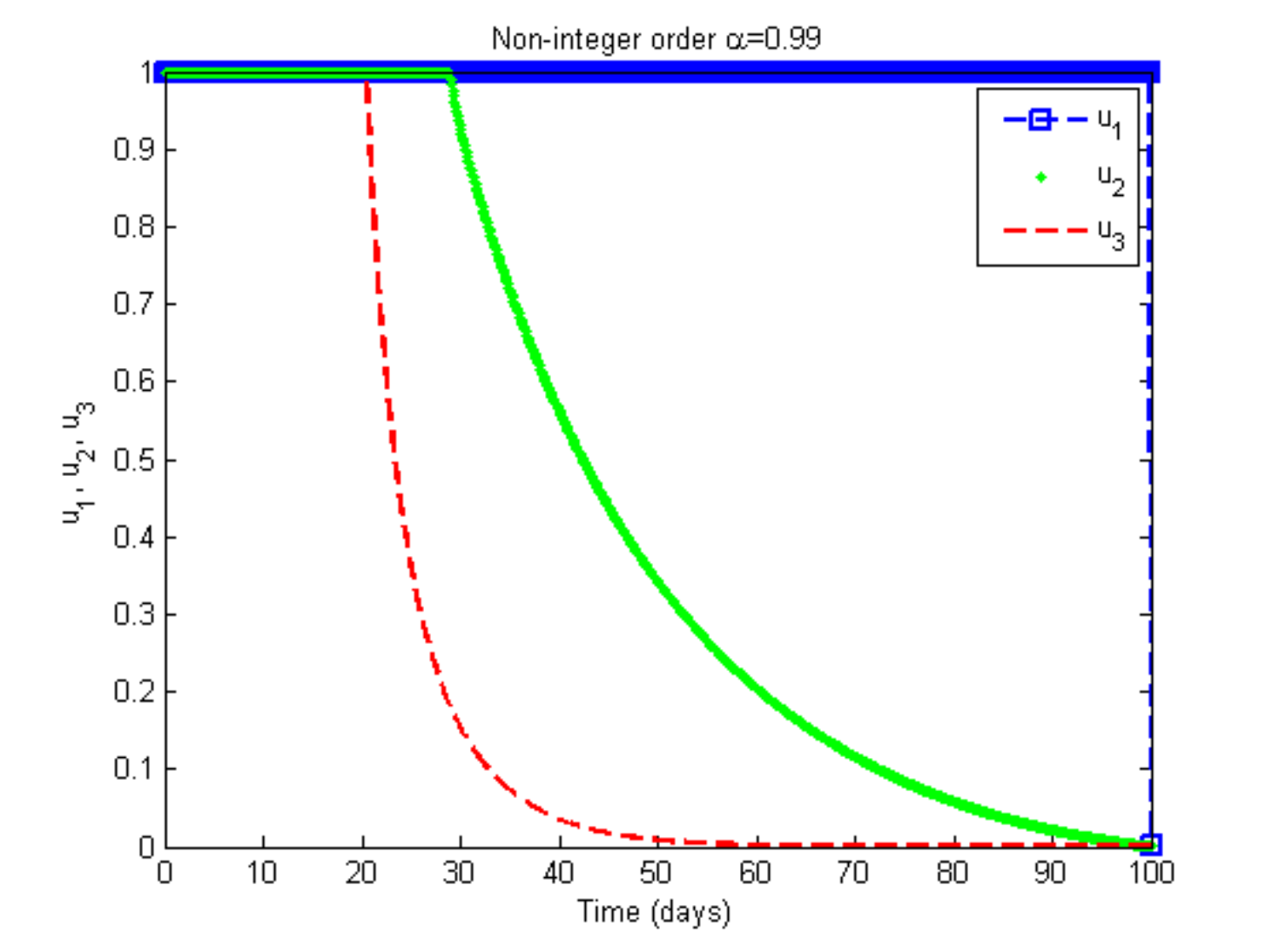}}\hfil
\subfigure[]{
\includegraphics[scale=0.5]{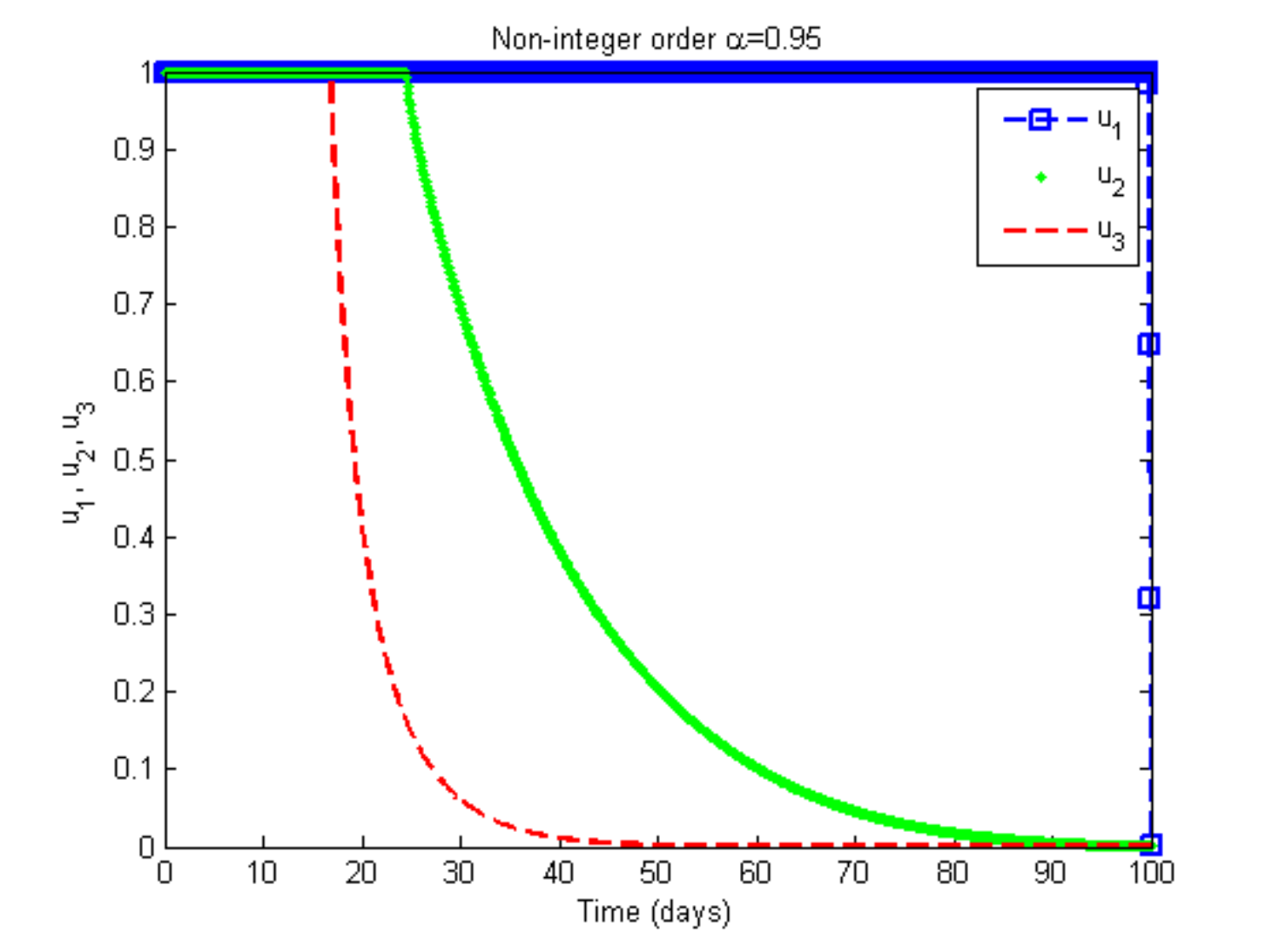}}\hfil
\subfigure[]{\includegraphics[scale=0.5]{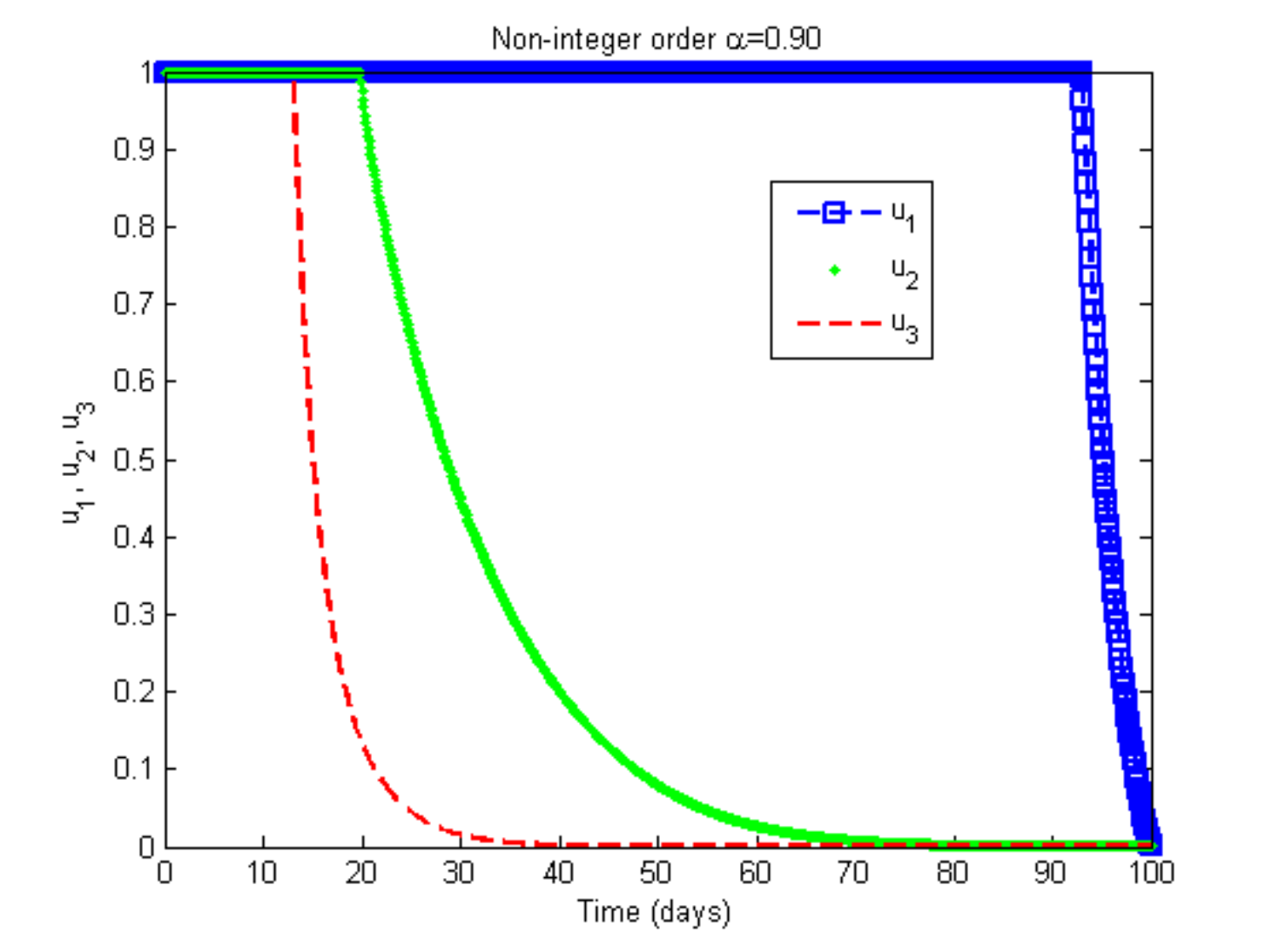}}
\caption{Optimal control functions $u_1, \ u_2 $ and $u_3$ for $\alpha=1, 0.99, 0.95, 0.90$}
\label{fg18}
\end{figure}

\section{Conclusion}\label{threesec}
This work considered an optimal control epidemiological model for malaria infection formulated in the sense of Caputo derivatives. The optimality of the fractional optimal control problem was solved using the forward-backward sweep method and the generalized euler method. We have numerically shown that when any of the time dependent controls or any combination of them are applied, the infected populations for both human and mosquitoes are reduced. Furthermore, when all the time dependent controls are used in the simulations, the reductions in the infected populations are more than any other control and prevention strategies. Our numerical solutions for the optimal control problem with fractional orders are much better than the integer order.
\bibliography{sample2}

\end{document}